\documentclass[a4paper,12pt,reqno]{amsart}
\usepackage{amsmath}
\usepackage{amssymb}
\usepackage{bbm}
\usepackage{graphicx}
\usepackage[all]{xypic}

\usepackage{graphicx}  
\usepackage{color}     

\newlength{\realsidemargin}
\newlength{\sidemargin}
\newcommand{\NewTheorem}[1]{
	\newtheorem{#1}[TheoremEnvironment]{#1}
}

\NewTheorem{Definition}
\NewTheorem{Remark}
\NewTheorem{Axiom}
\NewTheorem{Example}
\NewTheorem{Conventions}
\NewTheorem{Notations}
\NewTheorem{Setting}
\NewTheorem{Observation}
\NewTheorem{Question}
\NewTheorem{Conjecture}
\NewTheorem{Homework}

\NewTheorem{Proposition}
\NewTheorem{Lemma}
\NewTheorem{Theorem}
\NewTheorem{Corollary}

\newcommand{\mbZ}{\mathbb{Z}}

\newcommand{\mcA}{\mathcal{A}}
\newcommand{\mcB}{\mathcal{B}}
\newcommand{\mcC}{\mathcal{C}}
\newcommand{\mcD}{\mathcal{D}}
\newcommand{\mcE}{\mathcal{E}}
\newcommand{\mcF}{\mathcal{F}}

\newcommand{\mcH}{\mathcal{H}}
\newcommand{\mcI}{\mathcal{I}}
\newcommand{\mcJ}{\mathcal{J}}
\newcommand{\mcK}{\mathcal{K}}

\newcommand{\mcM}{\mathcal{M}}

\newcommand{\mcO}{\mathcal{O}}
\newcommand{\mcP}{\mathcal{P}}

\newcommand{\mcR}{\mathcal{R}}
\newcommand{\mcS}{\mathcal{S}}
\newcommand{\mcT}{\mathcal{T}}
\newcommand{\mcU}{\mathcal{U}}

\newcommand{\set}[2][]{\mathopen{#1\{}#2\mathclose{#1\}}}
\newcommand{\conditionalset}[3][]{\mathopen{#1\{}\,#2\mathrel{#1|}#3\,\mathclose{#1\}}}
\newcommand{\generatedset}[2][]{\mathopen{#1\langle}#2\mathclose{#1\rangle}}

\DeclareMathOperator{\Hom}{Hom}
\DeclareMathOperator{\End}{End}

\DeclareMathOperator{\Ob}{Ob}

\newcommand{\into}{\hookrightarrow}

\DeclareMathOperator{\QCoh}{QCoh}

\DeclareMathOperator{\Coh}{Coh}

\DeclareMathOperator{\Add}{Add}
\DeclareMathOperator{\add}{add}
\DeclareMathOperator{\smd}{smd}
\DeclareMathOperator{\Ker}{Ker}

\DeclareMathOperator{\Spec}{Spec}
\DeclareMathOperator{\Spc}{Spc}

\DeclareMathOperator{\Supp}{Supp}
\DeclareMathOperator{\supp}{supp}
\DeclareMathOperator{\Supph}{Supph}
\DeclareMathOperator{\supph}{supph}

\DeclareMathOperator{\DD}{\mathbf{D}}

\newcommand{\Dqc}{\DD_{\mathrm{qc}}}
\newcommand{\Dcoh}{\DD^{b}_{\mathrm{coh}}}
\newcommand{\Dperf}{\DD^{\mathrm{perf}}}
\newcommand{\Ld}{{\mathbf{L}}}
\newcommand{\Rd}{{\mathbf{R}}}
\newcommand{\one}{\mathbbm{1}}

\newcommand{\Ho}{
\mcS \mcH
}

\begin{document}

\title
[From Ohkawa to strong generation via
approximable triangulated categories]
{
From Ohkawa to strong generation via
\\
approximable triangulated categories 
\\
–-
a variation on the theme of 
\\
Amnon Neeman's Nagoya lecture series 
}

%
\author{Norihiko Minami 
}
\address{
Nagoya Institute of Technology  \\
Gokiso, Showa-ku, Nagoya 466-8555, JAPAN
}

 \email{nori@nitech.ac.jp} 
%
%
\thanks{This work was partially supported by JSPS KAKENHI Grant Number 15K04872. }






\maketitle

\begin{abstract}
This survey stems from 
Amnon Neeman's lecture series at
Ohakawa's memorial workshop.
Starting with Ohakawa's theorem,  
this survey intends to supply enough 
motivation, background and technical 
details  to read Neeman's recent 
papers on his \lq\lq approximable 
triangulated categories \rq\rq \  
and his $\Dcoh (X)$ strong generation 
sufficient criterion  
 via de Jong's regular 
alteration, even for non-experts.
%
%
%
%
%
%
\end{abstract}


\section{Introduction
}

This survey stems from Amnon Neeman's
lecture series at Ohakawa's memorial workshop. 
\footnote{
The author also would like 
to thank  Professors
Mitsunori Imaoka, Takao Matumoto, Takeo Ohsawa, Katsumi Shimomura, and 
Masayuki Yamasaki, for coorganizing the 
workshop.
}
%
The original lecture series started and
 ended with Ohkawa's theorem on the stable
 homotopy category. In the beginning 
 Ohkawa's theorem was presented in its 
 lovely, original
form. The lecture series then meandered 
 through some--definitely not all--of the 
 developments and generalizations made by 
 others in the years following Ohkawa's 
 paper. And at the end came what was then 
a recent result of Amnon Neeman's--and the
 relevance was that the Ohkawa set and its 
 properties, as developed in the years 
 following Ohkawa, turned out to be key to 
 the proof of the recent theorem.

Here, our presentation  significantly 
modifies 
Neeman's original 
presentation, partially fueled 
by  other distinguished submissions to this proceedings, 
mostly 
to motivate topologists
to get interested in this rich subject.
%
For this purpose, we have
reorganized and expanded the original 
framework of Amnon Neeman's
lecture series. 

Still, the underlying philosophy of 
Neeman's presentation to start with 
Ohkawas's theorem remains 
kept in this survey.  
And most significantly,
following a strong request of Professor 
Neeman, we reviewed 
Neeman's recent proof of:
\begin{quote}
\hspace{-9mm}
\underline{\em  
$\Dcoh (X)$ strong generation 
sufficient criterion via de Jong's 
regular alteration}\ 
\end{quote}
with enough background 
and technical details, 
expanding and sometimes even modifying 
parts of the original proof so as to 
make this review 
beginner-friendly from a homotopy 
theorist's point of view.
%
Actually, this proof of Neeman also 
makes critical use of, 
in addition to de Jong's regular alteration, a couple of Thomason's theorems:  
\begin{itemize}
\item
First, 
the fundamental theorem of Hopkins, Neeman, 
Thomason and others on the classification of 
thick tensor ideals of $\Dperf (X),$ the 
$\Dperf (X) =\Dqc(X)^c$ 
analogue of the Hopkins-Smith thick 
 subcategory theorem of $\Ho^{fin}=\Ho^c$ 
whose proof heavily depends upon 
the (Devinatz-)Hopkins-Smith nilpotency
theorem.  
\item
Second, 
Thomason's localization
theorem on 
$\Dperf (X \setminus Z),$ for which Neeman 
found a homotopy theoretical proof 
in the framework of Miller's finite 
 localiation. 
\end{itemize}
Considering these circumstance, 
we have also explained  
the role of (Devinatz-)Hopkins-Smith 
 nilpotency theorem in the proof of 
Hopkins-Smith thick subcategory theorem, 
as well as essentially all the details of 
Neeman's proof of 
Thomason's localization 
theorem.

Now the rest of this survey is 
organized as follows:

\begin{description}
\item[\underline{\S 2}:]
The first goal of this section is to recall
Ohkawa's theorem in 
stable homotopy 
theory.  Ohkawa's theorem claims the 
Bounsfield classes in the stable homotopy
category $\mcS \mcH$ form 
a set which is
very mysterious and beyond our imagination.
Then the second goal of this section is 
the fundamental theorem of 
Hopkins, Neeman, Thomason and others, which 
roughly states the analogue of the 
Bousfield classes in $\Dqc (X),$ in contrast to 
the Ohakawa's case of $\mcS \mcH,$ 
form a set with 
a clear algebro-geometric description.
For these purposes, standard facts about 
the Bousfield localization and triangulated 
categories are reviewed, including 
the existence of Bousfield localization 
for perfectly generated 
triangulated subcategories, 
Miller's finite localization for 
triangulated subcategories generated 
a set of compact objects, and 
the telescope conjecture.,

\item[\underline{\S 3}:]
In reality, Hopkins was not motivated by 
Ohkawa's Theorem~\ref{Ohkawa} for his 
influential paper in algebraic geometry  \cite{MR0932260} 
(Theorem~\ref{Ohkawa for D_{qc}(X)}).
Instead, Hopkins was motivated by 
his own theorem with Smith 
\cite{MR1652975} 
in the triangulated subcategory 
 $\mcS \mcH^c$ consisting of compact 
 objects, whose validity was already known 
 to them back around the time Hopkins wrote
\cite{MR0932260}.
In this section, we review this theorem 
of Hopkins-Smith, emphasizing the way how 
(Devinatz-)Hopkins-Smith nilpotency theorem
is used in its proof,  
In Theorem~\ref{The diagram for SH}, we 
summarize the main stories in 
$\Ho_{(p)}^c \subset \Ho_{(p)}$ 
(the Ohkawa theorem, the Hopkins-Smith theorem, 
Miller's version of the Ravenel telescope conjecture \ 
($C\circ I \overset{?}= Id_{\mathbb{T}\left( \Ho_{(p)}^{fin} \right)}$),
 and the conjectures of Hovey 
and Hovey-Palmieri) in the following succint 
commutative diagram: 

%
%
\begin{equation}
\label{THE DIAGRAM-SH-introduction}
\xymatrix{
{
\text{{\it mysterious} set}
}
\ar@{=}
[rr]^-{\text{Ohkawa Th.}
} & &
\mathbb{B}( \Ho_{(p)} )
\ar@{^{(}-_{>}}[rr]^{\text{Hovey Conj.}}
_{\overset{?}=}
 & &
\mathbb{L}( \Ho_{(p)} )
\\
\underset
{\cdots \subsetneq C_{n+1}
\cdots \subsetneq C_{n} \cdots}
{
{
\text{chromatic hierachy}
}
}
\ar@{^{(}-_{>}}[u]
\ar
@{=}[rr]^-{
\text{{\tiny Hopkins-Smith Th.}
}
} & &
\mathbb{T}\left( \Ho_{(p)}^{fin} \right)
\ar@<1mm>[rr]^-{I\ \text{(split inj.)}
} &  &
\ar@<1mm> [ll]^-{C\ \text{(split surj.)}
%
}
\mathbb{S}( \Ho_{(p)} )
\ar@{^{(}-_{>}}[u]
}
\end{equation}
%
%

We then review anlogues of the Hopkins-Smith 
theorem in the motivic setting 
by Joachimi and Kelly.
Also, inspired by this influence of 
Hopkins-Smith 
theorem to algebra and algebraic geometry, 
we briefly reviewed the 
couple of most prominent conjectures in 
homotopy theory, the telescope 
conjecture and the chromatic splitting 
conjecture, following a suggestion of 
Professor Morava.

\item[\underline{\S 4}:] From the previous
two sections, we are naturally led to 
investigate 
$\Dqc (X)^c.$  However, the story  
is not so simple.
Whereas there is a conceptually simple 
algebro-geometrical interpretation 
$\Dqc (X)^c = \Dperf (X),$ 
it is 
its close relative (actually equivalent if 
$X$ is smooth over a field) 
$\Dcoh (X)$ which traditionally has been 
 intensively studied
because of its rich geometric and physical 
information. 
So, we wish to understand both $\Dcoh (X)$ 
and $\Dperf (X).$
In this section, 
we start with brief,
 and so inevitably 
incomplete, 
summaries of  
$\Dcoh (X)$ and $\Dperf (X),$ 
focusing on their usages.
Still, we hope this would convince 
non-experts that $\Dcoh (X)$ and $\Dperf (X)$ 
are very important objects to study.
Amongst of all, we shall recall  
the fundamental theorem of Hopkins, Neeman, 
Thomason and others on the classification of 
thick tensor ideals of $\Dperf (X)$ and 
the  
Thomason's localization theorem 
on 
$\Dperf (X \setminus Z),$ 
both of which play critical roles in 
Neeman's proof of the strong generation 
of $\Dcoh (X)$ reviewed in \S 5. 
For the classification of 
thick tensor ideals of $\Dperf (X),$
we shall establish the 
following commmutative diagram 
\eqref{THE DIAGRAM} in 
Theorem~\ref{MAIN THEOREM}, which is the 
$\Dqc^c(X) = \Dperf(X)$ analogue of the 
Hopkins-Smith theorem, coupled with 
the fundamental theorem of 
Hopkins, Neeman, Thomason, and others, 
reviewed in \S 2, 
which is the $\Dqc (X)$ analogue of he
Ohkawa theorem:
\begin{equation}
\label{THE DIAGRAM-introduction}
\xymatrix{
2^{| X | } 
\ar@<1mm>[rr]^-{
\left\{ Q \in \Dqc (X) \  \mid \ 
\supp (Q) \subseteq - \right\}
} & &
\ar@<1mm> [ll]^-{\supp}
\mathbb{L}( \Dqc (X) )
\\
\operatorname{Tho}( | X | )
\ar@{^{(}-_{>}}[u]
\ar@<1mm>[r]^-{\Dperf_{-}(X)} & 
\mathbb{T}\left( \Dperf(X) \right)
\ar@<1mm> [l]^-{\supp}
\ar@<1mm>[r]^-{I_X} & 
\ar@<1mm> [l]^-{C_X}
\mathbb{S}( \Dqc (X) )
\ar@{^{(}-_{>}}[u]
}
\end{equation}
This commutative diagram is very important
because it encapsulates the story 
(of not only this article, but also of this 
procedings!).
In fact, this commutative diagram
in $\Dqc^c \subset \Dqc$, which is the 
analogue of 
the commutative diagram 
in 
$\Ho^c \subset \Ho$ 
(introduced in \S 3),  
leads us to 
extend these commutative diagrams
to other triangulated categories.
Furthermore, the mutually inverse arrows at the bottom right of the diagram yield a positive solution to the telescope 
conjecture 
(see Theorem~\ref{MAIN THEOREM} and 
Remark~\ref{THE REMARK} for more detail), 
unlike the original 
problematic telescope 
conjecture in $\Ho_{(p)}$ which shows up in 
the commutative diagram  \eqref{THE DIAGRAM-SH-introduction} (see the paragraph after 
Theorem~\ref{telescope conjecture}
).
Finally, to close this section, 
we shall review Neeman's recent 
result, which claims two close relatives
$\Dcoh (X)$ and $\Dperf (X)$ actually 
determine each other, and its 
main technical tool: approximable 
triangulated category whose principal 
example is $\Dqc (X),$ as well as 
$\mcS \mcH.$

\item[\underline{\S 5}:]

Having been convinced that
$\Dcoh(X)$ and $\Dperf(X)$ carry rich 
 information  and are intimately related to 
 each other in the previous section, 
we 
review here Neeman's recent
investigations of the important 
\lq\lq strong generation\rq\rq \ 
property,  
in the sense of Bondal and Van den Bergh 
\cite{{MR1996800}}, 
for $\Dcoh(X)$ and $\Dperf(X).$
The focus here (and in this paper) is 
Neeman's $\Dcoh (X)$ strong generation 
sufficent criterion via de Jong's 
regular alteration, for which we give 
a substantial part of its proof, 
including some modifications.

\begin{itemize}
\item 
Start with the 
$\Dqc (X)$ strong compact
generation sufficient criterion
Theorem~\ref{SCG criterion}, 
and give an outline 
of its proof, emphasizing where 
the approximability of $\Dqc (X)$ is 
used 

\item 
Apllying both 
the fundamental theorem of Hopkins, Neeman, 
Thomason and others on the classification of 
thick tensor ideals of $\Dperf (X)$ and 
the Thomason localization theorem on 
$\Dperf (X \setminus Z),$ both of which 
were reviewed in \S 4, we shall show 
how the 
$\Dqc (X)$ strong compact 
generation sufficient criterion
Theorem~\ref{SCG criterion}, reviewed above,  
implies the 
$\Dqc (X)$ strong bounded 
generation sufficient criterion via 
de Jong's regular alteration 
Theorem~\ref{SCG criterion}.
Here, we extend and partially modify 
Neeman's proof 
in order to make this review 
beginner-friendly.

\item Having the 
$\Dqc (X)$ strong compact 
generation sufficient criterion available, 
we can prove our desired $\Dcoh (X)$ 
strong generation sufficient crietrion 
via de Jong's regular alteration 
Theorem~\ref{GenerationCondition}. 
However, this proof is rather involoved, 
and requires, in addition to Christensen's 
theory of phantom masp, 
some algebrao-geometric result which 
we had to put in a black box.  
We have located this black box in 
Lemma~\ref{1.8.(i)+2.6.}
(ii).

\end{itemize}

\end{description}

Neeman's own results presented in
this survey are not exactly what he talked 
about 
at the workshop.  For instance, 
although the \lq\lq strong generation\rq\rq \ 
of $\Dcoh (X)$ and $\Dperf (X)$ was still a 
major issue in Neeman's lecture series, 
Neeman's theory 
of approximable triangulated category, which 
first appeared in Neeman's series of arxiv 
preprints in 2017, was not touched upon  
during 
2015 lectures.  
Likewise, nothing was mentioned 
from \S 3 and \S 4 in this survey
during  
2015 lectures.
In contrast, Neeman actually 
talked about other results of his own
, 
but they have been omitted in this survey.
All of these decisions were made  
in order to make 
this proceedings a \lq\lq coherent story,\rq\rq
with this survey at its philosophical core.
%
%
In fact, 
the author, who  
happened to be both an organizer of the
workshop and an editor of this follow-up  
proceedings, became confident that
 the mathematics 
 presented by 
Neeman at the workshop vividly interacts
with lots of other talks at the 
workshop and articles submitted to this 
proceedings.
So, the author repeatedly mentioned such 
interactions whenever appropriate.

In spite of such an excitement, 
the first version of this 
paper was just a twenty page short list of 
 results with no proof,
\footnote{
Actually, the author thought even such 
a short list is exciting.
}
but it was the requests and 
the suggestions by Professor Neeman and 
Professor Morava, which prompted the 
author to revise this article 
repeatedly to 
contain lots of useful results, 
including many proofs!



The author would like to express his
hearty thanks to Professor Amnon Neeman 
for his beautiful lecture series, his  
encouragement to write up his lecture 
series from the author's perspective
as a non-expert, and his request to 
write a beginner-friendly survey of 
his proof of the $\Dcoh (X)$ strong generation 
sufficient criterion, in such a way that
 the roles
of the homotopical ideas of Bousfield, 
Ohkawa, Hopkins-Smith and others in its 
proof become transparent.
Not only that, Professor Neeman kindly 
read a preliminary version of this 
survey and offered the author many 
many useful suggestions including 
locating author's confusions.

The author's thanks also goes to 
Professor Jack Morava for his suggestion
to emphasize the telescope conjecture 
and the chromatic splitting in this 
article, as well as 
many inspiring and useful comments, 
some of which emerged as 
footnotes of this paper.  

The author also thanks 
Dr. Tobias Barthel for his help with 
the chromatic splitting conjecture, 
Professor Mike Hopkins for his historical comment on an earlier version of this paper,  
Professors Srikanth B. Iyengar and Ryo Takahashi 
for their information of their work, 
and  Professor 
Peter May for his comments on the definition of the tensor triangulated category and supporting our emphasis of the conjecture(s) of Hovey 
and Hovey-Palmieri. 
The author also would like to thank 
Dr. Ryo Kanda for preparing a tex file
of Professor Neeman's lecture series for us.
%

Still, the author is solely responsible 
for any left over mistakes and confusions,
 as a matter of course.

Professor Haynes Miller informed the
author of  
interesting works of 
Ruth Joachimi and Tobias Barthel,
both of which have been 
incorporated in this survey and our 
proceedings, 
As an editor of this proceedings, the author would like to thank Professor 
Miller for these information and other 
valuable information, all of which 
were so crucial in organizing this 
proceedings. 



To conclude the introduction, the author 
dedicates this survey to Professor Tetsusuke 
Ohkawa, the author's former colleague at 
Hiroshima University.  Probably the author
should express his heartfelt gratitude to
Professor Tetsusuke Ohkawa with 
rhetorical flourish...  However, the author 
does not have such an ability, and, what is 
probably even more importantly, the 
author knows very well that Professor 
Ohkawa prefers interesting mathematics 
much more than 
such rhetorical flourish!   
So, the author would like to close this 
section with a homework on behalf of 
Professor Tetsusuke Ohkawa to be submitted 
to 
Professor Tetsusuke Ohkawa :

\begin{Homework}
Extend the commutative diagrams 
below to other triangulated categories:

\begin{equation*}
\label{THE DIAGRAM-SH-introduction}
\xymatrix{
{
\text{{\it mysterious} set}
}
\ar@{=}
[rr]^-{\text{Ohkawa Th.}
} & &
\mathbb{B}( \Ho_{(p)} )
\ar@{^{(}-_{>}}[rr]^{\text{Hovey Conj.}}
_{\overset{?}=}
 & &
\mathbb{L}( \Ho_{(p)} )
\\
\underset
{\cdots \subsetneq C_{n+1}
\cdots \subsetneq C_{n} \cdots}
{
{
\text{chromatic hierachy}
}
}
\ar@{^{(}-_{>}}[u]
\ar
@{=}[rr]^-{
\text{{\tiny Hopkins-Smith Th.}
}
} & &
\mathbb{T}\left( \Ho_{(p)}^{fin} \right)
\ar@<1mm>[rr]^-{I\ \text{(split inj.)}
} &  &
\ar@<1mm> [ll]^-{C\ \text{(split surj.)}
%
}
\mathbb{S}( \Ho_{(p)} )
\ar@{^{(}-_{>}}[u]
}
\end{equation*}
\begin{equation*}
\label{THE DIAGRAM-introduction}
\xymatrix{
2^{| X | } 
\ar@<1mm>[rr]^-{
\left\{ Q \in \Dqc (X) \  \mid \ 
\supp (Q) \subseteq - \right\}
} & &
\ar@<1mm> [ll]^-{\supp}
\mathbb{L}( \Dqc (X) )
\\
\operatorname{Tho}( | X | )
\ar@{^{(}-_{>}}[u]
\ar@<1mm>[r]^-{\Dperf_{-}(X)} & 
\mathbb{T}\left( \Dperf(X) \right)
\ar@<1mm> [l]^-{\supp}
\ar@<1mm>[r]^-{I_X} & 
\ar@<1mm> [l]^-{C_X}
\mathbb{S}( \Dqc (X) )
\ar@{^{(}-_{>}}[u]
}
\end{equation*}

\end{Homework}

\section{
 Ohkawa's theorem on Bousfield classes forming a set, and its shadows in algebraic geometry
}

%
The first goal of this section is to recall
Ohkawa's theorem in 
stable homotopy 
theory.  Ohkawa's theorem claims the 
Bounsfield classes in the stable homotopy
category $\mcS \mcH$ form 
a set which is
very mysterious and beyond our imagination.
\footnote{
Concerning this sentence, Professor Morava
communicated the following thoughts
to the author:
\lq\lq
When I read it I was reminded of a quotation from the English writer Sir Thomas Browne (from `Urn Burial', in 1658):

   What song the Sirens sang, or what name Achilles assumed when he hid
   himself among women, though puzzling questions, are not beyond all
   conjecture...

I believe understanding the structure of Ohkawa's set (perhaps by defining something like a topology on it) is very important, not just for homotopy theory but for mathematics in general. An analogy occurs to me, to other very complicated objects (like the Stone-\v{C}ech compactificatin of the rationals or the reals, or maybe the Mandelbrot set) which are very mysterious but can approached as limits of more comprehensible objects. 
Indeed I wonder if this is what Neeman's theory of approximable triangulated categories points toward.\rq\rq
}

Then the second goal of this section is 
the fundamental theorem of 
Hopkins, Neeman, Thomason and others, which 
roughly states the analogue of the 
Bousfield classes in $\Dqc (X),$ in contrast to 
the Ohakawa's case of $\mcS \mcH,$ 
form a set with 
a clear algebro-geometric description.

Since both $\mcS \mcH$  and 
$\Dqc (X)$ are triangulated categories, 
we start with 
recalling some basic terminologies of 
triangulated categories.

\subsection{Bousfield localizations}

Let $\mcT$ be a triangulated category. The suspension functor is denoted by $\Sigma$. In this article all triangulated categories are assumed to have small $\Hom$-sets, 
except Verdier quotients to be defined now.

In fact, to study highly rich objects like 
 triangulated categories, we should 
\lq\lq localize\rq\rq \ at various stages.
This is exactly what Verdier \cite{Ver77} 
did in the 
context of derived categories.

\footnote{
Let us briefly recall the localization in 
the abelian category setting:
\cite[III,1]{Gab62}
\cite[p.122,Exer.9]{GM03}.
Just as we may start with thick 
triangulated categories for Verdier 
quotients, which we will see in 
Remark~\ref{localizing subcategory is thick} (iii), to localize an abelian 
category $\mcA$ by its full 
subcategory $\mcB,$ we start with 
assuming $\mcB$ is a 
\underline{\em Serre subcategory}, i.e. 
\begin{equation*}
\text{for any exact sequene}\ 
0 \to B' \to B \to B'' \to 0 \ 
\text{in $\mcA$},\qquad
\left(
B\in \mcB \iff 
(
B' \in \mcB \ \text{and}\ B'' \in \mcB 
)
\right)
\end{equation*}
Then the 
\underline{\em 
quotient category $\mcA/\mcB,$
in the sense of Gabriel, Grothendieck, 
Serre}, is of the following form:
\begin{equation*}
\operatorname{Ob} \mcA/\mcB := 
\operatorname{Ob} \mcA;
\qquad
\text{\lq\lq $\Hom$\rq\rq}_{\mcA/\mcB} 
( A, A' ) := 
\varinjlim_{\underline{A},\underline{A'}
\ \text{s.t.}\ A/\underline{A} \in \mcB,
\underline{A'} \in \mcB
}
\ 
{\Hom_{\mcA}}
\left( \underline{A}, 
A'/\underline{A'} \right)
\end{equation*}
Thus, an element of 
$\Hom_{\mcA/\mcB}( A, A' )$ is of the 
following form:
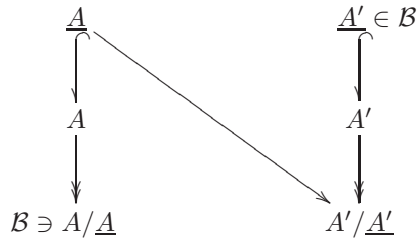
\begin{equation*}
\xymatrix{
\underline{A} 
\ar@{^{(}-_{>}}[d]  \ar[ddrr] & & 
\ \ \ \ \underline{A'} \in \mcB
\ar@{^{(}-_{>}}[d]
\\
A \ar@{->>}[d]  & & A' \ar@{->>}[d] 
\\
\mcB \ni 
A/ \underline{A} \ \ \ & & 
A' / \underline{A'} 
}
\end{equation*}
However, if we consider a similar diagram
in the setting of derived categories, 
we may take the  homotopy pullback 
$\widetilde{\underline{A}}$ as in the 
following diagram:
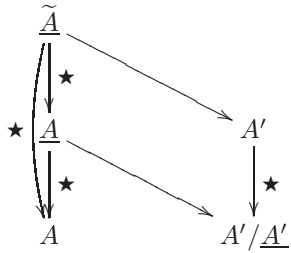
\begin{equation*}
\xymatrix{
\widetilde{ \underline{A} } 
\ar[d]^{\bigstar} \ar[drr] 
\ar@/_/[dd]_{\bigstar}
& & 
\\
\underline{A} 
\ar[d]^{\bigstar} \ar[drr] & & 
A'
 \ar[d]^{\bigstar}
\\
A & &  A'/  \underline{A'} 
}
\end{equation*}
Here, arrows with $\bigstar$ are local
maps, and so, this gives a pair of maps
$( A \xleftarrow{\bigstar} 
\widetilde{ \underline{A} }
\to A' ),$  which is a
typical element in the 
\lq\lq $\Hom$ \rq\rq class 
in the Verdier quotient.
%
}

\begin{Definition}[Verdier quotient 
(a.k.a. Verdier localization)]
{\rm \cite{Ver77} (see also 
\cite[Chapter 2]{MR1812507} )} 
For a triangulated category $\mcT$ and
its triangulated subcategory
\footnote{
\underline{WARNING!:} 
In this article, we follow the 
convention of 
\cite[Def.1.5.1]{MR1812507}
\cite[4.5]{MR2681709} for 
a \underline{\em triangulated  subcategory}, 
which is automatically full by this
convention.  On the other hand, 
it is not so in the convention of 
\cite[p.3,1.1]{MR1436741}.
}
 $\mcS,$  
the \emph{Verdier quotient (a.k.a. Verdier localization)}   
$\mcT/\mcS$ is a \lq\lq triangulated 
category\rq\rq
\footnote{
Verdier quotient does not necessarily 
have small $\Hom$-sets.
}
, which are characterized by the following
properties:

\begin{itemize}
\item $\Ob(\mcT/\mcS)=\Ob(\mcT).$ 
For $X,Y\in \Ob(\mcT/\mcS)=\Ob(\mcT),$ the
\underline{\emph{class}} of morphisms is 
given  by
\begin{equation*}
\hspace{-6mm}
''\Hom''_{\mcT/\mcS}(X,Y) = 
\frac{
\text{diagrams of the form}  \ 
( X \xleftarrow{l} Z \xrightarrow{f} Y )
\ 
\text{with}\ l,f \in \Hom_{\mcT},  
\mathrm{Cone}(l) \in \Ob(\mcS)
}
{
( X \xleftarrow{l_1} Z_1 \xrightarrow{f_1} Y_1 )
\ 
\simeq
 \ 
( X \xleftarrow{l_2} Z_2 \xrightarrow{f_2} Y_2 )
\ 
\ \iff \
\xymatrix{
& Z_1 \ar[ld]_{l_1} \ar[rd]^{f_1} & 
\\
X & \exists Z \ar[d] \ar[u] \ar[r] \ar[l] 
  & Y 
\\
& Z_2 \ar[lu]^{l_2} \ar[ru]_{f_2} & 
}
}
\end{equation*}

\item The \emph{Verdier localization functor} 
\begin{equation}
\label{Verdier localization functor}
\begin{split}
F_{univ} : \mcT &\to \mcT/\mcS  \\
       X &\mapsto X   \\
( X \xrightarrow{f} Y ) &\mapsto 
( X \xleftarrow{id_X} X \xrightarrow{f} Y )
\end{split}
\end{equation}
is universal for all triangulated functors 
$F : \mcT \to \mcT$ which sends all 
morphisms 
$(Z \xrightarrow{l} X)$ with  
$\mathrm{Cone}(l)\in \Ob(\mcS)$ to 
invertible morphisms.

\item The triangulated structure of 
$\mcT/\mcS$ is induced from 
that of $\mcT$ via 
the Verdier localization functor $F_{univ}$:
\begin{itemize}
\item The suspension $\Sigma_{\mcT/\mcS}$ 
of $\mcT/\mcS$ is induced from 
the suspension $\Sigma_{\mcT}$ 
of $\mcT$:
\begin{equation*}
\begin{split}
\Sigma_{\mcT/\mcS} : 
\mcT/\mcS &\to \mcT/\mcS 
\\
X &\mapsto \Sigma_{\mcT} X
\\
( X \xleftarrow{l} Z \xrightarrow{f} Y )
&\mapsto 
( \Sigma_{\mcT} X \xleftarrow{\Sigma_{\mcT}l} \Sigma_{\mcT}Z \xrightarrow{\Sigma_{\mcT}f} \Sigma_{\mcT}Y )
\end{split}
\end{equation*}

\item A distinguished triangle in 
$\mcT/\mcS$ is isomorphisc to the 
Verdier localization functor
$F_{univ}$ image of a 
distinguished triangle in $\mcT.$

\end{itemize}

\end{itemize}

\end{Definition}

As is always the case with such a 
localization procedure, the Verdier 
localization does not necessarily 
have small $\Hom$-sets. 
It was Neeman's insight \cite{MR1191736}
\cite{MR1308405}\cite{MR1812507} 
to make use of the Bousfield localization 
\cite{Bou79}, which was 
introduced in the context of stable 
homotopy theory, 
 to take case of this problem in 
general triangulated category theory.

To explain this theory of Neeman, 
we now prepare some definitions.

\begin{Definition}
{\rm
(\underline{WARNING!:}\
A 
\underline{\em triangulated  subcategory} 
is by definition 
\cite[Def.1.5.1]{MR1812507}
\cite[4.5]{MR2681709} automatically full.
)
}
\leavevmode
	\begin{enumerate}
        \item  A triangulated subcategory $\mcS$ of a triangulated category $\mcT$ with small coproducts is called \emph{localizing}, if it is closed under coproducts in $\mcT$.
		\item A triangulated 
subcategory $\mcS$ of 
$\mcT$ is called 
\emph{thick}, if it closed under direct summands in $\mcT$.
        \item 
{\rm \cite[p99,Rem.2.1.39]{MR1812507}} 
The \emph{thick closure}\ 
$\widehat{\mcS}$ of  a triangulated 
subcategory $\mcS$ of a triangulated 
category $\mcT$ is the triangulated 
subcategory of $\mcT$ consisting of 
direct summands in $\mcT$ of objects 
in $\mcS.$

\item 
{\rm \cite[1.4]{MR1436741} }
A triangulated 
subcategory $\mcS$ of a triangulated 
category $\mcT$ is called 
\underline{\bf dense}, if 
$\widehat{\mcS} = \mcT.$

	\end{enumerate}
\end{Definition}

\begin{Remark}
\label{localizing subcategory is thick}
{\rm (i)}
{
Every localizing triangulated subcategory is thick, for
any direct summand decomposition in
$\mcT:$ 
\begin{equation*}
S \ \ni \ x = ex \oplus (1-e)x\ 
\end{equation*}
can be realized using the cones in $\mcS:$
\begin{equation*}
\begin{cases}
ex &= \ \mathrm{Cone}
\left( \oplus_{\mathbb{N}} x \ \to \
\oplus_{\mathbb{N}} x : \  
(\xi_n)_{n\in\mathbb{N}} \mapsto
(\xi_n - e\xi_{n-1})_{n\in\mathbb{N}} 
\right)
\\
(1-e)x &= \ \mathrm{Cone}
\left( \oplus_{\mathbb{N}} x \ \to \
\oplus_{\mathbb{N}} x : \  
(\xi_n)_{n\in\mathbb{N}} \mapsto
(\xi_n - (1-e)\xi_{n-1})_{n\in\mathbb{N}} 
\right)
\end{cases}
\end{equation*}
}
{\rm (ii)} A triangulated subcategory 
$\mcS$ of a triangulated category 
$\mcT$ is thick if and only if 
$\mcS = \widehat{\mcS}.$
\newline
{\rm (iii)} 
{\rm \cite[p99,Rem.2.1.39]{MR1812507}} 
The thick clusure is 
nothing but the kernel of the 
Verdier localization functor: 
For a triangulated subcategory 
$\mcS$ of a triangulated category 
$\mcT,$ $\widehat{\mcS} = 
\Ker \left( F_{univ} : \mcT \to 
\mcT/\mcS \right).$
\newline
{\rm (iv)} {\rm 
\cite[p.148,Cor.4.5.12]{MR1812507}}
If $\mcS$ is a dense triangulated 
subcategory of a triangulated category 
$\mcT,$ then,
\begin{equation}
\forall x \in \mcT,\quad 
x\oplus \Sigma x \in \mcS.
\end{equation}

To see this,
\footnote{
If $\mcT$ is essentially small, this result 
also follows immedaitely from a general result
reviewed later in 
Proposition~\ref{dense triangulated category}.
}
 since $\exists y\in \mcT\ 
\text{s.t.}\ x\oplus y \in \mcS,$ 
form a triangle:
\begin{equation*}
x\oplus 0 \oplus y 
\xrightarrow{0\oplus 0 \oplus id_Y} 
0\oplus x \oplus y 
\xrightarrow{0\oplus id_X \oplus id_Y} 
\Sigma x \oplus x \oplus 0,
\end{equation*}
where the first and the second terms 
are contained in $\mcS$: 
$x\oplus 0 \oplus y \cong 
0\oplus x \oplus y  \cong x\oplus y
\in \mcS,$ and so is the third term: 
$\Sigma x \oplus x  \cong 
\Sigma x \oplus x \oplus 0 \in \mcS,$
as desired.

\end{Remark}

From Remark~\ref{localizing subcategory is thick} (iii),  to search for 
criteria which guarantee the 
Verdier quotient to have small 
$\Hom$-sets, we may  
start with a thick triangulated 
subcategory $\mcS$ of 
$\mcT.$ Also, while the 
original Bousfield localization \cite{Bou79}  
require $\mcT$ to have small coproducts, 
there are many cases where we wish 
Verdier quotients $\mcT/\mcS$ to have 
samall $\Hom$-sets, even when 
$\mcT$ does not have small coproducts, 
Now, Neeman \cite{MR1812507}
proposed the following general definition 
for Bousfield localization:

\begin{Definition} \label{Neeman's Bousfield localization}
{\rm
\cite[Def.9.1.1,Def.9.1.3,Def.9.1.4,Def.9.1.10]{MR1812507} 
\cite{MR2681709}  (i) }
Let $\mcS$ be a thick subcategory of 
a triangulated category $\mcT.$
\footnote{
We dot not require $\mcT$ to have 
small coproducts in this definition.
} 
Then the pair $\mcS \subset \mcT$ is said 
to possese a \underline{\em 
Bousfield localization functor} when 
the  Verdier localization  functor 
$F_{univ} : \mcT \to \mcT/\mcS$
has a right adjoint
$G:  \mcT/\mcS \to \mcT,$ 
which is called the 
\underline{\em Bousfield localization 
functor}.
The resulting composite
\begin{equation*}
L := G\circ F_{univ} : \mcT 
\xrightarrow{F_{univ}} 
\mcT/\mcS \xrightarrow{G} \mcT
\end{equation*}
is also called the 
\underline{\em Bousfield localization 
functor}
by an abuse of terminology.
\newline
{\rm (ii)}
$\mcS \subset \mcT$ is, by 
definition, the full subcategory of 
\underline{\em $\mcS$-colocal objects}.
\newline
{\rm (iii)}
$\mcS^{\perp} \subset \mcT$ is, by 
definition, the full subcategory of 
\underline{\em $L$-local objects} 
or \underline{\em $\mcS$-local objects}.
\end{Definition}

An adjoint functor between triangulated 
categories showed up in the above definition, but such an adjoint functor
actually becomes a triangulated functor:

\begin{Lemma}
\label{adjoint functor between triangulated categories}
{\rm 
\cite[Lem.5.3.6]{MR1812507}
}
Suppose a pair of adjoint functors 
between triangulated categories are given:
\begin{equation*}
\xymatrix{
\mcS \ar@<1mm>[r]^F & \mcT
\ar@<1mm>[l]^G
}
\end{equation*}
If either one of $F$ or $G$ is a 
triangulated functor, then so is 
the other.
\end{Lemma}

We shall freely use this useful fact
for the rest of this article.

Still, readers might worry that 
the more existence of a right adjoint 
$G : \mcT / \mcS \to \mcT$ 
in the definition of the above Bousfield 
localization too weak. However, in 
this particular case, we have a very 
 special property that the 
natural map from the category of 
fractions 
$\mcT \left[ \Sigma(F_{univ})^{-1} \right]
$
to the Verdier quotient $\mcT/\mcS$ 
becomes an equivalence:
\begin{equation*}
\mcT \left[ \Sigma(F_{univ})^{-1} \right] \ 
\xrightarrow{\cong} \  \mcT/\mcS,
\end{equation*} 
where $\Sigma(F_{univ})$ is the collection 
of morphisms in $\mcT$ whose image in 
$\mcT / \mcS$ is invertible, i.e. those 
maps in $\mcT$ whose mapping cone is in 
$\mcS.$ And, using this useful fact, we
can see any right adjoint 
$G : \mcT\/\mcS \to \mcT$ is fully 
faithful by applying the following useful
fact:

\begin{Lemma}
\label{fully faithful right adjoint criterion}
{\rm ( see \cite[I,Prop.1.3]{GZ67}
\cite[Prop.2.3.1]{MR2681709}).}
For an adjoint pair:
\footnote{
This is an adjoint pair of functors 
between ordinary categories, and we 
are not considering any triangulated 
structure.
}
$
\xymatrix{
\mcC \ar@<1mm>[r]^F & \ar@<1mm>[l]^G \mcD
},
$
the following conditions are equivalent:

\begin{itemize}
\item The right adjoint $G$ is fully 
faithful.
\item The adjunction $F\circ G \to 
\operatorname{Id}_{\mcD}$ is an
isomorphism.
\item The functor $\overline{F} :
\mcC \left[ \Sigma(F)^{-1} \right] 
\to \mcD$ satisfying $F = \overline{F}
\circ Q_{\Sigma (F) }$ is an equivalence,
where $\Sigma (F)$ is the collection of 
morphisms in $\mcT$ whose images in 
$\mcT'$ by $F$ becomes invertible, and 
\newline
$Q_{\Sigma (F) } : \mcC \to 
\mcC \left[ \Sigma(F)^{-1} \right]$ is 
the canonical quotient functor to the 
category of fractions.
\end{itemize}

\end{Lemma}

 Thus, from Lemma~\ref{fully faithful right adjoint criterion} and
Lemma~\ref{adjoint functor between triangulated categories}, we obtain
the following:

\begin{Proposition}
Any right adjoint $G : \mcT/\mcS \to \mcT$ 
in Neeman's definition of the Bousfield 
localiztion Definition~\ref{Neeman's Bousfield localization} is automatically 
a fully faithful triangulated functor.
\end{Proposition}

In fact, as is well known, 
if a triangulated functor 
$F : \mcT \to \mcT'$ enjoys good 
properties listed in 
Lemma~\ref{fully faithful right adjoint criterion}, then we 
have the following very useful result:
\footnote{Goes back at least to Verdier.} 
\footnote{
Let us recall the following precursor of 
this result in the setting of abelian 
categories, which goes back at least to
Gabriel (see also \cite[Lem.3.2]{Rou10}):  
If an exact functor $F : \mcA \to \mcB$ 
between abelian categories has a 
fully faithful right adjoint $G$ 
(i.e. the adjunction $F\circ G \to
\operatorname{Id}_{\mcB}$ is an 
 isomorphism, then $\Ker F$ is  Serre 
subcategory of $\mcA,$ and $F$ induces the
following equivalence of abelian
categories: $\mcA/ \Ker F \ 
\xrightarrow{\cong} \ \mcB,$
where the left hand side is the abelian
quotient category in the sense of 
Gabriel, Grothendieck, Serre.
}

\begin{Proposition}
\label{equivalence by fully faithful right adjoint or adjunction isomorphism}
{\rm (see e.g. \cite[Lem.3.4]{Rou10})}
If a triangulated functor 
$F : \mcT \to \mcT'$ has a fully faithful 
 right adjoint $G$ or a right adjoint $G$ 
with its adjunction an isomorphism 
$F\circ G \xrightarrow{\cong} 
\operatorname{Id}_{\mcD},$ then
$\Ker F$ becomes a thick triangulated
subcategory of $\mcT,$ and 
$F$ induces the following equivalence
of triangulated categories:
\begin{equation*}
\mcT/ \Ker F \ \xrightarrow{\cong} \ \mcT'
\end{equation*}
\end{Proposition}

Going back to Bousfield localization, 
we prepare some more difinitions 
to state its basic properties.

\begin{Definition}
\leavevmode
	\begin{enumerate}
        \item
{\rm ( \underline{WARNING!:} \ 
These conventions are those of 
\cite[4.8]{MR2681709}, which are the 
\underline{\em opposite}\ of 
\cite[Def.9.1.10;Def.9.1.11]{MR1812507}! )
}
	For a full subcategory $\mcA$ of $\mcT$, define the full subcategory $\mcA^{\perp}$ of $\mcT$ by
	\begin{equation*}
		\mcA^{\perp}=\conditionalset{t\in\mcT}{\Hom_{\mcT}(\mcA,t)=0}.
	\end{equation*}
	
	Dually, ${^{\perp}}\mcA$ is defined by
	\begin{equation*}
		{^{\perp}}\mcA=\conditionalset{t\in\mcT}{\Hom_{\mcT}(t,\mcA)=0}.
	\end{equation*}
		\item For full subcategories $\mcA$ and $\mcB$ of $\mcT$, denote by $\mcA*\mcB$ the full subcategory of $\mcT$ consisting of all objects $y\in\mcT$ for which there exists a triangle $x\to y\to z\to \Sigma x$ with $x\in\mcA$ and $z\in\mcB$.
		\end{enumerate}
\end{Definition}


\begin{Proposition} 
{\rm 
\cite[Prop.9.1.18;
Th.9.1.16;
Th.9.1.13;
Cor.9.1.14]{MR1812507} \cite[Prop.4.9.1]{MR2681709}
}

\label{Bousfield localization in triangulated category}
	Let $\mcS$ be a 
thick 
subcategory of a triangulated 
category $\mcT$.
\footnote{ 
We do not require
$\mcT$ to have small coproducts.
}
Then the following assertions are equivalent.

	\begin{enumerate}

\item The inclusion functor $I : \mcS 
\hookrightarrow \mcT$ has a right adjoint 
$\widetilde{\Gamma} : 
\mcT \to \mcS.$

\item $\mcT=\mcS*\mcS^{\perp}$.

\item $\mcS \subset \mcT$ posseses a
\underline{\em Boundfield localization functor}, i.e.
the Verdier localization functor
 $F_{univ} : \mcT \to \mcT/\mcS$ 
has a right adjoint $G : 
\mcT/\mcS \to \mcT.$

\item The composite
$E : \mcS^{\perp} \hookrightarrow 
\mcT \to \mcT/\mcS$ 
is an equivalence.

	\item The inclusion $J : \mcS^{\perp}\into\mcT$ has a left adjoint
$\mcT \to 
\mcS^{\perp}$ 
and  ${^{\perp}}(\mcS^{\perp})=\mcS$.

	
	\end{enumerate}
These equivalent conditions can be
succinctly expressed, via the standard 
adjoint functor notation, 
\footnote{An arrow above is left adjoint 
to the arrow below.}
as follows:
\begin{equation}
\label{Bousfield localization visualization}
\xymatrix{
\mcS  \ar@<1mm>[r]^-{I} & 
\mcT  \ar@<1mm>[r]^-{F_{univ}} 
\ar@<1mm>[l]^-{\widetilde{\Gamma}}  &  
\mcT/\mcS   \ar@<1mm>[l]^-{G}
}
\end{equation}

\end{Proposition}

\begin{Remark}\label{LocalizationFunctor}
	Assume that the inclusion $I\colon\mcS\into\mcT$ has a right adjoint $\widehat{\Gamma}$ as in 
Proposition~\ref{Bousfield localization in triangulated category}.1. Then, 
for each $t\in\mcT$, embed the counit of adjunction $\Gamma(t)=I
\widetilde{\Gamma}(t)\to t,$ 
where $\Gamma : \mcT \to \mcT$ is 
called the 
\underline{\em Bousfield colocalization functor}\ for the pair $\mcT \to \mcT/\mcS,$
\footnote{
A Bousfield colocalization functor means its 
opposite functor is a Bousfield localization functor \cite[Def.3.1.1]{HPS97}
\cite[2.8]{MR2681709}. 
\underline{\em WARNING:} This terminology 
is not consistent with that of Bousfield 
\cite{Bou79} (see 
\cite[Rem.3.1.4]{HPS97} ).

}
into a triangle
	\begin{equation*}
		\Gamma(t)\to t\to L(t)\to\Sigma\Gamma(t),
	\end{equation*}
which yields a functor $L : \mcT \to \mcT.$ 
	Then we see $L(t)\in\mcS^{\perp},$ which

\begin{itemize}
\item implies  $\mcT=\mcS*\mcS^{\perp}$ 
in Proposition~\ref{Bousfield localization in triangulated category}.2; 

\item yields a left adjoint 
$\widetilde{L} : \mcT \to \mcS^{\perp}$ 
to 
the inclusion $J : \mcS^{\perp}
\hookrightarrow \mcT,$ 
stated in 
Proposition~\ref{Bousfield localization in triangulated category}.5,, and 
$\widetilde{L}$ yields the 
\underline{\em Bousfield localization 
functor}, recovering the above 
functor $L$ by the composition
\begin{equation} \label{Bousfield localization as adjoint composite}
L = J\circ \widetilde{L} : 
\mcT \to \mcS^{\perp} \to \mcT.
\end{equation}

\item yields a left adjoint 
$G : \mcT/\mcS \to \mcT$  to the 
Verdier localization functor 
$F_{univ} ; \mcT \to \mcT/\mcS$ 
as the composition 
$G : \mcT/\mcS \xrightarrow{\overline{
\widetilde{L}}} 
 \mcS^{\perp} 
\hookrightarrow \mcT$ 
stated in 
Proposition~\ref{Bousfield localization in triangulated category}.3, and,

\item 
assuming 
Proposition~\ref{Bousfield localization in triangulated category}.4, 
$\widetilde{L}$ is equivalent to 
$E^{-1}\circ F_{univ} : 
\mcT \to \mcT/\mcS \to \mcS^{\perp}.$

\end{itemize}

\end{Remark}

{

\begin{Remark}
Actually, the property in 
Proposition~\ref{Bousfield localization in triangulated category}.2 is exactly what
Bondal-Orlov \cite[Def.3.1]{MR1957019} 
call \emph{semiorthogonal decomposition} 
and denote by 
\begin{equation} \label{semiorthogona decomposition}
\mcT =
 \langle \mcS^{\perp}, \mcS \rangle .
\end{equation}

\end{Remark}
}


Of course, the fundamental question is
when Bounsfield location exists.
Now, Neeman's insight \cite[Th.8.4.4]{MR1812507} 
is to apply Brown representability 
to construct Bousfield localization.
We now review this development 
following mostly Krause 
\cite{Kra02}\cite{MR2681709}.

\begin{Definition} 
\label{generation when there are small coproducts}
Let $\mcT$ be a 
triangulated category with small 
coproducts.
\newline
{\rm (i)\ \cite[Def.6.2.8]{MR1812507}} A set $G$ of objects in 
$\mcT$ is said to 
\emph{generate} $\mcT,$ if
$(\bigcup_{n\in\mbZ}\Sigma^{n}G)^{\perp}=0,$
i.e., 
\newline
given $t\in \mcT,$
\begin{equation*}
\forall g\in G, \forall n\in \mathbb{Z},
\ \Hom_{\mcT}(\Sigma^ng,t)=0
\qquad \implies \qquad t=0.
\end{equation*}
{\rm (ii)}
An element $t\in \mcT$ 
is called \emph{compact}\ if, for every
set of objects 
$\{ t_{\lambda} \}_{\lambda\in\Lambda}$ 
in $\mcT,$ the natural map
\begin{equation*}
\oplus_{\lambda\in\Lambda}
\Hom_{\mcT}\left(t,t_{\lambda}\right)
\to
\Hom_{\mcT}\left(t,
\oplus_{\lambda\in\Lambda}
t_{\lambda}\right)
\end{equation*}
is an isomorphism.
\newline
{\rm (iii)} \  
$\mcT$ is called \emph{compactly generated},
\ if $\mcT$ is generated by a set 
of compact objects in $\mcT.$
\newline
{\rm (iii)\ 
(c.f.\cite[Def.1]{Kra02}\cite[5.1]{MR2681709}
\footnote{
Strictly speaking, the definition here 
is slightly differently from 
Krause's, but essentially the same.
}
(see also \cite[Def.8.1.2]{MR1812507})
)
}
A set of objects $P$ in 
$\mcT$ is said to 
\emph{perfectly generate} \ $\mcT,$ if,
\begin{enumerate}
\item $P$ generates $\mcT,$
\item for every countable set of
morphisms $x_i \to y_i$ in $\mcT$
such that 
$\mcT( p, x_i ) \to \mcT(p, y_i)$ is 
surjective for all $p\in P$ and $i,$
the induced map
\begin{equation*} 
\mcT\left( p, \coprod_i x_i \right) \to 
\mcT\left( p, \coprod_i y_i \right)
\end{equation*}
is surjective.
\end{enumerate}
$\mcT$ is called \emph{perfectly generated},
\ if $\mcT$ is perfectly generated by 
a set $P$ of objects in $\mcT.$

\end{Definition}

\begin{Remark}
Any compactly generated triangulated 
category is perfectly generated.
\end{Remark}

%
%

\begin{Theorem}[Brown representability] 
\label{Brown representability}
{\rm \cite[Th.A]{Kra02}\cite[Th.5.1.1]{MR2681709} (\cite{MR1308405}  \cite{MR1812507}) }
Suppose a 
triangulated 
category $\mcT$ is perfectly generated.
\begin{enumerate}
\item A functor $F : \mcT^{op} \to Ab,$ 
the category of abelian groups,  
is cohomological and sends coproducts 
in $\mcT$ to products in $Ab$ if and 
only if 
\begin{equation*}
F \cong \mcT( - , t )
\end{equation*}
for aome object $t$ in $\mcT.$
\item A triangulated functor 
$\mcT \to \mcU$ 
preserves small coproducts
if and only if it has a right adjoint.
\end{enumerate}

\end{Theorem}

From the second part of this theorem 
and the second characterization of 
Bousfield localization in 
Proposition~\ref{Bousfield localization in triangulated category}, 
we immediately obtain the following:

\begin{Corollary}[Existence of 
Bousfield localization 
]
\label{Existence of Bousfield localization}
{\rm \cite[Prop.5.2.1]{MR2681709}
\cite[Prop.9.1.19]{MR1812507}
} 
Bousfield localization exists for any 
perfectly generated triangulated 
subcategory $\mcS$ of $\mcT,$ 
a triangulated category with small 
coproducts.
\end{Corollary}

\begin{Corollary} 
\label{Bousfield localization for a compactly generated subcategory}
Bousfield localization exists for any 
compactly generated triangulated 
subcategory $\mcS$ of $\mcT,$ 
a triangulated category with small 
coproducts.

\end{Corollary}

To be precise, the 
\lq\lq compactly generated\rq\rq \  
assumption adapted in \cite[Lem.1.7]{MR1191736} meant the smallest 
localizing triangulated subcategory containing the 
generating set is the entire triangulated
category. 
But this can be reconciled by the following 
corollary of Corollary~\ref{Existence of Bousfield localization}:

\begin{Corollary}
{\rm
\cite[Th.8.3.3;Prop.8.4.1]{MR1812507}
}
Suppose $\mcT$ is perfectly generated 
by a set $P$ of objects in $\mcT,$ then
\begin{equation*}
\mcT = 
\text{the smallest localizing 
triangulated 
subcategory containing $P$}.
\end{equation*}

\end{Corollary}

For a special case of 
Corollary~\ref{Bousfield localization for a compactly generated subcategory}, 
Neeman and Miller 
gave a simple explicit 
homotopy theoretical construction
of Bousfield localization with a nice
property:

\begin{Theorem}
\label{Nee92b, Lem.1.7}
{\rm \cite{MR1317588}
\cite[Lem.1.7]{MR1191736} }
For any localizing triangulated subcategory 
$\mcR$ of a compactly generated triangulated 
category with small coproducts $\mcT$ 
such that $\mcR$ 
is the smallest ocalizing triangulated subcategory containing a set 
$R$ 
consisting of 
compact objects in $\mcT,$ 
\begin{enumerate}
\item
Bousfield localization exists,
\footnote{This claim itself is a special 
case of Corollary~\ref{Bousfield localization for a compactly generated subcategory}.
}
given explictly by Miller's 
\underline{\em finite localization}
\cite[p.554,Proof of Lem.1.7]{MR1191736} \cite[From p.384,-6th line to p.385, 1st line]{MR1317588}:
for $x\in \mcT,$ proceed inductively 
as follows:

\begin{itemize}
\item $x_0 := x,$
\item Suppose $x_n$ has been defined, then
set
\begin{equation*}
x_{n+1} := \operatorname{Cone}
\left(
\oplus_{r\in R}
\oplus_{f_r\in \Hom_{\mcT}(r,x_n)}
r \
\xrightarrow{  \oplus_{r\in R}
\oplus_{f_r\in \Hom_{\mcT}(r,x_n)} f_r}
\ x_n
\right)
\end{equation*}

\item Then Miller's finite localization of 
$x\in \mcT$ is simply given by 
the mapping telescope:
\begin{equation*}
x \ \to \  Lx := 
\operatorname{hocolim} ( x_n ).
\end{equation*}

\end{itemize}

\item Miller's finite localization 
is \underline{\em smashing}, i.e. $L$ 
preserves arbitrary coproducts.

\end{enumerate}

\end{Theorem}

Let us record the above definition of 
\lq\lq smashing\rq\rq , because this 
definition of \lq\lq smashing\rq\rq \ 
without smash (tensor) product is 
not the traditional Ravenel's definition 
\cite{Rav84}:

\begin{Definition} 
\label{smashing localization}
{\rm \cite[5.5]{MR2681709}}
A Bousfield localization $L : \mcT \to 
\mcT$ is \underline{\em smashing} if 
$L$ preseves arbitrary coproducts in $L.$
Then, $\mcS = \Ker L$ is also called 
 \underline{\em smashing}.

\end{Definition}

We have the following equivalent
characterizations of smashing Bousfield 
localization without smash (tensor) product:

\begin{Proposition}
\label{smashing localization characterizations}
{\rm \cite[Prop.5.5.1]{MR2681709}}
For a thick subcategory $\mcS$ of a 
triangulated category with small coproducts, 
suppose there is a Bousfield localization 
$L=G \circ F_{univ} : 
\mcT \to \mcT$ 
for the pair $\mcS \to \mcT$ 
in the following set-up:
(see \eqref{Bousfield localization visualization}):
\begin{equation}
\label{Bousfield localization visualization - Re}
\xymatrix{
\mcS  \ar@<1mm>[r]^-{I} & 
\mcT  \ar@<1mm>[r]^-{F_{univ}} 
\ar@<1mm>[l]^-{\widetilde{\Gamma}}  &  
\mcT/\mcS   \ar@<1mm>[l]^-{G}
}
\end{equation}
Then the following conditions are equivalent:

\begin{enumerate}
\item Bousfield localization 
$L=G \circ F_{univ}$ is \underline{\em smashing}, i.e.  $L=G \circ F_{univ} : 
\mcT \to \mcT$ preserves coproducts 
(see 
Definition~\ref{smashing localization}
).
\item Bousfield colocalization 
$\Gamma =I\circ \widetilde{\Gamma} : 
\mcT \to \mcT$ preserves coproducts.
\item The right adjoint $G : \mcT/\mcS 
\to \mcT$ of the Verdier quotient 
$F_{univ} :  \mcT \to \mcT/\mcS$ 
preserves coproducts.
\item The right adjoint 
$\widetilde{\Gamma} : \mcT \to \mcS$ of 
the canonical inclusion $I: \mcS \to \mcT$ 
preserves coproducts.

\item 
The full subcategory $\mcS^{\perp}$ of all 
 $L$-local ($\mcS$-local) objects is 
localizing.
\end{enumerate}

If $\mcT$ is perfectly generated, 
\footnote{
This \lq\lq perfectly generated\rq\rq \ 
condition is used to apply 
Brown representability (Theorem~\ref{Brown representability}) to construct 
two right adjoints in the recollement.
}
in addition the following is equivalent.

\begin{enumerate}
\item[6.]  In the set-up 
\eqref{Bousfield localization visualization - Re}, both $\widetilde{\Gamma}$ and $G$ 
have right adjoints and 
\eqref{Bousfield localization visualization - Re} is amplified to a 
\underline{\em recollement} 
\footnote{
For the precise definition of 
recollement, consult 
\cite[1.4]{BBD82}.
}
of the following form:
\begin{equation}
\label{Bousfield localization visualization - Re}
\xymatrix{
\mcS  \ar@<5mm>[r]^-{I} 
\ar@/_/[r]
& 
\mcT  \ar@<5mm>[r]^-{F_{univ}} 
\ar@/_/[r]
\ar@<-3mm>[l]^-{\widetilde{\Gamma}}  &  
\mcT/\mcS   \ar@<-3mm>[l]^-{G}
}
\end{equation}

\end{enumerate}

\end{Proposition}

Later in 
Proposition~\ref{implication of Ravenel's smashing condition}, 
 all of these conditions are 
 shown to be equivalent to Ravenel's \cite{Rav84}, 
when $\mcT$ is a rigidly compactly generated 
 tensor triangulated category.
Smashing localization is frequently 
referred in the context of the 
telescope conjecture, which 
asks whether the converse of the second 
claim in Theorem~\ref{Nee92b, Lem.1.7} 
holds or not
\footnote{
Strictly speaking, this is the 
telescope conjecture without 
smash (tensor) product, but coincides 
with the original Ravenel's 
telescope conjecture for $\mcT=\mcS\mcH,$ 
and more generally for rigidly compactly 
generated tensor triangulated categories
\cite[Def.3.3.2,Def.3.3.8]{HPS97}
(see also
Proposition~\ref{implication of Ravenel's smashing condition}).
}
:

\begin{Conjecture}[Telescope conjecture 
 without smash (tensor) product
]
\label{telescope conjecture without product}
\ \newline
{\rm
\cite[Def.3.3.2,Def.3.3.8]{HPS97}
(see also
Proposition~\ref{implication of Ravenel's smashing condition}
)
 }
In a rigidly compactly generated 
tensor triangulated category $\mcT,$ 
a smashing 
localization $L: \mcT \to \mcT$ is 
a finite localization, 
i.e. $\Ker L$ is generated by a set of 
compact objects in $\mcT.$
\end{Conjecture}

After we take into account the tensor product strucure, we shall revisit 
the fininite localization and the 
telescope conjecture in
Theorem~\ref{telescope conjecture}.
For now, we record another easy consequnce
of Miller's finite localization 
 construction 
presented in 
Theorem~\ref{Nee92b, Lem.1.7}:

\begin{Proposition} 
\label{compact objects observation}
{\rm (See 
\cite[p.556, from 7th to 10th lines]{MR1191736}) }

Let $R$ be a set of compact objects in 
a triangulated category with 
small coproducts $\mcT,$ and  
$\mcR$ 
 be the 
smallest localizing triangulated 
subcategory 
containing $R.$ 

Then, every element in $\mcR^c$ is 
isomorphic in $\mcR^c$ to a direct summand 
of a finite extensions of finite coproducts 
of elements in $R.$
In particular, $\mcR^c$ is 
essentially small.

\end{Proposition}

In fact, for any $X\in \mcR,$ 
the Bousfield localization
with respect to the pair 
$\langle R\rangle = \mcR \subset \mcR,$
is trivial for any $x\in \mcR$:
\begin{equation*}
 x \to Lx  := \operatorname{hocolim} (x_n)
\simeq 0.
\end{equation*}
Then, if $x\in \mcR^c,$ this map 
becomes trivial at some \lq\lq finite\rq\rq 
\ stage, which implies $x$ is a direct 
 summand 
of a finite extensions of finite coproducts 
of elements in $R,$ as claimed.

\subsection{Bousfield classes and Ohkawa's theorem}

Now we focus on a special case: let $\mcT=\Ho$ be the homotopy category of spectra. Then $\mcT$ is a triangulated category with coproducts. It has the smash product $\wedge\colon\mcT\times\mcT\to\mcT$ and the unit object $S^{0}\in\mcT$ which make $\mcT$ a tensor triangulated category. 
\footnote{
For a serious treatment of the definition of 
\lq\lq tensor triangulated category,\rq\rq\  
consult \cite{May01}.
}

The smash product preserves coproducts in each variable. $\mcT$ is generated by $\set{S^{0}}$, and $\mcT$ satisfies Brown representability.

For each $H\in\mcT$, put $H_{*}=H\wedge (-)$. We consider the localizing 
triangulated subcategory
\begin{equation*}
	\Ker H_{*}=\conditionalset{t\in\mcT}{H\wedge t=0},
\end{equation*}
which is called the \emph{Bousfield class} of $H$.

\begin{Theorem}[Bousfield \cite{Bou79}]
	Let $\mcT=\Ho$ be the homotopy category of spectra.
	\begin{enumerate}
		\item If $\mcS\subset\Ho$ is a localizing triangulated
subcategory which is generated by a \emph{set} of objects, then a Bousfield localization exists for $\mcS$.
		\item For every $H\in\Ho$, there exists a \emph{set} of objects which generates $\Ker H_{*}$. Therefore a Bousfield localization exists for $\Ker H_{*}$.
	\end{enumerate}
\end{Theorem}

Now, we can state the truly surprising 
theorem of Ohkawa
\footnote{Somewhat surprisingly, 
Ohkawa's theorem had been elusive from researchers' attention for more than a decade. 
It was the paper of Dwyer and Palmieri 
\cite{DP01} which drew researchers' attention 
to Ohkawa's surprising theorem.}
\footnote{
For a concise summary of the academic life
of Professor Tetsusuke Ohkawa, see 
\cite{MatPr} in this proceedings.
}
:

\begin{Theorem}[Ohkawa \cite{MR1035147}] \label{Ohkawa}
	$\conditionalset{\Ker H_{*}}{H\in\Ho}$ is a \emph{set}.
\end{Theorem}

We note that no explicit structure of this 
set is known.

For more detials, including a
proof, of the Ohakawa theorem, see
 the survey \cite{Carles} in this proceedings.

\subsection{
Casacuberta-Guti\'{e}rrez-Rosick\'{y} 
theorem, 
motivic analogue of Ohkawa's theorem }

Ohkawa's theorem is a statement in 
the stable homotopy category  
$\mcS \mcH,$ 
which is \lq\lq a part\rq\rq\ 
 of the Morel-Voevodsky stable homotopy 
 category $\mcS \mcH (k)$ when 
$k\subseteq \mathbb{C},$ via the 
retraction of the following form:

\begin{equation}
\label{MotivicDominates}
\xymatrix{
\mcS\mcH 
\ar[r]
\ar@/_/@<-2ex>[rr]_-{id} & 
\mcS\mcH (k)
\ar[r]^{R_k}
&
\mcS\mcH 
}
\end{equation}

So, a natural question here is whether 
there is a shadow of Ohkawa's theorem in this
algebro-geometrical setting, i.e. whether  
there is a motivic analogue of Ohkawa's 
theorem or not. 

Now, 
Casacuberta-Guti\'{e}rrez-Rosick\'{y} 
\cite{CGR14}
answered this question affirmatively 
under some very mild assumption.

\begin{Theorem} \label{MotivicOhkawa}
{\rm \cite[Cor.3.6]{CGR14}} 
For each Noetherian scheme $S$ of finite Krull dimension, there is only a set of
distinct Bousfield classes in the stable motivic homotopy category $\mcS \mcH (S)$ 
with base scheme $S.$
\end{Theorem}

Once again,  no explicit structure of this 
set is known.

{
For  various generalizations of 
Ohkawa's theorem, see afore-quoted 
\cite{CGR14}, also \cite{KOSSPr} and 
the review \cite{CRPr}, 
both in this proceedings.

\subsection{Localizing tensor ideals of derived categories and the fundamental theorem of
Hopkins, Neeman, Thomason and others
}

In both Ohkawa's Theorem~\ref{Ohkawa} and 
its algebro-geometric shadow 
Theorem~\ref{MotivicOhkawa}, the 
resulting sets
are completely mysterious and beyond our 
imagination.  However, if we take a look at 
 the algebro-geometrical 
shadow of Ohkawa's theorem from a different 
angle, i.e. by 
considering $\Dqc (X)$ for a fixed 
Noetherian scheme instead of 
$\mcS \mcH (k),$ then we see an explicit set 
representing clear algebro-geometric 
information.  This is the fundamental theorem
of Hopkins, Neeman, Thomason, and others, 
which has been the guiding principle of 
the area.

Now, the tensor structure is essential 
for this fundamental theorem, and we 
must start with some review of 
fundamental facts about 
general tensor triangulated categories 
and Bousfield localization from the
tensor triangulated category point of view.

\begin{Definition}
\label{tensor ideal, prime, strongly dualizable}
	Let $\mcT$ be a tensor triangulated category. 
\begin{enumerate}
\item
A triangulated subcategory $\mcI$ of $\mcT$ is called a 
{
$\begin{cases}
\text{\emph{tensor ideal} } \\
\text{\emph{prime} }
\end{cases}$
if \ 
$$
\begin{cases}
\mcT\otimes\mcI\subset\mcI; \\
\text{it is a tensor ideal and}\ 
\left(\mcT \setminus \mcI \right) \otimes
\left(\mcT \setminus \mcI \right) \subset
\left(\mcT \setminus \mcI \right) \neq
\emptyset.
\end{cases}$$
}

\item {\rm \cite[Chapter III]{LMS86} \ 
(see also \cite[App.A,]{HPS97} 
\cite[p.1163]{MR2806103}
)} 
An element $x$ in a closed symmetric 
 monoidal triangulated category 
$(\mcT, \otimes, \underline{\Hom})$ 
is called 
\underline{\em strongly dualizable}  
or simply \underline{\em rigid},
\footnote{
If $x\in \mcT$ is strongly dualizable, i.e. rigid, the natural map $x \to D^2x$ 
is an isomorphism \cite[Chapter III]{LMS86}
\cite[Th.A.2.5.(b)]{HPS97}. 
}
if the natural map
$Dx \otimes y \to
\underline{\Hom}(x,y),$ where
$Dx :=  \underline{\Hom}(x, \one),$ 
is an isomorphism for all $y\in \mcT.$

\item {\rm \cite[Def.1.1.4]{HPS97} \
(see also \cite[Hyp.1.1]{MR2806103})} 
A closed symmetric 
 monoidal triangulated category 
$(\mcT = \langle G\rangle, \otimes, \underline{\Hom})$ 
is called a
\newline
\underline{\em unital
algebraic stable homotopy category} \ or a 
\newline
\underline{\em rigidly compactly 
 generated tensor triangulated category}, if
$\one$ is compact and 
$\mcT = \langle G\rangle$ for a set $G$ 
of rigid and compact objects.
\footnote{
In a rigidly compactly generated tensor 
triangulated category, any compact object 
is rigid, for, by 
Proposition~\ref{compact objects observation}, any compact object 
is seen to be isomorphic to a direct summand 
of a finite extensions of finite
coproducts of rigid elements.
In particular, in 
 a rigidly compactly generated tensor 
triangulated category, $\one$ is both 
rigid and compact.
}

\end{enumerate}

\end{Definition}

%



Now, we are ready to reconcile our 
previous definition (Definition~\ref{smashing localization}) 
of smashing localization with Ravel's original definition in \cite{Rav84} for 
rigidly compactly generated tensor triangulated categories:

\begin{Proposition}
\label{implication of Ravenel's smashing condition}
{\rm \cite[Def.3.3.2]{HPS97}} 
For a thick subcategory $\mcS$ of a 
closed symmetric 
 monoidal triangulated category with 
small coproducts 
$(\mcT = \langle G\rangle, \otimes, \underline{\Hom})$ 
\footnote{
Recall in this case $\mcT$ becomes
\underline{\em distributive}, because 
for any objects 
$x_{\lambda}\ (\lambda\in\Lambda), y, z$ 
in $\mcT,$
$
\Hom\left( ( \oplus_{\lambda} x_{\lambda} )
\otimes y, z \right) \cong 
\Hom\left( \oplus_{\lambda} x_{\lambda}, 
\underline{\Hom}(y,z) \right)
\cong \prod_{\lambda} \Hom \left( 
x_{\lambda}, 
\underline{\Hom}(y,z) \right) \cong 
\prod_{\lambda} \Hom ( 
x_{\lambda}\otimes y, z ) \cong 
\Hom \left( \oplus_{\lambda}
x_{\lambda}\otimes y, z \right).$ 
}, 
suppose there is a Bousfield localization 
$L : \mcT \to \mcT$ for the pair 
$\mcS \to \mcT.$ Consider the 
following \lq\lq smishing\rq\rq\ 
conditions:
\begin{itemize}
\item[{\rm (S):}] \quad (Ravenel's original definition 
of smashing localization \cite{Rav84}): 
\begin{equation*}
L\cong L(\one)\otimes-, \ 
\text{
where $\one$ 
is the unit onject of $(\mcT, \otimes).$
}
\end{equation*}

\item[{\rm (C):}]\quad 
(The definition of smashing localization 
in Definition~\ref{smashing localization}):
\begin{equation*}
\text{$L$ preserves arbitrary coproducts.}
\end{equation*}

\end{itemize}

Then, the implication 
{\rm (S)} $\implies$ {\rm (C)} always holds.
If $\mcT$ is also a 
rigidly compactly generated tensor triangulated category, 
the converse 
{\rm (C)} $\implies$ {\rm (S)} also holds,
and so, {\rm (C)} and {\rm (S)} become 
equivalent.

\end{Proposition}

\begin{proof}
The implication 
{\rm (S)} $\implies$ {\rm (C)} is easy:
\begin{equation*}
L\left( \oplus_{\lambda} x_{\lambda} \right)
\overset{(S)}\cong 
 L(\one)\otimes 
\left( \oplus_{\lambda} x_{\lambda} \right)
\cong \oplus_{\lambda}
\left(
 L(\one)\otimes x_{\lambda}
\right) \overset{(S)}\cong  
\oplus_{\lambda} Lx_{\lambda}.
\end{equation*}
For the converse 
{\rm (C)} $\implies$ {\rm (S)}, first note
that {\rm (C)} implies those $x\in \mcT$ 
which satisfies $L\one\otimes x \cong Lx$ 
form a localizing triangulated subcategory 
of $\mcT,$ even without the rigidly 
compactly gnerated assumption. For instance,
if $L\one\otimes x_{\lambda} \cong 
 Lx_{\lambda}\ \forall \lambda\in\Lambda,$
then
\begin{equation*}
 L(\one)\otimes 
\left( \oplus_{\lambda} x_{\lambda} \right)
\cong \oplus_{\lambda}
\left(
 L(\one)\otimes x_{\lambda}
\right) \cong  
\oplus_{\lambda} Lx_{\lambda}
\overset{(C)}\cong 
L\left( \oplus_{\lambda} x_{\lambda} \right).
\end{equation*}
Now, we are reduced to showing 
$L\one\otimes g \cong Lg$ for any rigid 
element $g.$ For this, we start with the 
tensor product of the localization 
distinguished sequence for $\one$ with $g$:
\begin{equation*}
\Gamma ( \one )\otimes g 
\to \left( g \cong \one \otimes g \right)
\to L(\one)\otimes g,
\end{equation*}
and apply the Bousfield localization $L$ 
to drive the equivalence 
$L(\one)\otimes g \cong Lg$ as follows:

\begin{equation*}
\begin{split}
&\left( * \overset{(TI)}{\cong} 
L( \Gamma (\one) \otimes g ) 
\right) 
\to 
\left( L g \cong 
L ( \one \otimes g) \right) 
\\
\xrightarrow[\because) (TI)]{\cong}
&\left( L( L(\one) \otimes g) 
\overset{(R)}\cong 
L \underline{\Hom}( Dg, L (\one) )
\overset{(L)}\cong 
 \underline{\Hom}( Dg, L (\one) ) 
\overset{(R)}\cong 
L(\one)\otimes g
\right),
\end{split}
\end{equation*}
where $(TI)$ holds because $\Ker L$ is 
a tensor ideal, (R) holds because 
$g$ is rigid, and (L) holds because 
$\underline{\Hom}( Dg, L (\one) )$ is 
$L$-local. 

\end{proof}

In general, when we talk about  
smashing Bousfield localization in  
tensor triangulated setting, we adopt 
the following equivalent conditions, 
where the localizing tensor ideal 
$\mcI$ is  called a
\underline{\em smashing ideal} \ 
\cite[Def.2.15]{MR2806103}:

\begin{Proposition}[See 
{\cite[Th.2.13]{MR2806103}}]
\label{smashing-BF11}
	Let $\mcT$ be a tensor triangulated category with coproducts, and let $\mcI$ be a localizing 
tensor ideal of $\mcT$ for which a Bousfield localization exists. Define the Bousfield localization functor $L\colon\mcT\to\mcI^{\perp}$ as in Remark~\ref{LocalizationFunctor}. Then the following assertions are equivalent.
	\begin{itemize}
		\item[(TI)] \quad $\mcI^{\perp}$ is a tensor ideal. That is, $\mcT\otimes\mcI^{\perp}\subset\mcI^{\perp}$.
		\item[(S)]  \quad 
$L$ is smashing in Ravenel's sense: \  
$L\cong L(\one)\otimes-$.
	\end{itemize}
\end{Proposition}

\begin{Remark}

{\rm (TI)} is a tensor triangulated 
analogue of Proposition~\ref{smashing localization characterizations}.5.

\end{Remark}

\begin{proof}[Proof of 
Proposition~\ref{smashing-BF11}]
Now, for the implication 
{\rm (TI)} $\implies$ {\rm (S)}, 
consider the 
tensor product of the localization 
distinguished sequence for $\one$ with $x
\in \mcT$:
\begin{equation} \label{product of LDT with x}
\Gamma ( \one )\otimes x 
\to \left( x \cong \one \otimes x \right)
\to L(\one)\otimes x,
\end{equation}
where $\Gamma ( \one )\otimes x \in \mcI$ 
 because $\mcI$ is a tensor ideal by 
 assumption, and 
$L(\one)\otimes x \in \mcI^{\perp}$ 
 because $\mcI^{\perp}$ is also 
a tensor ideal by (TI).
Then, from the uniqueness of the 
localization distinguished sequence for 
$x\in \mcT,$ we find $Lx \cong 
L(\one)\otimes x,$ which implies (S).

The converse 
{\rm (S)} $\implies$ {\rm (TI)} is easy;  for, if $l = L(l)\in \mcI^{\perp}$ be 
a $\mcI$-local object and $x\in \mcT,$ then
\begin{equation*}
l\otimes x = L(l)\otimes x 
\overset{(S)}= ( L(\one)\otimes l )\otimes x
= L(\one)\otimes ( l\otimes x ) 
\overset{(S)}= L( l\otimes x ) \in \mcI^{\perp}.
\end{equation*}

\end{proof}

In the above proposition, we started with 
a localizing tensor ideal for which a 
Bousfield localization exists. 
However, we have the following example of 
a localizing tensor ideal for which an
existence of the Bousfield localization 
is problematic:

\begin{Example}
	Let $\mcT=\Ho$ be the homotopy category of spectra. For every $H\in\mcT$, its Bousfield class $\Ker H_{*}$ is a localizing tensor ideal. The subcategory
	\begin{equation*}
		\Ker H^{*}=\conditionalset{t\in \mcT}{\text{$\Hom(t,\Sigma^{i}H)=0$ for all $i\in\mbZ$}},
	\end{equation*}
called the 
\emph{
cohomological Bousfield class}
of $H,$ 
	is also a localizing tensor ideal.
Actualy, as was noticed by Hovey 
\cite[Prop.1.1]{Hov95b}, any Bousfield 
class is a cohomological Bousfield class:
\begin{equation*}
\Ker H_{*} = \Ker (IH)^{*},
\end{equation*}
where $IH$ is the \emph{Brown-Comenetz 
dual} of $H,$ charaxterized by:
$
(IH)^*(t) = \Hom \left(H_*(t), 
\mathbb{Q}/\mathbb{Z} \right),\ 
\forall t\in \mcT.
$
\end{Example}

Here, Hovey \cite{Hov95b} and 
Hovery-Palmieri \cite{HP99} proposed 
the following conjectures, any one of which
implies that an arbitrary localizing tensor 
 ideal $\Ker H^{*}$ admits a Bousfield
 localization: 
\footnote{
Conjecture~\ref{localizing conjecture}  should
be taken more seriously.
In fact,
Professor Peter May is very glad to see 
Conjecture~\ref{localizing conjecture} 
is advertised here.   
}

\begin{Conjecture} \label{localizing conjecture}
{\rm (i)} 
{\rm
\cite[Conj.1.2]{Hov95b}
}
Every cohomological Bousfield class is
a Bousfield class.
\newline
{\rm (ii)}
Every localizing tensor ideal
 is a Bousfield class.
\newline
{\rm (iii)} 
{\rm
\cite[Conj.9.1]{HP99}
}
Every localizing triangulated subcategory 
 is a Bousfield class.
\end{Conjecture}
Of course,
$(iii) \ \implies \ (ii) \ \implies \ (i),$
for we have an obvious inclusions of classes:
\begin{equation*}
\begin{split}
&\text{Bousfield-Ohkawa set} :=
\text{The class of Bousfield classes}
\\
\subseteq \ 
&\text{The class of cohomological 
Bousfield classes}
\\
\subseteq \ 
&\text{The class of localizing tensor ideals}
\\
\subseteq \ 
&\text{The class of 
localizing triangulated subcategories},
\end{split}
\end{equation*}
where all the inclusings become $=$ if 
the above conjecture (iii) holds.
However, even (i) is still open, and so 
it is still unknown even whether the 
the class of cohomological Bousfield 
classes becomes a set or not. 
Similarly, it is still unknown 
even whether any cohomological Bousfield
class admits a Bousfield localization
or not.
Here, we shall show an analogue of (ii)
 holds holds with an explicit geometric 
description of its set structure for 
$\Dqc (X).$ 

For a scheme $X$, 
$\Dqc(X)$ 
is the derived category of complexes of arbitrary modules on $X$ whose cohomologies are quasi-coherent. If $X$ is quasi-compact and separeted, then $\Dqc(X)$ is equivalent to $\DD(\QCoh X)$, where $\DD(\QCoh X)$ is the derived category of complexes of quasi-coherent sheaves on $X$ (\cite[Corollary~5.5]{MR1214458}). 
Here we have the nice theorem of 
Gabriel \cite{Gab62} and 
Rosenberg \cite{Ros04}: 

\begin{Theorem} \label{Gabriel-Rosenberg}
Any quasi-compact and separated scheme $X$ 
can be reconstructed from $\QCoh X.$
\end{Theorem}

Glancing at this theorem of Gabriel and Rosenberg, we naturally hope 
$\Dqc(X) \cong \DD(\QCoh X)$ would carry 
rich information of $X.$

Now $\Dqc(X) \cong \DD(\QCoh X)$ is a 
 tensor triangulated category  with coproducts, with respect to 
the derived tensor product $-\otimes_{X}^{\Ld}-,$ which is 
defined using flat resolutions 
(see e.g. \cite[(2.5.7)]{Lip09}), and 
the unit object given by the structure sheaf $\mcO_{X}$.
Let us also recall the following standard
facts about derived functors:

\begin{Proposition} 
\label{derived functors in derived categories}
\ 
\newline
{\rm (i) (see e.g. \cite[(2.1.1)(2.7.2)(3.1.3)(3.9.1)(3.6.4)${}^*$]{Lip09} )}
For any map of schemes$f : X\to Y,$ 
we can define the derived pullback 
triangulated functor
\begin{equation*}
\Ld f^* : \Dqc (Y) \to \Dqc (X),
\end{equation*}
via flat resolutions. 

Furthermore, 
we have a natural functorial 
isomorphism
\begin{equation*}
\Ld f^* \Ld g^* \ \xrightarrow{\sim} \ 
\Ld (gf)^*
\end{equation*}
{\rm (ii) (see e.g. \cite[(2.1.1)(2.3.7)(3.1.2)(3.9.2)(3.6.4)${}_*$]{Lip09})} 
For any  quasi-compact and 
quasi-separated map of schemes
$f : X\to Y,$ we can define the
derived direct image (a.k.a. derived 
pushforward) triangulated functor
\begin{equation*}
\Rd f_* : \Dqc (X) \to \Dqc (Y),
\end{equation*}
via injective resolutions.
Furthermore, 
we have a natural functorial 
isomorphism
\begin{equation*}
\Rd (gj)_* \ \xrightarrow{\sim} \ 
\Rd g_* \Rd f_* ,
\end{equation*}
when both $f$ and $g$ and quasi-compact
and quasi-separated maps.
\newline
{\rm (iii) 
(see e.g. \cite[(3.6.10)]{Lip09})} 
For any  quasi-compact and 
quasi-separated map of schemes
$f : X\to Y,$ $(\Ld f^*, \Rd f_*)$ 
gives an afjunction pair:
\begin{equation*}
\xymatrix{
\Dqc (Y ) \ar@<1mm>[r]^-{\Ld f^*} &
\Dqc (X) \ar@<1mm>[l]^-{\Rd f_*}
}
\end{equation*}
{\rm (iv) (see e.g. \cite[(3.2.1)(3.9.4)]{Lip09})}
For any  quasi-compact and 
quasi-separated map of schemes
$f : X\to Y,$ the projection formula holds, 
i.e. we have natural isomorphisms
for any $F\in \Dqc (X), G\in \Dqc (Y)$: 
\begin{equation*}
(\Rd f_*F)\otimes^{\Ld} G 
\ \xrightarrow{\cong} \ 
\Rd f_*\left( F \otimes^{\Ld} \Ld f^* G \right),\quad
G\otimes^{\Ld} \Rd f_* F 
\ \xrightarrow{\cong} \ 
\Rd f_* \left( \Ld f^*G\otimes^{\Ld} F
\right)
\end{equation*}

\end{Proposition}

To investiate an analogue of Ohkawa's
theorem for $\Dqc (X),$ we must consider
localizing tensor ideals of $\Dqc(X).$ 
However, those smashing (localizing tensor
ideals) are sometimes, more important.
To stufy such (smashing) localizing tensor 
ideals of $\Dqc (X),$ an appropriate 
concept of \lq\lq stalk\rq\rq \ 
becomes crucial:

\begin{Definition}
\label{support definitions}
{\rm (compare with \cite[Proof of Th.4.12]{MR2071654} \cite[App.A]{ILN15})
\footnote{
Our presentation of \lq\lq supports\rq\rq \ 
in this definition and next proposition 
is  somewhat different 
from those given in \cite[Proof of Th.4.12]{MR2071654} \cite[App.A]{ILN15}, but 
the author hopes this would be more 
transparent to the reader.
}
} 
Let $x\in X$ be a point in a scheme.
Then we have the following canonical maps 
involving the local ring $\mcO_{X,x}$ and
the residue field $k_x$ at $x\in X$:
\begin{equation}
\label{local diagram}
\hspace{-5mm}
\xymatrix{
\Spec k_x   \ar@/^2pc/[rr]^{i_x} 
\ar[r]_-{r_x} & \Spec \mcO_{X,x} 
\ar[r]_-{l_x}^-{\text{\rm flat}} & X  
}
\ \implies \ 
\xymatrix{
\Dqc ( \Spec k_x )  
 & 
\Dqc \left( \Spec \mcO_{X,x} \right)
\ar[l]^-{\Ld (r_x)^*}
 & \Dqc (X )  
\ar[l]^-{\Ld (l_x)^* = (l_x)^*}
 \ar@/_2pc/[ll]_-{\Ld (i_x)^*} 
}
\end{equation}
Then, for $E\in \Dqc  (X),$ we have 
four notions of \lq\lq supports\rq\rq : 
\begin{equation} \label{two supports}
\begin{split}
&\supp (E) := 
\left\{ x \in X \ \mid \ 
\Ld ( i_x )^* E \neq 0 
\in \Dqc ( \Spec k_x ) \right\}
\\
\subseteqq \quad  
&\Supp (E) := 
\left\{ x \in X \ \mid \ 
 ( l_x )^* E \neq 0 
\in \Dqc ( \Spec \mcO_{X,x}  ) \right\};
\\
&\supph (E) := 
\left\{ x \in X \ \mid \ 
\Ld ( i_x )^* \left( 
\oplus_{\bullet \in \mathbb{Z}}
\mcH^{\bullet}E \right) \neq 0 
\in \Dqc ( \Spec k_x ) \right\}
\\
\subseteqq \quad  
&\Supph (E) := 
\left\{ x \in X \ \mid \ 
\oplus_{\bullet \in \mathbb{Z}}
(\mcH^{\bullet}E)_x =
( l_x )^* \left( 
\oplus_{\bullet \in \mathbb{Z}}
\mcH^{\bullet}E \right)
 \neq 0 
\in \QCoh ( \Spec \mcO_{X,x}  ) \right\}
\end{split}
\end{equation}
where:
\begin{itemize}
\item $\mcH^{\bullet}E$ is the associated
homology sheaves, regarded as a chain
complex with trivial boundries, of $E.$ 

\item the
inclusive relations follow from 
$\Ld ( r_x )^* 
 ( l_x )^*
=
\Ld ( r_x )^* \Ld ( l_x )^*
\xrightarrow{\sim} \Ld ( l_xr_x ) = 
\Ld ( i_x )^*,$ where the former 
equality follows rrom $\Ld (l_x)^* = (l_x)^*,$
a consequence of the flatness of $l_x,$ 
and the latter isomorphism is a direct 
consequence of Proposition~\ref
{derived functors in derived categories}(i).

\item these inclusing relations become
equalities when $E\in \Dcoh (X)$
because of Nakayama's lemma.

\item If it becomes necessary to  distinguish 
these four concepts,
we call $\supp (E)$ the 
\underline{\em small support of $E$},
$\Supp (E)$ the \underline{\em 
large support of $E$}. 
$\supph (E)$ the 
\underline{\em small homology support of $E$},
$\Supph (E)$ the 
\underline{\em 
large homology support of $E$}. 
Otherwise, we simply call
$\supp (E)$ the 
\underline{\em support of $E$}, because 
this is the most essential object, 
and $\Supph (E)$ he 
\underline{\em homology support of $E$}, 
because this is a tractible ordinary 
sheaf theoretical support for the 
associated homology sheaves 
$\oplus_{\bullet\in \mathbb{Z}} 
\mcH^{\bullet}E.$

\end{itemize}

\end{Definition}

Then the following useful fact will be
used later:

\begin{Proposition}
\label{support definitions conincide for coherent complexes}
{\rm (i)} Given $E \in \Dqc (X),$ 
we have for any $x\in X$ and 
$\bullet \in \mathbb{Z},$
\begin{equation*}
\mcH^{\bullet}
\left( (l_x)^*E \right) \ \cong \ 
(l_x)^* \left( \mcH^{\bullet} E \right) 
\ \in \QCoh 
\left( \Spec \mcO_{X,x} \right).
\end{equation*}
Consequently, for any $E \in \Dqc (X),$ 
\begin{equation*}
\Supp E = \Supph E.
\end{equation*}
{\rm (ii)} The commutative diagram of
quasi-coherent sheaves in 
\eqref{local diagram} restricts to 
coherent sheaves, and 
for any  $E \in \Dcoh (X),$ 
all the four concepts of supports in 
Definition~\ref{support definitions} 
coincide:
\begin{equation*}
\supp E = \Supp E = \Supph E = \supph E.
\end{equation*}

\end{Proposition}

\begin{proof}
In view of Definition~\ref{support definitions}, 
we only have to verify the first claim
in (i):
$\mcH^{\bullet}
\left( (l_x)^*E \right) \ \cong \ 
(l_x)^* \left( \mcH^{\bullet} E \right) 
\ \in \QCoh 
\left( \Spec \mcO_{X,x} \right).$ 
However, this follows immediately from
the flatness of $l_x$ which implies
$(l_x)^*$ preserves exactness at the
cochain level.
\end{proof}

Now, the fundamental theorem 
of Hopkins, Neeman, Thomason and others
classify (smashing) localizing tensor
ideals of $\Dqc (X)$ under a mild assumption
of $X$:

\begin{Theorem}
\label{Ohkawa for D_{qc}(X)}
{\rm
(\cite{MR0932260}
\cite[Th.2.8,Th.3.3]{MR1174255}
\cite{MR1436741},
\newline
\cite[Cor4.6;Cor.4.13;Th.5.6
]{MR2071654}
\cite[Cor.6.8]{MR2806103}
\cite[Cor.6.8;Ex.6.9]{DS13}
\cite[Th.B]{MR3797596})
}
	Let $X$ be a 
{N}oetherian scheme. Then every localizing tensor ideal of $\Dqc(X)$ is of the form
	\begin{equation*}
		\Ker H_{*}=
\conditionalset{Q\in 
{
\Dqc (X)
}
}{ \supp Q\subseteq S},
	\end{equation*}
for some $S\subset X$. 
	
	The subcategory $\Ker H_{*}$ is smashing if and only if the corresponding $S\subset X$ is closed under specialization.
\end{Theorem}

{
Note those $S$'s with
$S \subset X$ clearly form a set.
So, we see an analogue of Ohkawa's theorem, 
however 
with a clear algebro-geometrical 
interpretation of \lq\lq the 
Bousfield-Ohkawa set\rq\rq \  
in contrast to the case of Ohkawa's theorem. 
Furthermore, Theorem~\ref{Ohkawa for D_{qc}(X)} 
solves
Conjecture~\ref{localizing conjecture} (ii) 
affirmatively for the case $\Dqc(X).$

Also note that, 
in the special case when $S$ in 
Theorem~\ref{Ohkawa for D_{qc}(X)} is
$Z = X \setminus U \subset X,$
the complement of 
a quasi-compact Zarisiki open immersion
$j : U \hookrightarrow X,$
we have the following equivalence
for not only 
noetherian, but also 
more general 
quasicompact, quasiseparated schemes
(in which case, as
$\Ld j^*$ has a right adjoint $\Rd j_*$ 
with 
$\epsilon : \Ld j^* \Rd j_*  
\to \operatorname{id}$ 
an isomorphism, 
we may  apply
Proposition~\ref{equivalence by fully faithful right adjoint or adjunction isomorphism}
):
\footnote{So, should had been known to 
Verdier.}
\footnote{
Let us recall the following precursor of 
this result in the setting of abelian 
category of quasi-coherent sheaves, 
which should go back at 
least to Gabiriel (see e.g. 
\cite[In the proof of Prop.3.1]{Rou10}):
%
$
\QCoh (X) \big/ \QCoh_Z (X) 
\
\xrightarrow[\cong]{\overline{j^*}} 
\
\QCoh (U),
%
$
where the left hand side is the abelian
quotient category in the sense of 
Gabriel, Grothendieck, Serre.
}
%
\begin{equation}
\label{BL of TC}
\Dqc (X) \big/ \left( {\Dqc} \right)_Z (X) 
\
\xrightarrow[\cong]{\overline{\Ld j^*}}
\
\Dqc (U),
\end{equation}
where 
$\left( {\Dqc} \right)_Z (X) := 
\left\{ Y \in \Dqc (X) \ \mid \ 
\Supp Y \subseteq Z \right\} = 
\mathrm{Ker}\ \Ld j^*.$
\footnote{
Unlike Theorem~\ref{Ohkawa for D_{qc}(X)}
stated under the noetherian assumption, 
\eqref{BL of TC} is stated under more 
general quasicompact, quasiseparated 
assumption. Therefore, in this equality
$\left( {\Dqc} \right)_Z (X) := 
\left\{ Y \in \Dqc (X) \ \mid \ 
\Supp Y \subseteq Z \right\} = 
\mathrm{Ker}\ \Ld j^*,$ we may not replace
$\Supp$ with $\supp.$ In fact, without
the noetherian hypothesis, 
Theorem~\ref{Ohkawa for D_{qc}(X)} becomes 
very bad as was shown in \cite{Nee00}.  
The author is grateful to Professor Neeman 
for this reference.
}
In this generality
of quasicompact, separated schemes, 
Bousfield localization $L$ 
is smashing 
(see Proposition~\ref{smashing-BF11}), 
given explicitely as follows:
\begin{equation}
\label{smashing Bousfield localization of Dqc(X)}
\hspace{-7mm}
L 
\overset{\eqref{Bousfield localization as adjoint composite}}= 
\Rd j_* \Ld j^* 
=
\left( \Rd j_*\mathcal{O}_U \right)
\otimes_{\mathcal{O}_X}^{\Ld} - :
 \ 
\Dqc (X) \ \to \ 
\Dqc (X) \big/ \left( {\Dqc} \right)_Z (X) 
\
\xrightarrow[\cong]{\overline{\Ld j^*}}
\
\Dqc (U)
\xrightarrow{\Rd j_*}
\Dqc (X).
\end{equation}

\section{
Hopkins-Smith theorem and its motivic analogue
}

In reality, Hopkins was not motivated by 
Ohkawa's Theorem~\ref{Ohkawa} for his 
influential paper in algebraic geometry  \cite{MR0932260} 
(Theorem~\ref{Ohkawa for D_{qc}(X)}).
Instead, Hopkins was motivated by 
his own theorem with Smith 
\cite{MR1652975} 
in the
 sub stable homotopy category 
 $\mcS \mcH^c, $ 
consisting of compact objects,  
whose validity was already known to them 
back around the time Hopkins wrote
\cite{MR0932260}.


\begin{Theorem}
\label{Hopkins-Smith}
{\rm
\cite{MR1652975}
}
For any prime $p,$  any
thick (\'{e}paisse) subcategories of
the subtriangulated category 
$\Ho_{(p)}^c$ consisting
of compact objects
\begin{equation*}
\Ho_{(p)}^c = 
\Ho_{(p)}^{fin}
 = 
\text{the homotopy category 
of $p$-local finite spectra}
\end{equation*}
is of the form
\begin{equation} 
\label{Hopkins-Smith thick subcategory}
\begin{split}
\mcC_n &:=  \Ker E(n-1)_*\big|_{\Ho_{(p)}^{fin}}
 =
\left\{ X \in 
\Ho_{(p)}^{fin}
 \mid 
E(n-1)\wedge X = 0 \right\}
\\
&= \Ker K(n-1)_*\big|_{\Ho_{(p)}^{fin}} =
\left\{ X \in 
\Ho_{(p)}^{fin}
 \mid 
K(n-1)\wedge X = 0 \right\}.
\end{split}
\end{equation}
Furthermore, these form a decreasing
filtration of $\mcF_{(p)}$:
\begin{equation}
\label{Hopkins-Smith filtration}
\{ * \} \subsetneq \cdots \subsetneq
\mcC_{n+1} \subsetneq \mcC_n 
\subsetneq \mcC_{n-1} \subsetneq \cdots
\subsetneq \mcC_1 \subsetneq \mcC_0 = 
\Ho_{(p)}^{fin}
.
\end{equation}

\end{Theorem}

In this Hopkins-Smith classification of 
thick triangulated subcategories of 
$\mcS \mcH^c,$ the first step is an easy
observation that any thick triangulated 
subcategory of $\mcS \mcH^c$ is a 
thick (tensor) ideal
\footnote{
Such a property is not 
usually satisfied for general 
triangulated categories. So, most effort 
to generalize the Hopkins-Smith theorem 
for a general triangulated category $\mcT$ 
 aim at
a classification of thick (tensor) ideals 
of $\mcT^c.$
}, 
Furthermore, 
$E(n-1)$ and $K(n-1)$ are the $(n-1)$-st 
Johnson-Wilson spectrum and  Morava 
$K$-theory, respectively, and 
the equality 
$\Ker E(n-1)_*\big|_{\Ho_{(p)}^{fin}} = 
\Ker K(n-1)_*\big|_{\Ho_{(p)}^{fin}}$ in 
\eqref{Hopkins-Smith thick subcategory} and  
the inclusions 
\eqref{Hopkins-Smith filtration} are 
consequences of the following 
results found in Ravenel's paper 
\cite{Rav84}:

\begin{Theorem} 
\label{Ravenel's theorems}
{\rm (i)
\cite[Th.2.1(d)]{Rav84} } \quad 
$\Ker E(n-1)_* = \Ker 
\left( \vee_{0\leq i\leq n-1} 
K(i) \right)_*$
\newline
{\rm (ii)
\cite[Th.2.11]{Rav84}} 
For $X \in \Ho_{(p)}^{fin},$ if 
$K(i)_*X = 0,$ then $K(i-1)_*X = 0.$
\end{Theorem}

By the Hopkins-Smith work \cite{MR1652975}, 
the smashing conjecture for $E(n)$ 
\cite{Rav84} also holds \cite{Rav92}, and 
so, $\Ker E(n-1)_*$ in 
\eqref{Hopkins-Smith thick subcategory} is 
a smashing tensor ideal. Actually, the 
first equality in 
\eqref{Hopkins-Smith thick subcategory} 
is a part of the following 
elegant reformulation of the telescope
conjecture \cite{Rav84}\cite[Def.3.3.8]{HPS97}  
(see also Conjecture~\ref{telescope conjecture without product}) 
by 
Miller \cite{MR1317588} 
\cite[Th.3.3.3]{HPS97} (here we 
follow more recent formulations of 
\cite[Th.4.1;Def.4.2]{MR2806103} 
\cite[Cor.2.1;Def.3.1]{MR3797596}.):


\begin{Theorem} \label{telescope conjecture}
{\rm (Miller's finite localization and
the Ravenel telescope conjecture)}\ 
Let $\mcT$ be a rigidly
compactly generated tensor triangulated category.
Let $\mathbb{S}(\mcT)$ denote the collection of all smashing 
localizing tensor ideals of $\mcT$, and let 
$\mathbb{T}(\mcT^c)$ denote the collection of all thick tensor ideals of 
$\mcT^c.$
\newline
{\rm (i)} 
{\rm
\cite[Cor.6;Prop.9]{MR1317588} 
\cite[Th.1.7]{MR3797596}
}
For any $C \in \mathbb{T}(\mcT^c),$ 
the smallest localizing 
triangulated
subcategory 
$\langle C \rangle$ 
containing $C$ 
in $\mcT$ is  smashing, i.e. 
$
\in \mathbb{S}( \mcT ).$
Thus, we obtain the 
\underline{\em inflation map}:
$$
I : \mathbb{T}( \mcT^c ) \to 
\mathbb{S}( \mcT ).
$$
{\rm (ii)} {\rm \cite[Th.3.3.3]{HPS97}}\ 
There is also the 
\underline{\em contraction map}:
$$
C: \mathbb{S}( \mcT ) 
\to 
\mathbb{T}( \mcT^c );\quad 
S 
\mapsto S\cap \mcT^c,
$$
which enjoys:
\begin{equation*}
\xymatrix{
C\circ I = id_{\mathbb{T}(\mcT^c)} :\ 
\mathbb{T}(\mcT^c) 
\ar@<1mm>[r]^-{I} & 
\ar@<1mm> [l]^-{C}
\mathbb{S}(\mcT )
}
\end{equation*}
{\rm (iii)} 
{\rm
\cite[Cor.6;Prop.9;Cor.10]{MR1317588} 
}
The \underline{\em telescope conjecture
for $S \in \mathbb{S}(\mcT)$} holds
 if and only if, in addition to 
$C\circ I = id_{\mathbb{T}(\mcT^c)}$ 
stated in (ii), 
the following also holds:
\begin{equation*}
I\circ C (S) = S \in \mathbb{S}(\mcT)
\end{equation*}
{\rm (iv)}
{\rm
\cite[Cor.6;Prop.9;Cor.10]{MR1317588} 
}
The \underline{\em telescope conjecture
for $\mcT$}
\footnote{
This telescope conjecture is equivalent to
the {telescope conjecture without product}
Conjecture~\ref{telescope conjecture without product} via 
Proposition~\ref{implication of Ravenel's smashing condition} and 
Proposition~\ref{smashing-BF11}.
}
 holds
 if and only if $I$ and $C$ give mutually 
inverse equivalence:
\begin{equation*}
\xymatrix{
C\circ I = id_{\mathbb{T}(\mcT^c)} :\ 
\mathbb{T}(\mcT^c) 
\ar@<1mm>[r]^-{I} & 
\ar@<1mm> [l]^-{C}
\mathbb{S}(\mcT )
\ : id_{\mathbb{S}(\mcT)} = I \circ C.
}
\end{equation*}

\end{Theorem}

However, the telescope conjecture of this
generality has been shown to be false 
\cite{Kel94}, and even the original 
 telescope conjecture for $\mcS \mcH$ 
is now believed to be false by many 
experts, including Ravenel himself 
\cite{MRS01}. Still, algebraicists have shown the 
validity of its various algebraic 
analogues  
(e.g.
\cite{MR2806103}
\cite{KS10} \cite{BIKP18})
as we shall review an algebraic analogue
of the Hopkins-Smith theorem, 
in conjunction with the above 
telescope conjecture, later in 
Theorem~\ref{MAIN THEOREM}.
Furthermore, Krause \cite{Kra00} showed 
the underlying philosophical message of 
the telescope conjecture that smashing 
tensor ideals are completely characterized 
by their restrictions to 
compact objects 
In fact, whereas the original telescope 
conjecture only concerns local compact objects, Krause proves his characaterization 
of smashing tensor ideals via 
\lq\lq local maps\rq\rq between compacts objects. For more details, consult 
Krause's own paper \cite{Kra00}.

Going back to the Hopkins-Smith theorem, 
a major part of its proof was to show:

\begin{Theorem}
\label{any thick subcategory is of the form Cn}
{\rm \cite[Th.7]{MR1652975}}
Any thick subcategory of 
$\Ho_{(p)}^{fin}$ is of the form 
$\mcC_n$ for some 
$n \in \mathbb{Z}_{\geq 0}.$
\end{Theorem}

To show this, Hopkins-Smith prepared the 
following version of the niloptency 
theorem \cite[Cor.2.5]{MR1652975}, 
building upon their earlier collaboration 
work with Devinatz \cite{DHS88}: 

\begin{Theorem}
\label{Hopkins-Smith nilpotency}
{\rm \cite[Cor.2.5.ii)]{MR1652975}}
For a map $f: F\to A$ between finite 
$p$-local spectra and another finite
$P$-local spectra $Y,$ the following 
conditions are equivalent:
\begin{itemize}
\item $\exists m \gg 0$ such that
\quad
$
f^{\wedge m}\wedge I_Y : F^{\wedge m}\wedge Y \to 
A^{\wedge m}\wedge Y 
$
\quad 
is null.

\item $0\leq \forall n < \infty,\quad 
K(n)_* \left( f\wedge I_Y \right) = 0.$

\end{itemize}

\end{Theorem}

Now, to prove Theorem~\ref{any thick subcategory is of the form Cn}, it 
suffices to prove the following:

\begin{Lemma}
\label{Hopkins-Smith (2.9)}
{\rm \cite[(2.9)]{MR1652975}}
Let $\mcC$ be a thick subcategory of 
$\Ho_{(p)}$ and $X, Y$ be $p$-local 
finite spectra.   Then, if
$X \in \mcC$ and $\{ n\in \mathbb{Z}_{\geq 0} 
\mid K(n)_*Y \neq 0 \} 
\subseteq 
\{ n\in \mathbb{Z}_{\geq 0} 
\mid K(n)_*X \neq 0 \},$ then $Y\in \mcC.$

\end{Lemma}

Actually, if Lemma~\ref{Hopkins-Smith (2.9)} 
is shown to be correct, together 
with Ravenel's 
Theorem~\ref{Ravenel's theorems}(ii),  
it would imply 
\begin{equation*}
\mcC = \mcC_{m},\ \text{where}\ 
m = \min \{ n \in \mathbb{Z}_{\geq 0} \mid 
\mcC_n \subseteq \mcC \}.
\end{equation*}
Then, the proof of Lemma~\ref{Hopkins-Smith (2.9)} in \cite{MR1652975} proceeds 
as follows (see also \cite{Rav92}:

\begin{itemize}
\item Starting with $X,$ let $e : S^0 \to 
X\wedge DX$ be the $S$-dual of the 
identity map $: I_X : X\to X,$ and extend 
it to a triangle with a map between 
$p$-local finite spectra 
$f : F \to S^0$ as the fiber 
 as follows:
\begin{equation}
\label{S-duality triangle}
F \xrightarrow{f} S^0 
\xrightarrow{e} X\wedge DX
\simeq C_{f}, 
\text{the cofiber of $f.$}
\end{equation}

\item Applying the smash product with 
$Y$ to \eqref{S-duality triangle}, 
we obtain:
\begin{equation}
\label{F tensored S-duality triangle}
F\wedge Y \xrightarrow{f\wedge I_Y} 
 S^0\wedge Y \cong Y 
\xrightarrow{e\wedge I_Y} 
X\wedge DX\wedge Y \simeq C_{f}\wedge Y,
\end{equation}
for which, we claim
\begin{equation}
\label{induced homomorphisms of Morava K-theories are trivial}
0\leq \forall n<\infty,\quad 
K(n)_*(f\wedge I_Y)=0.
\end{equation}

\begin{itemize}
\item If $K(n)_*Y =0$ then 
$K(n)_*(I_Y)=0,$ which implies the 
triviality of 
\eqref{induced homomorphisms of Morava K-theories are trivial}, by the
Kunneth theorem for 
Morava $K$-theories: 
\begin{equation}
\label{Morava Kunneth}
K(n)_*(X\wedge Y) \cong K
(n)_*X\otimes_{K(n)_*} K(n)_*Y \quad
\text{for any $p$-local spectra $X,Y$}
\end{equation}

\item If $K(n)_*Y \neq 0$ then 
 $K(n)_*X \neq 0$ by the assumption 
of Lemma~\ref{Hopkins-Smith (2.9)}.
Then, by the duality isomorphism for Morava 
$K$-theories:
\begin{equation}
\hspace{-2mm}
\begin{split}
\label{duality isomorphism for Morava 
K-theorie}
\Hom_{K(n)_*}( K(n)_*X, K(n)_*Y) &= 
\Hom_{K(n)_*}\left( K(n)_*, 
K(n)_*(Y\wedge DX) \right) =
K(n)_*(Y\wedge DX)
\\
&\text{
for any $p$-local spectra $X,Y,$ 
}
\end{split}
\end{equation}
we also find the non-triviality: 
$K(n)_*(e\wedge I_Y) \neq 0.$
But, this in turn implies  the triviality:
$K(n)_*(f\wedge I_Y)=0$ from the 
Morava $K$-theory exact sequence 
associated to
\eqref{F tensored S-duality triangle}, 
making use of the Morava Kunneth isomorphism  \eqref{Morava Kunneth} again,

\end{itemize}

\item Since \eqref{induced homomorphisms of Morava K-theories are trivial}, we may apply
the {Hopkins-Smith nilpotency} 
Theorem~\ref{Hopkins-Smith nilpotency} to 
$f\wedge I_Y$ in 
\eqref{F tensored S-duality triangle} to
find $m \gg 0$ such that
$
f^{\wedge m}\wedge I_Y : F^{\wedge m}\wedge Y \to 
(S^0)^{\wedge m}\wedge Y \cong Y
$
\ 
is null.  This implies:

\begin{equation}
\label{direct summand}
\text{
 $Y$ is a direct 
summand of $C_{f^{\wedge m}\wedge I_Y} 
\cong C_{f^{\wedge m}}\wedge Y.$}
\end{equation}

\item By the assumption, $X\in \mcC,$ but 
as the thick subcategory $\mcC$ of $\Ho_{(p)}^{fin}$ is also a thick ideal, 
this implies 
$C_f
\overset{\eqref{F tensored S-duality triangle}}\cong X\wedge DX
\in 
\mcC.$

\item For any $n\in \mathbb{N},$ 
consider the commutative diagram:
\begin{equation*}
\xymatrix{
F^{\wedge n}\wedge F  
\ar[r]^-{f^{\wedge n}\wedge I_F
} \ar@{=}[d]  &
(S^0)^{\wedge n}\wedge F  \cong F
\ar[d]^{f} \ar[r] &
C_{f^{\wedge n}} \wedge F 
 \ar[d]
\\
F^{\wedge (n+1)}  \ar[d]
\ar[r]^-{f^{\wedge (n+1)}
} &
(S^0)^{\wedge (n+1)} \cong S^0
\ar[r] \ar[d] &  C_{f^{n+1}} \ar[d]
\\
\bullet  \ar[r] &  C_f
\ar@{=}[r]   &  C_f
}
\end{equation*}
From this, we obtain a triangle
\begin{equation*}
%
C_{f^{\wedge n}} \wedge F 
\to
 C_{f^{n+1}}
\to
C_f 
\end{equation*}
Since $C_f \in \mcC$ and $\mcC$ is a 
tensor ideal, we see inductively 
from this triangle that
\begin{equation} 
\label{C_f^n is in C}
C_{f^m} \in \mcC \qquad 
(\forall m\in \mathbb{N})
\end{equation}

\item Since $\mcC$ is a thick ideal, 
we conclude from 
\eqref{direct summand} and 
\eqref{C_f^n is in C} that 
$Y \in \mcC.$  This complete the 
proof of Lemma~\ref{Hopkins-Smith (2.9)}.

%
%
\end{itemize}
{\flushright
\qed
}

Now, the basic philosophy 
underlying the above picture of 
Hopkins-Smith was already perceived by 
Morava much earlier (see 
the \lq exercises\rq \ 
in \S 0.5 of 
\cite{Mor85}, whose preprint version 
was circulated nearly a decade ago 
before its publication).
%
%
For a modern development of 
Morava $K$-theory,
consult Morava's own paper \cite{MorPr} 
in this proceedings.

The author believes the 
Hopkins-Smith theorem 
(Theorem~\ref{Hopkins-Smith}) and 
the Ohkawa theorem (Theorem~\ref{Ohkawa})
are best understood, 
when they are appreciated simultaneously 
in a single commutative diagram.  
Since this commutative diagram can be 
drawn for more general rigidly
compactly generated tensor triangulated 
category $\mcT,$ 
let us first set up our notations of our 
interests in this generality:

\begin{itemize}
\item
$\mathbb{L}(\mcT)$: the collection of 
{\em
localizing
tensor ideals
}
of 
$\mcT.$

\item
$\mathbb{S}(\mcT)$: 
the collection of 
{\em smashing 
localizing tensor ideals}
of $\mcT.$ 

\item
$\mathbb{T}(\mcT^c)$: 
the collection of 
{\em thick 
tensor ideals} of 
$\mcT^c.$ 

\item
$\mathbb{B}(\mcT)$:  
the collection of 
{\em Bousfield classes}, i.e. those of 
the form 
\newline
$\Ker (h\otimes - ) 
\subseteq \mathbb{L}(\mcT)
\  (h\in \mcT).$

\end{itemize}

Now let us specialize to the case 
$\mcT = \Ho_{(p)}$:

\begin{Theorem}
\label{The diagram for SH}
In $\Ho_{(p)},$ 
the Ohkawa theorem, the Hopkins-Smith theorem, 
Miller's version of the Ravenel telescope conjecture \ 
($C\circ I \overset{?}= Id_{\mathbb{T}\left( \Ho_{(p)}^{fin} \right)}$),
 and the conjectures of Hovey 
and Hovey-Palmieri can be simultaneously 
expressed in the following succint 
commutative diagram:
\begin{equation}
\label{THE DIAGRAM-SH}
\xymatrix{
{
\text{{\it mysterious} set}
}
\ar@{=}
[rr]^-{\text{Ohkawa Th.}
} & &
\mathbb{B}( \Ho_{(p)} )
\ar@{^{(}-_{>}}[rr]^{\text{Hovey Conj.}}
_{\overset{?}=}
 & &
\mathbb{L}( \Ho_{(p)} )
\\
\underset
{\cdots \subsetneq C_{n+1}
\cdots \subsetneq C_{n} \cdots}
{
{
\text{chromatic hierachy}
}
}
\ar@{^{(}-_{>}}[u]
\ar
@{=}[rr]^-{
\text{{\tiny Hopkins-Smith Th.}}
} & &
\mathbb{T}\left( \Ho_{(p)}^{fin} \right)
\ar@<1mm>[rr]^-{I\ \text{(split inj.)}
} &  &
\ar@<1mm> [ll]^-{C\ \text{(split surj.)}
%
}
\mathbb{S}( \Ho_{(p)} )
\ar@{^{(}-_{>}}[u]
}
\end{equation}

\end{Theorem}

For more on the Hopkins-Smith theorem and 
related \lq\lq chromatic mathematics,\rq\rq 
\ 
see \cite{Rav92} and, for some of the latest 
developments, 
\footnote{
A trend here is to apply the higher algebra 
technique of Lurie \cite{Lur09}\cite{Lur16} 
to understand chromatic phenomena 
\cite{BRPr} and \cite{TorPr}, where the 
latter contains a concise review of higher 
 algebra technology.
Different kinds of applications of 
Lurie's higher algebra 
technique can be seen in \cite{M1Pr,M2Pr}.
}
see   
\cite{BarPr}\cite{BRPr}\cite{TorPr} in this proceedings.
Actually, Bartel's survey 
\cite{BarPr} focuses upon
the telescope \cite{Rav84}\cite{Rav92} 
and chromatic splitting conjectures \cite{Hov95a}, which are major directions of research, not only in 
chromatic homotopy theory, but also 
 in stable homotopy theory as a whole.
Considering the traditional influence of 
stable homotopy theory, initiated by 
Hopkind, Rickard, 
Neeman, Thomason and others, 
 to the represetation
theory of finite dimensional algebras and 
the derived category theory in algebraic, 
as is highlighed by 
Brown representability, Bousfield localization, 
Hopkins-Smith theorem, 
researchers in these areas might 
better to keep this fact in mind.


Comparing with the telescope conjecture, 
the chromatic splitting conjectre appears 
to be elusive  
for them.  
%
In short, the chromatic splitting 
 conjecture predicts, 
for a $p$-completed finite spectrum $F,$ 
the first map in 
the canonical cofiber sequence
\begin{equation} \label{HCSC}
\underline{\Hom}
\left( L_{E(n-1)}S^0, L_{E(n)} F \right) 
 \to L_{E(n-1)}F \to L_{E(n-1)}L_{K(n)}F
\end{equation}
is trivial; stated differently, the 
second map in \eqref{HCSC} is split
injective. 
\footnote{
This splitting conjecture 
implies, for any $p$-completed 
fintie spectrum $F$ and any infinite  
subset $\{ n_i \}_{i=1}^{\infty} \subseteq
\mathrm{N},$ the natural map
$F \to \prod_{i=1}^{\infty} L_{K(n_i)}F$ 
is split injective.  
For this and much more, consult 
\cite{Hov95a}\cite{BarPr}.
}

In fact, Hopkins 
\cite[Conj.4.2(iv)]{Hov95a} 
further predicted,  presumably hoping 
to provide a program to prove the 
 triviality of the first map in 
 \eqref{HCSC}, 
 an explicit decomposition of 
$\underline{\Hom}
\left( L_{E(n-1)}S^0, L_{E(n)} F \right),$ 
inspired by Morava's old observation 
\cite[Rem.2.2.5]{Mor85}. 
The strucutre of 
$\underline{\Hom}
\left( L_{E(n-1)}S^0, L_{E(n)} F \right)$ 
is highly reflected by its divisible 
homotopy group elements.  
In general, divisible homootopy group
elements of a spectrum $X$ can be 
isolated in the spectrum 
$\underline{\Hom}( L_0S^0, X),$ which is
in the current case:
\begin{equation*}
\begin{split}
\underline{\Hom}\left( L_0S^0, 
\underline{\Hom}
\left( L_{E(n-1)}S^0, L_{E(n)} F \right)
\right)
&\cong
\underline{\Hom}\left( L_0S^0\wedge 
L_{E(n-1)}S^0, L_{E(n)} F \right)
\\
&\cong 
\underline{\Hom}\left( L_0S^0, L_{E(n)} F \right)
\end{split}
\end{equation*}
To understand this, Morava \cite{Mor14} 
suggested to consider the 
following cohomology theory $L_n^*$:
\begin{equation*}
X \mapsto L_n^*(X) := 
\Hom \left( 
\pi_{-*} 
\underline{\Hom}\left( L_0S^0, L_{E(n)} X \right), \mathbb{Q} \right)
\end{equation*}
Actually, 
Morava \cite{Mor14} noticed the validity 
of the Hopkins' prediction on the explicit
structure of \newline
$\underline{\Hom}
\left( L_{E(n-1)}S^0, L_{E(n)} F \right)$ 
would imply the cohomology theory $L_n^*$ 
is represented by the $p$-adic 
rationalization of the spectrum:
\footnote{It appears that 
\cite[p.4,Corollary]{Mor14} should 
be modified as in
\eqref{L_n-rationally-representing-spectrum}}

\begin{equation} \label{L_n-rationally-representing-spectrum}
\Sigma^{2n} 
\left(
\bigvee_{
\left\{ n_i \in 
\mathbb{Z}_{\geq 0}
\right\}_{i=1}^{\infty}; \ 
%
\sum_{i=1}^{\infty} in_i = n
}
\frac{ \left( \sum_{i=1}^{\infty} n_i \right)! }
{ \prod_{i=1}^{\infty} (n_i!) }
\ 
\left(
\prod_{i=1}^{\infty} U(i-1)^{n_i}
\right)_+
\right)
\end{equation}
While Hopkins' prediction 
\cite{Hov95a} above of the 
explicit decomposition of 
\newline
$\underline{\Hom}
\left( L_{E(n-1)}S^0, L_{E(n)} F \right),$ 
which the above work of Morava 
\cite{Mor14} is based upon, 
is known to hold for $n=1$ or $n=2$ and 
$p\geq 3,$ Beaudry \cite{Bea17} has 
recently shown it to fail for the case 
$n=2$ and $p=2.$ Still, as was pointed 
out to the author by Tobias Barthel, 
The above formula 
\eqref{L_n-rationally-representing-spectrum}, which was derived from Morava's 
calculation, still holds even for 
this troublesome case of
$n=2$ and $p=2,$ because the 
descrepancy found by Beaudry \cite{Bea17} 
is $p$-torson and so vanishes rationally.
Thus, it could well 
be the case \eqref{L_n-rationally-representing-spectrum} holds 
for any pair of a prime $p$ and a natural
number$n.$ 

Furthermore, 
it could be the case that Hopkins' 
prediction of the 
explicit decomposition of \newline
$\underline{\Hom}
\left( L_{E(n-1)}S^0, L_{E(n)} F \right)$ 
still holds, 
consequently so does 
Morava's deduction
\eqref{L_n-rationally-representing-spectrum} above , when the base prime $p$ 
is sufficiently large comparing with the 
height $n.$

It would be fantastic, if, as  
Professor Morava dreams of, there 
hold formulae analogous to the predicted 
Hopkins' and Morava's.
in algebraic 
examples like $\Dqc (X),$ where the 
fundamental theorem of Hopkins, Neeman, 
Thomason and others gave us an explict
 \lq\lq Bousfield-Ohkawa set\rq\rq ,\ 
not only for Bousfield classes, 
but also for localized tensor ideals, 
 whereas the original Ohkawa's set for   
$\mcS \mcH$ only takes into account  Bousfield classes 
and is not explicit at all.
Furthermore, as we mentioned before, 
while the telescope conjecture is 
now believed to be false by many 
experts, algebraicists have shown the 
validity of its various algebraic 
analogues. 
So, why not for the chromatic splitting 
conjeture, as Professor Morava dreams of!

Actually, restricring to the conjectured  
splitting of the second map in 
\eqref{HCSC}, recent effort of 
Beaudry-Goerss-Henn \cite{BGH17} 
has shown its validity 
even for the case $n= p =2,$ which is the
 case \cite{Bea17} showed Hopkins' 
conjectural decomposition of 
$\underline{\Hom}
\left( L_{E(n-1)}S^0, L_{E(n)} F \right)$ 
is false. 
Furthermore, Barthel-Heard-Valenzuela 
\cite{BHV18} has recently proved an
algebraic analogue of 
the conjectural 
splitting of the second map in 
\eqref{HCSC}.
For this and much more, consult 
Bartel's survey \cite{BarPr}.


Going back to the Hopkins-Smith theorem, 
it is natural to look after 
its  
motivic analogue 
\eqref{MotivicDominates}
(This means efforts to classify thick 
(tensor) ideals of $\mcS \mcH (k)^c.$).

In this regard, Ruth Joachimi \cite{Joachimi} 
 constructed some motivic thick ideals in 
$\mcS \mcH (k)^c$ for $k \subseteq 
\mathbb{C}$:

\begin{Theorem}
{\rm
\cite[Th.13]{Joachimi}
}
\begin{enumerate}
\item
If $k\subseteq \mathbb{C},$
then
$\left( \mcS\mcH (k)^c \right)_{(p)}$ 
contains at least an infinite chain
of different thick ideals, given by
$\overline{R}_k^{-1}(\mcC_n), 
0\leq n\leq \infty,$ 
where $\overline{R}_k$ denotes the
$p$-localisation of the restriction
of $R_k$ to $\mcS\mcH (k)^c$:

\begin{equation*}
\xymatrix{
\left( \mcS\mcH^c 
\right)_{(p)} \ar[d]
\ar[r]_{c_k}
\ar@/^/@<2ex>[rr]^-{id} & 
\left( \mcS\mcH (k)^c 
\right)_{(p)}  \ar[d]
\ar[r]_{\overline{R}_k}
&
\left( \mcS\mcH^c 
\right)_{(p)}  \ar[d]
\\
\left( \mcS\mcH 
\right)_{(p)}
\ar[r]
\ar@/_/@<-2ex>[rr]_-{id} & 
\left( \mcS\mcH (k)
\right)_{(p)}
\ar[r]^{R_k}
&
\left( \mcS\mcH 
\right)_{(p)}
}
\end{equation*}

Here,
\begin{itemize}
\item $c_k$ is induced from the
constant presheaf functor 
\cite[Th.10]{Joachimi}, which
resticts to the compact objects 
\cite[Rem.53, Prop.58, Prop.61]{Joachimi}.

\item The existence of $\overline{R}_k$ follows since
$R_k$ preserves compactness 
\cite[Prop.61]{Joachimi}.
\end{itemize}

\item If $k\subseteq \mathbb{R},$
then $\left( \mcS\mcH (k)^c 
\right)_{(p)}$ contains at least a
two-dimensional lattice of different
thick ideals, given by 
$\left( \overline{R}_k' \right)^{-1}
(\mcC_{m,n}),$ for all $(m,n) \in
\Gamma_p$ ( see \cite[Def.35]{Joachimi} for the definition of 
$\Gamma_p$ and more detail):

\begin{equation*}
\xymatrix{
\left( \mcS\mcH (\mathbb{Z}/2)^c 
\right)_{(p)} \ar[d]
\ar[r]
\ar@/^/@<2ex>[rr]^-{id} & 
\left( \mcS\mcH (k)^c 
\right)_{(p)}  \ar[d]
\ar[r]_{\overline{R}_k'}
&
\left( \mcS\mcH (\mathbb{Z}/2)^c 
\right)_{(p)} \ar[d]
\ar@< 2pt> [r]^-{\phi^{\{1 \}}} 
\ar@< -2pt> [r]_-{\phi^{ \mathbb{Z}/2}} 
& 
\left( \mcS\mcH )^c
\right)_{(p)} \ar[d]
\\
\left( \mcS\mcH (\mathbb{Z}/2) 
\right)_{(p)}
\ar[r]^{c_k'}
\ar@/_/@<-2ex>[rr]_-{id} & 
\left( \mcS\mcH (k)
\right)_{(p)}
\ar[r]^{R_k'}
&
\left( \mcS\mcH (\mathbb{Z}/2) 
\right)_{(p)}
\ar@< 2pt> [r]^-{\phi^{\{1 \}}} 
\ar@< -2pt> [r]_-{\phi^{ \mathbb{Z}/2}}  & 
\left( \mcS\mcH )
\right)_{(p)}
}
\end{equation*}

Here, 
\begin{itemize}

\item {\rm \cite[Th.11]{Joachimi}}
$c_k' : 
\left( \mcS\mcH (\mathbb{Z}/2) 
\right)_{(p)}
\to \left( \mcS\mcH (k)
\right)_{(p)}$ 
is induced by
\begin{equation*}
\begin{split}
c' : sSet ( \mathbb{Z}/2 )
&\to sPre ( 
\operatorname{Sm} / \mathbb{R} 
)
\\
M &\mapsto
\left( 
\coprod_{M^{\mathbb{Z}/2}}
\right)
\coprod
\left(
\coprod_{(M\setminus 
M^{\mathbb{Z}/2} ) / 
 ( \mathbb{Z}/2 ) }
\Spec \mathbb{C} \right), 
\end{split}
\end{equation*}
which
resticts to the compact objects 
\cite[Rem.53, Prop.58, Prop.61]{Joachimi}.

\item {\rm (Strickland's theorem 
\cite[Cor.34]{Joachimi}
\footnote{
Strickland's theorem for $G=\mathbb{Z}/2$ 
has recently been generalized to arbitrary 
 finite group $G$ by 
Balmer-Sanders \cite{BS17}.
}
 )}
 Any thick ideal in the category 
$\left( \mcS\mcH (\mathbb{Z}/2)^c 
\right)_{(p)}$ is of the form
\begin{equation*}
\mcC_{m,n} = 
\{ X \mid 
\text{
$\phi^{ \{ 1 \} }(X)
\in \mcC_m$ and 
$\phi^{\mathbb{Z}/2}(X) \in \mcC_n$
} \},
\end{equation*}
where $m,n \in [0, \infty].$

\end{itemize}

\end{enumerate}

\end{Theorem} 

Jjust like 
the nilpotency theorem 
Theorem~\ref{Hopkins-Smith nilpotency} was 
crucial in the proof of Hopkins-Smith 
theorem Theorem~\ref{Hopkins-Smith}, 
the above theorem of Strickland is shown 
by first proving an appropriate nilpotency 
theorem \cite[Th.3]{Joachimi}.  
At the same time, Joachimi \cite{Joachimi} 
explains various difficulties in 
proving an appropriate nilpotency theorem 
in the motivic setting.
Furthermore, the above Joachmi's 
 construction of motivic thick ideals in 
$\mcS \mcH (k)^c$ for $k \subseteq 
\mathbb{C}$ is so far limited to 
importing the Hopkins-Smith stable homotopy thick ideals in $\mcS \mcH^c.$ 
Thus, constructions of motivic thick ideals 
of truly algebro-geometric origin is 
highly desired.  For details and much more 
of Joachimi's work, construct her own 
exposition \cite{Joachimi} in this 
proceeding.

For a case of $k \not\subseteq \mathbb{C},$
Kelly \cite{Kelly}
obtained the following 
surprisingly  simple 
description of the set of prime thick tensor 
ideals $\operatorname{Spc} \left(\mcS \mcH (
\mathbb{F}_q)^c_{\mathbb{Q}} \right)$ 
\footnote{
See Defjnition~\ref{tt-geometry} for this
concept.
}, 
up to a couple of widely believed conjectures:

\begin{Theorem} {\rm \cite[Th.1.1]{Kelly} }
Let $\mathbb{F}_q$ be a field with a prime power, $q,$ number of
elements. Suppose that for all connected smooth projective varieties $X$ we have:
\begin{equation*}
\begin{split}
CH^i(X; j)_{\mathbb{Q}} = 0; \quad
\forall j \neq 0; i \in \mathbb{Z}\quad 
 \text{(Beilinson-Parshin conjecture)},
\\
CH^i(X)_{\mathbb{Q}} \otimes
CH_i(X)_{\mathbb{Q}} \to 
CH_0(X)_{\mathbb{Q}}
\ \text{is non-degenerate.}\quad
\text{(Rat. and num. equiv. agree)}
\end{split}
\end{equation*}
Then
\begin{equation*}
\operatorname{Spc}\left(
\mcS \mcH (\mathbb{F}_q)^c_{\mathbb{Q}}
\right) \cong 
\operatorname{Spec}( \mathbb{Q} ).
\end{equation*}

\end{Theorem}

For details, consult 
Kelly's own 
exposition \cite{Kelly} in this 
proceeding.

}

\section{
$\Dcoh(X)$ and $\Dperf(X)$
}
%
%
%
%
%

In the last two sections, we reviewed:
\begin{itemize}
\item 
Ohkawa's theorem in $\mcS \mcH$, 
which states the Bousfield classes form a 
somewhat mysterious set.
\item
Its analogue in $\Dqc(X)$ is explicitly computable: the fundamental theorem of Hopkins,
Neeman,..., identifies the set of 
Bousfield classes with the set of of 
 localizing tensor ideals, which turns out 
 to have a concrete and algebro-geometric description.
%
%

\item  
Hopkins' motivation of  
his fundamental theorem 
in $\Dqc (X)$ was his own theorem with 
Smith in $\mcS \mcH^c.$

\end{itemize}

Thus, we are naturally led to 
investigate 
$\Dqc (X)^c.$  However, the story  
is not so simple.
Whereas there is a conceptually simple 
categorical interpretation 
$\Dqc (X)^c = \Dperf (X),$
it is 
its close relative (actually equivalent if 
$X$ is smooth over a field) 
$\Dcoh (X)$ which traditionally has been 
 intensively studied
because of its rich geometric and physical 
information. 
\footnote{
Or, researchers might prefer 
\lq\lq $\heartsuit$-felt\rq\rq \ 
$\Dcoh (X) \cong \mcD^b( \Coh (X) )$
(although separated, not mere 
quasi-separated, assumption is needed
for this equivalence)
over simply formal 
$\Dperf(X) \cong \Dqc (X)^c$...
}

}
So, we wish to understand both 
$\Dcoh (X)$ and $\Dperf (X).$

In this section, 
we start with brief,
 and so inevitably 
incomplete, 
summaries of  
$\Dcoh (X)$ and $\Dperf (X),$ 
focusing on their usages.
Still, we hope this would convince 
non-experts that $\Dcoh (X)$ and $\Dperf (X)$ 
are very important objects to study.

Then, we shall review Neeman's recent 
result, which claims these two close relatives
$\Dcoh (X)$ and $\Dperf (X)$ actually 
determine each other, and its 
main technical tool: approximable 
triangulated category.

\subsection{$\Dcoh (X)$}

\begin{itemize}

\item There is a classical functoriality result of Grothendieck:

\begin{Theorem} 
\label{proper morphism preserves coherence}
{\rm \cite[Th.3.2.1]{Gro61}}
Let $f : X\to Y$ be a proper 
morphism with $Y$ locally noetherian.
Then
\begin{equation*}
\Rd f_* \Dcoh (X) \ \subset \ \Dcoh (Y).
\end{equation*}

\end{Theorem}

Actually, there is a sharp converse 
(i.e. we do not have to check 
$\Rd f_* \Dcoh (X) \ \subset \ 
\Dcoh (Y)$)
to Theorem~\ref{proper morphism preserves coherence} \cite[Cor.4.3.2]{LN07} 
\cite[Lem.0.20]{1703.04484}:

\begin{Theorem}
\label{preserving coherence implies proper}
{\rm \cite[Lem.0.20]{1703.04484}}
Let $f : X \to Y$ be a separated, 
finite-type morphism of noetherian schemes 
such that
\begin{equation*}
\Rd f_* \Dperf (X) \ \subset \ 
\Dcoh (Y).
\end{equation*}
Then $f$ is proper.
\end{Theorem}

%

\item
For an essentially small triangulated category $\mcT,$ 
its Grothedieck $K_0$-group $K_0(\mcT)$ 
is defined by generators and relations as 
 follows
\cite[Def.4.5.8]{MR1812507}
\cite[Def.1]{Nee05}:
\begin{equation} \label{K_0-definition}
\hspace{1mm}
K_0(\mcT) :=
\frac{
\mathbb{Z}
\left\{
[X]  \mid 
\text{
$[X]$ is an isomorphism class of 
$X\in \mcT$}
\right\}
}{\mathbb{Z}
\left\{
[X] - [Y] + [Z] \mid 
\text{there is a distiguished triangle}\ 
X \to Y \to Z \to \Sigma X  
\right\}
}
\end{equation}

\begin{itemize}
\item Having defined $K_0(\mcT),$ 
we should not be too optimistic to hope  
 $K_0(\mcT)$ always carries a 
rich information of $\mcT,$   In fact, if $\mcT$ 
contains an arbitrary countable direct sum 
(coproduct)
\footnote{
Having arbitrary small coproducts was 
an indispensqble 
assumption for Brown representability and 
Bousfield localization
(Theorem~\ref{Brown representability},
Corollary~\ref{Existence of Bousfield localization}).
}
, then, for any $X\in \mcT,$ 
we have a distinguished triangle of the 
following form:
\begin{equation*}
\oplus_{n\in \mathbb{N}} X
\ \xrightarrow{\text{index shift}} \
\oplus_{n\in \mathbb{N}} X \ \to \ X \ 
\to \ 
\Sigma 
\left(
\oplus_{n\in \mathbb{N}} X
\right)
\end{equation*}
From the defining relation of $K_0(\mcT)$ 
\eqref{K_0-definition}, this implies
$[X] = 0 \in K_0( \mcT)$ for any  
$X\in \mcT .$ By the definition 
\eqref{K_0-definition}, this means
$K_0( \mcT ) = 0$ whenever  $\mcT$ 
contains an arbitrary countable direct sum 
(coproduct).  As a very important special 
case, we emphasize:
\begin{equation*}
K_0( \Dqc (X) ) = 0.
\end{equation*}

\item Grothendieck $K_0$-group is useful
to classify dense subcategories of an
essentially small triangulated subcategory.

\begin{Proposition}
\label{dense triangulated category}
{\rm \cite[p.5,Lem.2.2,p.6,Cor.2.3]{MR1436741}
\cite[Prop.4.5.11]{MR1812507}
}
Suppose a triangulated 
subcategory $\mcS$ of an essentially 
small triangulated 
category $\mcT$ is 
\underline{\bf dense}, i.e. 
$\widehat{\mcS} = \mcT.$ Then,

\begin{enumerate}
\item The induced map
$K_0(\mcS ) \to K_0(\mcT)$ is a 
monomorphism.
\item For any $X \in \mcT,$ 
\begin{equation*}
X \in \mcS \quad \iff \quad 
[X] \in \operatorname{Im}
( K_0(\mcS ) \to K_0(\mcT) ).
\end{equation*}

\end{enumerate}

\end{Proposition}

\begin{Theorem} 
{\rm \cite[p.5,Th.2.1]{MR1436741} }
For an essentially small triangulated 
category $\mcT,$ there is a one-to-one correspondence between the
dense triangulated subcategories 
of $\mcT$ and 
the subgroups of $K_0(\mcT)$:
\begin{equation*}
\begin{split}
\left\{
\text{dense triangulated subcategories
of $\mcT$}
\right\}
\
&\overset{\cong}
{\rightleftarrows}
\ 
\left\{
\text{
subgroups of $K_0(\mcT)$
}
\right\}
\\
\mcS \ 
&|\hspace{-1.2mm}\rightarrow
\
\operatorname{Im}
\left( 
K_0( \mcS ) \to K_0( \mcT )
\right)
\\
{
\text{
$\triangle$ subcategory 
consisting of 
$X \in \mcT$ 
with
$[X] \in H \subseteq K_0(\mcT)$}
}
\
&\leftarrow\hspace{-1.2mm}|
\
H
\end{split}
\end{equation*}

\end{Theorem}

\item For any small abelian category 
$\mcA,$ the 
  functor $\mcD^b$ comes with 
the canonical
embedding $\mcA \to \mcD^b(\mcA),$ 
which
induces an equivalence of Grothendieck 
$K$-groups of an abelian category $\mcT$ 
and a triangulated category 
$\mcD^b(\mcA)$:
\begin{equation} \label{K-isomorphism}
K_0\left( \mcA \right) \xrightarrow{\cong}
K_0 \left( \mcD^b(\mcA) \right),
\end{equation}

\item Whenever a bounded $t$-structure 
is given on $\mcT,$ if we denote by 
$\mcT^{\heartsuit}$ its heart, then 
we have another isomorphism of 
$K_0$-groups of an abelian category and 
a triangulated category:
\begin{equation}
\label{another K-isomorphism}
K_0( \mcT^{\heartsuit} ) \ 
\xrightarrow{\cong} \
K_0( \mcT ).
\end{equation}
Applying 
\eqref{another K-isomorphism} 
 to $\mcT = \Dcoh (X), 
\mcT^{\heartsuit}= \Coh (X),$ 
\footnote{
If we apply \eqref{K-isomorphism} 
in order to obtain the isomorphism 
\eqref{K_0(X)-isomorphism}, we must
require the extra \lq\lq separated \rq\rq \ 
assumption, for then we should also use
the isomoprhism:
\begin{equation} \label{Dbcoh(X)}
\Dcoh (X) = 
\mcD^b( \Coh (X) ),
\end{equation}
which requires the \lq\lq separated\rq\rq \
assumption of $X.$ This fact, and the
above approach to use 
\eqref{another K-isomorphism} was 
communicated to the author by Professor 
Neeman.
%
%
%
%
}
we find the canonical 
isomorphism:

\begin{equation}
\label{K_0(X)-isomorphism}
K\left( \Coh (X) \right) \xrightarrow{\cong}
K \left( \Dcoh (X) \right)
\end{equation}

\end{itemize}

\item The sheaf theory has its origin
in Oka-Cartan theory of complex functions
of several variables (see e.g. 
\cite{MR3443603} for a general picture, 
and \cite{OhsPr} for a review of  
the $L^2$-technique in complex geometry,
both by Ohsawa
\footnote{
Professor Takeo Ohsawa is the 
AMS Stefan Bergman Prize 2014 recipient.
His survey paper \cite{OhsPr} in this 
proceedings is a concise summary of 
his work for which this prize was awarded. 
It was his Bergman Prize money which enabled 
us to invite distinguished lecturers to 
 Ohkawa's memorial conference at Nagoya 
University in the 
summer of 2015.
Takeo Ohsawa was also 
Tetsuske Ohkawa's highschool classmate at 
Kanazawa University High School in 
Kanazawa, Japan.
}
).
The pivotal achievement at the time was
\emph{Oka's Coherence Theorem}, which 
states that the structure sheaf 
$\mcO_M$ of a complex manifold $M$ is
coherent (for a proof, see e.g. 
\cite{MR3526579}). From the viewpoint of 
algebraic geometry, interest of complex
manifolds emerge through the GAGA theorem
of Serre \cite{MR0082175}, which, 
for a proper scheme $X$ over $\Spec 
\mathbb{C},$ can
be stated as an equivalence
of abelian categories of
coherent modules 
\cite[XII,Th.4.4]{SGA1}:
\begin{equation*}
\phi^* : \Coh (X) \xrightarrow{\cong}
\Coh (X^{\text{an}}),
\end{equation*}
where
$
\phi : X^{\text{an}} \to X
$
is the canonical morphism from the
associated analytic space 
$X^{\text{an}}$ of $X$ 
\cite[XII,1.1]{SGA1}, and 
$\phi^*$ 
consequently induces isomorphisms of 
resulting derived categories:
\footnote{
$X$ being proper over $\Spec(\mathbb{C})$ 
 implies (as part of
the definition of properness) that it is separated, hence 
$\DD^b (\Coh (X)) = \Dcoh (X).$
Hence, these two isomorphisms are trivial
consequences of the isomorphism 
$\phi^* : \Coh (X) \xrightarrow{\cong}
\Coh (X^{\text{an}}).$
These two isomorphism are supplied just 
for reader's information.
}

\begin{equation*}
\DD_{\Coh} (X ) \xrightarrow{\cong} 
\DD_{\Coh} (X^{\text{an}} );\quad 
\Dcoh (X ) \xrightarrow{\cong} 
\Dcoh (X^{\text{an}} );\ \cdots
\end{equation*}

Recently, Jack Hall \cite{1804.01976}
proposed a unified treatment of
\lq\lq GAGA type theorems,\rq\rq 
in which, a prominent role of Oka's coherence theorem became transparent 
in his deduction
of the classical GAGA theorem 
\cite[Example 7.5]{1804.01976} 
(also consult the updated version of 
\cite[Remark 1.7 and Appendix A]{1804.02240} to 
appreciate how short and simple the Jack 
 Hall's new proof is.).

\item Derived categories in the complex 
analytic setting shows up in the 
Kontsevich homological
mirror symmetry \cite{MR1403918}
\footnote{
Of course, there are many other mathematical approaches to physics.
For instance, some of Costello's approach 
to quantum field theory via 
Lurie's higher algebra \cite{Lur09,Lur16} 
point of vew are touched upon in
Matsuoka's surveys \cite{M1Pr,M2Pr} in this 
proceedings.
}
which in the Calabi-Yau setting
is of the following form:
%
\begin{equation} \label{Kontsevich homological mirror conjecture}
\Dcoh (X) \cong \DD^b Fuk (X^{\vee} ),
\end{equation}
where $X$ is expected to be
a mirror of  $X^{\vee},$  
given by a sigma model: 
\begin{equation*} \label{sigma model}
(M, I, \omega, B),
\end{equation*}
where we only note $I$ is the complex
structure of $M,$ and that
whose category of $D$-branes of type $B$ 
(B-model) is the left side of 
\eqref{Kontsevich homological mirror conjecture} : 
\begin{equation*}
DB (M, I, \omega, B) \ \cong \ 
\Dcoh (M, I)
\ \cong \
\Dcoh (X).
\end{equation*}

On the other hand,
$\DD^b Fuk (X^{\vee} ),$  the 
\emph{derived Fukaya category} consisting
of Langrangian submanifolds of the
mirror $X^{\vee},$ is not a derived
category of an abelian category 
(but of an $A_{\infty}$ category;  
see \cite{FOOO1,FOOO2} for more details).

\item 
Recall that 
$\Dcoh (X)$ is given by the composite of functors:
\begin{equation} \label{reconstruction}
\Dcoh : 
X \overset{\Coh}\mapsto \Coh (X) 
\overset{\mcD^b}\mapsto 
\mcD^b( \Coh (X) ) = \Dcoh (X).
\end{equation}

It is instructive to keep reconstruction 
problems arising from these functors  
in mind.  
For instance, 
Theorem~\ref{Gabriel-Rosenberg} 
of Gabriel-Rosenberg 
can be specialized to the following
(which is essentially the original 
theorem of Gabriel \cite{Gab62})
reconstruction theorem with respect to
$\Coh$:
\footnote{
Theorem~\ref{Gabriel} is reduced to 
Theorem~\ref{Gabriel-Rosenberg}
for $\QCoh (X) \cong 
\operatorname{Ind} \Coh (X)$ under the 
Noetherian hypothesis \cite[Lem.3.9]{Lur04}. See also \cite[p.2]{CG15} \cite{Per09}.
}

\begin{Theorem} \label{Gabriel}
Any Noetherian and separated scheme $X$ 
can be reconstructed from $\Coh (X).$

\end{Theorem}

\item Glancing at this theorem of Gabriel,
we naturally hope $\Dcoh (X)$ would 
carries rich information of $X.$ 
Concerning the reconstruction problem 
associated with \eqref{reconstruction},
any smooth connected projective 
variety with either $K_X$ ample
or $-K_X$ ample can be 
reconstructed from $\Dcoh (X)$
(the Bondal-Orlov  reconstruction 
theorem \cite{MR1818984} ).

\item On the other hand, among those $X$ with
trivial $K_X$ like an abelian variety 
or Calabi-Yau, 
 many examples of 
so-called Fourier-Mukai partners,
i.e. non-isomorphic smooth projective
varieties with equivalent 
$\Dcoh,$ have been produced, 
 starting with Mukai \cite{MR0607081},
Thus, the restruction for the composite
$\Dcoh: X  \mapsto \Dcoh (X)$ in
\eqref{reconstruction} does not
hold in general.
Considering the Gabriel reconstruction 
Theorem~\ref{Gabriel}, we find this 
failure results from that 
of the reconstruction of $\mcD^b$ among 
those $X$ with trivial $K_X.$
This suggests an existence of a of 
 moduli of hearts of $\Dcoh (X)$ for these $X.$ 
\footnote{
As we shall briefly review later, 
Bridgeland's space of stability conditions 
is a kind of moduli space of 
\lq\lq enriched hearts\rq\rq \ of 
a triangulated category.
}

\item 

If $X$ is affine locally regular
and finite-dimensional, then we have 
the following canonical equivalence
(which is a local assertion):
\begin{equation*}
\Dcoh (X) \xrightarrow{\simeq} \Dperf (X)
\end{equation*}

This, in turn, suggests the 
Verdier quotient
\begin{equation*}
\DD_{Sg}(X) := \Dcoh (X) / \Dperf (X)
\end{equation*}
reflects singular information of $X,$ and 
is consequently called the derived category
of singularities \cite[Def.1.8]{MR2101296}.

In the Kontsevich homological mirror
symmetry, a mirror of varieties other than
Calabi-Yau is not expected to be given
by a sigma model.  
For a variety with either $K_X$ ample
or $-K_X$ ample, its mirror is
expected to be given by a
Landau-Ginzburg model
\begin{equation*} \label{Landau-Ginzburg model}
(Y, I, \omega, B, W),
\end{equation*}
where $W : Y \to \mcA^1$ is a regular
function called the superpotential. 
In this case, the category of $D$-branes
of type $B$ is, via its identification with the category of 
matrix factorizations, shown to be
of the following form 
\cite{MR2039036}\cite{MR2101296}\cite{MR2910782} :
\begin{equation} \label{oracle of physics}
DB(Y, I, \omega, B, W) \ \cong \
\prod_{\lambda \in \mcA^1}
\DD_{Sg}\left( W^{-1}(\lambda) \right).
\end{equation}

\item This oracle of physics 
\eqref{oracle of physics}, which 
highlights essentially only the 
singular part, might appear surprising 
for mathematicians.  However, 
in the development of the minimal model 
program in birational geometry, it has 
become clear that we should take into 
account singular information even if
we are only interested in smooth ones 
\cite{MP97}\cite{KM98}\cite{Mat02}.

Now, close relationship between 
$\Dcoh$ and 
birational geometry
have been observed 
\cite{MR1957019, MR1949787}. A central
problem here is the 
\emph{Kawamata DK-hypothesis}:

\begin{Conjecture}
{\rm
\cite[Conj.1.2]
{1710.07370}
}
For birationally equivalent smooth
projective varieties $X,Y,$ suppose there
 exists a smooth projective
variety $Z$ with birational morphisms
$f : Z \to X,\ g : Z \to Y.$

\begin{description}
\item[
\underline{
$K$-equivalence $\implies$
$D$-equivalence}:
]
\begin{equation*}
\begin{split}
&\text{ $K$-equivalence  } \quad
\bigg(i.e.\quad 
f^*K_X \sim g^*K_Y  \quad
(\text{linearly equivalent})
\bigg)
\\
\text{{\rm implies}}\qquad
&\text{ $D$-equivalence } \quad
\bigg(i.e.\quad
\Dcoh (X) \cong \Dcoh (Y) \bigg)
\end{split}
\end{equation*}

\item[
\underline{
$K$-inequality $\implies$ 
fully faithful triangulated funtor}:
]

\begin{equation*}
\begin{split}
&\text{ $K$-inequality } \quad
\bigg(
\substack{
i.e.\quad 
\text{there exists an effective divisor
$E$ on $Z$ s.t.}
\\
f^*K_X + E\sim g^*K_Y  \quad
(\text{linearly equivalent})
}
\bigg)
\\
\text{{\rm implies}}\qquad
&
\bigg(
\substack{
 \text{there is a fully faithful functor of triangulated categories}
\\
\Dcoh (X) \to \Dcoh (Y).
}
\bigg)
\end{split}
\end{equation*}
%
%

\end{description}

\end{Conjecture}

While the converse
($D$-equivalence $\implies$
$K$-equivalence) does not hold
in general
\cite{MR2067481}, 
if there is a fully faithful functor 
$\Psi : \Dcoh (X) \to \Dcoh (Y),$
then we obtain a semi-orthogonal
decomposition \eqref{semiorthogona decomposition} 
\cite{MR1996800}:
\begin{equation}
\Dcoh (Y) = \langle
\Psi \left( \Dcoh (X) \right)^{\perp},
 \Psi \left( \Dcoh (X) \right) \rangle
\end{equation}

\item Motivated by the Kontsevich
homologial mirror symmetry, 
some previously unexpected structures 
of $\Dcoh (X)$ have been discovered:

\begin{itemize}
\item Motivated by the generalized
Dehn twist associated with the
Lagrangian spheres of the
(hypothetical) mirror $X^{\vee},$ 
Seidel-Thomas \cite{MR1831820}
constructed a braid group $B_{m+1}$ action 
under the presence of  the spherical
$A_m$-configuration, i.e.
there are $\mcE_i\in \Dcoh (X) \ 
(1\leq i\leq m)$ such that the
following two conditions are satisfied:
\ \vspace{3mm}

\begin{description}
\item[(sphericality):]
For $1\leq i \leq m,\ 
\mcE_i\otimes \omega_X \cong
\mcE_i$ and 
\begin{equation*}
\Hom_{\Dcoh (X)}
\left( \mcE_i, \mcE_i [r] \right)
=
\begin{cases}
\mathbb{C} \ &\text{if} \ r = 0, \dim X
\\
0 \ &\text{if} \ r \neq 0, \dim X
\end{cases}
\end{equation*}

\item[($A_m$-configuration):]

\begin{equation*}
\operatorname{dim}_{\mathbb{C}} 
\oplus_r \Hom_{\Dcoh (X)}
\left( \mcE_i, \mcE_j [r] \right)
=
\begin{cases}
1 \quad & | i - j | = 1 
\\
0 \quad & | i - j | \geq 2.
\end{cases}
\end{equation*}

\end{description}


\item
Going back to the reconstruction problem 
of $\mcD^b$ in \eqref{reconstruction}, 
existence of Fourier-Mukai partners 
suggests an existence of a 
 moduli of hearts of $\Dcoh (X)
= \mcD^b( \Coh (X) ).$
%
%
%
%
%

To begin with, we recall a related toy 
model for 
$\Coh (X),$ where we can 
construct moduli spaces, 
$M_{\mcO_X(1)}(P)$ for a fixed 
Hilbert polynomial, by restricting to 
(Gieseker-Maruyama-Simpson)
(semi)-stable sheaves 
\cite[Th.4.3.4]{MR2665168}.

Thus, its not suprising that 
some kind of stability 
condition is needed to construct 
a moduli in of hearts of $\Dcoh (X)
= \mcD^b( \Coh (X) ).$
In fact, axiomatizing Douglas' study 
\cite{MR1957548} of the $\Pi$-stability of
D-branes, Bridgeland \cite{MR2373143} 
proposed a way of constructing a moduli space of \lq\lq enriched hearts,\rq\rq \ 
space of stability conditions, 
out of certain triangulated categories.
Bridgeland \cite{MR2373143} 
defined 
a \emph{stability condition} on 
a triangulated category $\mcD$ to be
a data $(Z,\mcA)$ such that:
\begin{itemize}
\item $\mcA \subset \mcD$ is the 
heart of a bounded t-structure on $\mcD.$

\item $Z : K(\mcA) \to \mathbb{C}$ is a 
\emph{stability function}, i.e.
\begin{itemize}
\item $Z : K(\mcA) \to \mathbb{C}$ is a
group homomorphism.
\item For any 
$E \in \mcA \setminus \{ 0 \},$ \ 
\begin{equation*}
\begin{split}
Z(E)  &:= 
 r(E) \exp \left(i\pi \phi (E) \right)
 \quad ( r(E) > 0, 0 < \phi(E) \leq 1 ) 
\\
&\in 
\overline{\mathbb{H}} := \ 
\left\{ r \exp (i\pi \phi ) \mid
r > 0, 0 < \phi \leq 1 \right\}.
\end{split}
\end{equation*}
\end{itemize}

\item This stability function 
$Z : K(\mcA) \to \mathbb{C}$ is 
furthermore a 
\emph{stability condition}, i.e.
any $E \in \mcA$ admits a 
\emph{Harder-Narasimhan filtration}:
\begin{equation*}
0 = E_0 \subset E_1 \subset \cdots
\subset E_n = E,
\end{equation*}
such that
\begin{itemize}
\item each $F_i = E_i / E_{i-1}$ is
\emph{$Z$-semistable}, i.e. 
for all nonzero subobjects $F_i' \subset 
F_i$ we have 
$$\phi (F_i') \leq \phi (F_i).$$
%

\item
$\phi (F_1 ) > \phi (F_2) > \cdots > 
\phi (F_n ).$

\end{itemize}

Since $Z$ is a homomorphism, we can easily
verify:
\begin{equation*}
\text{$E,F$:\ $Z$-semistable s.t.
$\phi (E) > \phi (F)$
} \ \implies \
\operatorname{Hom}_{\mcA}(E,F)=0.
\end{equation*}

Thus, topologists should recognize
a similarity between the Harder-Narasimhan 
 filtration and the (finite) Postnikov
tower with the following analogy
\begin{equation*}
\begin{split}
&\text{$K(\pi_1,n_1), K(\pi_2,n_2)$:\ 
Eilenberg-MacLane spectra s.t.
$n_1 >  n_2$
} \\
 \implies \
&\operatorname{Hom}_{\Ho}
\left( K(\pi_1,n_1), K(\pi_2,n_2) \right)
= H^{n_2}\left( K(\pi_1,n_1), \pi_2 \right)
= 0.
\end{split}
\end{equation*}

\end{itemize}

Here, we wish to vary the heart 
$\mcA=\mcD^{\heartsuit}$ 
while fixing the amibient
triangulated category $\mcD.$
For this purpose, in view of 
\eqref{K-isomorphism}, 
we impose an extra structure on the
stability function, i.e.
\begin{equation*}
\hspace{-4mm}
\xymatrix{
 K\left( 
\mcD^{\heartsuit} \right)
\ar[r]^{\cong} &
K\left( \mcD  \right)
\ar[r]^-Z \ar[d]_{\operatorname{cl}} 
& \mathbb{C}
\\
& \Gamma \ar@{-->}[ur]_{\exists} &
}, 
\text{where} \ 
\begin{cases}
\text{
$\Gamma$ is a finitely generated 
free abelian group,
}
\\
\text{\qquad s.t. 
$\Gamma\otimes_{\mathbb{Z}}\mathbb{R}$
is equipped with a norm}
\\
\text{\qquad 
(which allows us to define 
$\Vert \operatorname{cl}(E) \Vert$ 
for $E\in K(\mcD)$
).
}
\\
\operatorname{cl} : \Gamma \to
\mathbb{C} \
\text{is a homomorphism}
\\
\end{cases}
\end{equation*}
We further impose the 
\emph{support property} \cite{arXiv:0811.2435}:
\begin{equation*}
\left\{ \frac{ | Z(E) | }{ 
\Vert  \operatorname{cl}(E) \Vert } \ 
\bigg| \ 
E \in 
\left( \cup_{i\in \mathbb{Z} } 
\mcD^{\heartsuit} [i] \right) \setminus 0 \right\}
\ \text{is bounded.}
\end{equation*}

When we fix $\mcD$ with such a
homomorphism $K(\mcD) \to \Gamma,$
Bridgeland \cite{MR2373143} showed the
set of such stability conditions
can be topologized and 
becomes a complex manifold 
$\operatorname{Stab}_{\Gamma}(\mcD).$

However, for the case of our interest
$\mcD = \Dcoh (X),$ as soon as
$\dim X\geq 3,$ there is no
stability condition on $\mcD=\Dcoh (X)$
with $\mcD^{\heartsuit} =
\Coh (X)$ 
\cite[Lem.2.7]{MR2541209}, and 
even 
the existence of such a stability 
condition is problematic, i.e. 
the possibility of 
$\operatorname{Stab}_{\Gamma}(\mcD)
= \emptyset$ is yet to be excluded.

\end{itemize}

\end{itemize}

\subsection{$\Dperf (X)$}

\begin{itemize}

\item The functoriality results for 
$\Dcoh$ reviewed in 
Theorem~\ref{proper morphism preserves coherence}
and 
Theorem~\ref{preserving coherence implies proper}
have the following analogue for $\Dperf$:

\begin{Theorem}
{\rm \cite[Th.1.2]{LN07}
\cite[Ill.0.19]{1703.04484}}
For a separated, finite-type morphism of 
 noetherian schemes $f : X \to Y,$
\begin{equation*}
\begin{split}
&\Rd f_* \Dperf (X) \ \subset \ 
\Dperf (Y) \ (\text{i.e. 
\underline{{perfect}} } )
\\
\iff \quad 
&\text{$f$ is proper and 
of finite Tor-dimension}
\end{split}
\end{equation*}

\end{Theorem}

\item
$\Dperf (X)$
can be  directly recovered from
$\Dqc(X):$  

\begin{Theorem} 
{\rm ( \cite{MR1308405}
 \cite{MR1996800} ) }
\label{compact objects}
The canonical functor
\begin{equation*}
\Dperf (X) \to \Dqc (X)
\end{equation*}
identifies 
$\Dperf (X)$ as the full triangulated
subcategory $\Dqc (X)^c$ of compact objects 
 in $\Dqc (X):$
\begin{equation*}
\Dperf (X) = \Dqc (X)^c
\end{equation*}
\end{Theorem}

\item 
%
Thomason-Trobaugh \cite[App.F]{MR1106918} 
proved $\Dperf (X) = \Dqc (X)^c$ is 
essentially small (i.e. equivalent to a 
small category) for 
any quasi-compact and quasiseparated
scheme $X$ (e.g. for any noetherian scheme). 
Starting with this, Thomason 
\cite[Th.3.15]{MR1436741} classified 
thick tensor triangulated 
ideals of  $\Dperf (X) = \Dqc (X)^c$ 
for any quasi-compact and quasiseparated
scheme $X.$
Here, we review Paul Balmer's generalization \cite{MR2196732} 
of such a classification to 
certain essentially small 
tensor triangulated categories.


\begin{Definition} \label{tt-spectrum}
For 
a tensor triangulated category $\mcK,$ 
\begin{itemize}

\item
{\rm
\cite[Def.4.1]{MR2196732} 
\cite[Def.7]{Bal10}
}
A thick tensor ideal $\mcI \subset 
\mcK$ is called 
\underline{\em radical} if
\begin{equation*}
\mcI = \sqrt{\mcI} := 
\left\{ a\in \mcK \ \mid \ 
\exists n\geq 1 \ 
\text{such that}\ a^{\otimes n}
\in \mcI \right\}.
\end{equation*}
The collection of 
radical thick 
tensor ideals of $\mcK$ is denoted by 
$\mathbb{R}(\mcK).$

\item
{\rm
\cite[Def.2.1]{MR2196732} 
\cite[Con.8]{Bal10},\ 
(see also 
Definition~\ref{tensor ideal, prime, strongly dualizable})
}\ 
A proper thick tensor ideal 
$\mcP \subsetneq \mcK$ is called 
\underline{\em prime}, if
\begin{equation*}
a\otimes b \in \mcP \quad \implies \quad 
a\in \mcP \ \text{or} \ b\in \mcP.
\end{equation*}
\item 
{\rm
\cite[Def.2.1]{MR2196732} 
\cite[Con.8]{Bal10}
}\ 
If $\mcK$ is further 
essentially small, 
its
\underline{
\emph{spectrum}
} $\Spc (\mcK)$ 
is given by the following 
(set, by the \lq\lq essentially small\rq\rq \
assumption):
\begin{equation*}
\Spc (\mcK) = 
\left\{ \mcP \subsetneq \mcK \ \mid \ 
\text{$\mcP$ is a proper prime
thick tensor ideal 
 of $\mcK$} \right\},
\end{equation*} 
which is endowed with the topology 
whose open subsets are of the form
\begin{equation*}
U(\mcE) := \left\{ \mcP \in 
\Spc (\mcK)
\ \mid \
\mcE \cap \mcP \neq \emptyset 
\right\} \ \quad
(\mcE \subseteq \mcK) ;
\end{equation*}
in other words, given by the closed basis
$\{ \supp (a) \}_{a \in \mcK},$ where
\begin{equation*}
\supp (a) = \{ \mcP \in 
\Spc (\mcK) \ \mid \ a \not\in \mcP \}
\end{equation*}
is the \underline{\em support} of 
$a\in \mcK.$
\footnote{
\underline{WARNING!} We had already 
 introduced  
the same notation $\supp$ back in 
Definition~\ref{support definitions}. 
However, from Proposition~\ref{support definitions conincide for coherent complexes}
Theorem~\ref{topological reconstruction}, 
these two usages
of $\supp$ coincide 
for the most fundamental example of 
$\mcK = \Dperf (X).$
}
\item
{\rm
\cite[Rem.12]{Bal10}
}
For a general topological space $T$
(we are particularly interested in the 
case $T=\Spc (\mcK)$), 
a subset $Y\subset T
$ 
of the form
\begin{equation*}
Y = \cup_{i\in I}Y_i \quad
\text{with each complement 
$X\setminus Y_i$ open and 
quasi-compact}
\end{equation*}
is called a 
\underline{\em Thomason subset}\ 
of $T
.$ 
The set of Thomason subsets of $T$ 
is denoted by 
$\operatorname{Tho}(T).$

\end{itemize}

\end{Definition}

\begin{Theorem} 
\label{Tho=R,T}
{\rm (i) \cite[Th.4.10]{MR2196732} 
\cite[Th.14]{Bal10} 
\cite[Th.5.9]{MR2806103}
 } 
For an essentially small tensor 
triangulated category $\mcK,$ 
there are mutually inverse 
isomorphisms between radical 
thick tensor ideals of $\mcK$ and 
Thomason subsets of $\Spc (\mcK)$:
\begin{equation}
\begin{split}
\mcK_{-} : 
\operatorname{Tho}\left( \Spc (\mcK) \right)
\
&\overset{\cong}
{\rightleftarrows}
\ 
\mathbb{R}(\mcK) : \supp
\\
Y \ 
&|\hspace{-1.2mm}\rightarrow
\
\mcK_Y := \{ a \in \mcK \mid 
\supp (a) \subset Y \}
\\
\supp ( \mcR ) := 
\cup_{a\in \mcR} \supp (a)
\
&\leftarrow\hspace{-1.2mm}|
\
\mcR
\end{split}
\end{equation}
\newline
{\rm (ii) \cite[Prop.2.4]{Bal07}}
Suppose further $\mcK$ is rigid, then 
every thick tensor ideal is radical, 
and so,
$\mathbb{R}(\mcK) = \mathbb{T}(\mcK).$
Consequently, the mutually inverse 
isomorphisms in (i) becomes the following:
\begin{equation*}
\mcK_{-} : \operatorname{Tho}\left(\Spc (\mcK) \right) \rightleftarrows 
\mathbb{T}(\mcK) : \supp
\end{equation*}

\end{Theorem}

\begin{Theorem}
\label{topological reconstruction}
{\rm
\cite{MR1436741}
\cite[Cor.5.6]{MR2196732} 
\cite[Cor.5.2]{BKS07}
\cite[Th.16]{Bal10}
}
For a quasi-compact and 
quasi-separated scheme $X,$ its underlying
topological space $| X |$ is homeomorphic 
to the spectrum $\Spc ( \Dperf (X) )$ via
\begin{equation*}
\begin{split}
| X | &\xrightarrow{\cong} 
\Spc \left( \Dperf (X) \right)
\\
x &\mapsto \mathfrak{P}(x) :=
\left\{ P \in \Dperf (X) \ \mid \
P_x \cong 0 \right\}.
\end{split}
\end{equation*}
For any $P \in \Dperf (X),$ this 
homeomorphism restricts to the homeomorphism
\begin{equation*}
\Supph(P) \xrightarrow{\cong} \supp(P),
\end{equation*}
where $\Supph(P) \subseteq X$ is the 
\underline{\em homological support} of 
$P \in \Dperf(X),$ i.e. 
the usual sheaf theoretical support of 
the total homology of $P$ given in 
Definition~\ref{support definitions}
Proposition~\ref{support definitions conincide for coherent complexes}.

\end{Theorem}

From Theorem~\ref{topological reconstruction}, Theorem~\ref{Tho=R,T}
(ii) yields the following theorem of 
Thomason, which is a $\Dqc (X)$ analogue
of the Hopkins-Smith 
Theorem~\ref{Hopkins-Smith}:

\begin{Theorem} 
\label{Thomason's theorem of Thomason sets}
{\rm \cite[Th.3.15]{MR1436741}}
For a quasi-compact and 
quasi-separated scheme $X,$ 
there are mutually inverse 
isomorphisms between 
thick tensor ideals of $\Dperf(X)$ and 
Thomason subsets of $| X |$:
\begin{equation}
\begin{split}
\Dperf_{-}(X) : 
\operatorname{Tho}\left( | X | \right)
\
&\overset{\cong}
{\rightleftarrows}
\ 
\mathbb{T}\left( \Dperf (X) \right) : 
\supp
\\
Y \ 
&|\hspace{-1.2mm}\rightarrow
\
\Dperf_{Y}(X) := \{ P \in \Dperf (X) \mid 
\Supph (P) \subset Y \}
\\
\supp ( \mcR ) := 
\cup_{a\in \mcR} \supp (a)
\
&\leftarrow\hspace{-1.2mm}|
\
\mcR.
\end{split}
\end{equation}

\end{Theorem}

\begin{Remark}
\label{thick tensor ideal generated by a single element}
{\rm \cite[Lem.3.1]{1703.04484}}

For an object $H$ of a tensor triangulated 
category $\mcT,$ denote by 
$\langle H \rangle_{\otimes}$ the 
thick tensor ideal (tensor)
generated by $H.$ 
Then we easily see: 
\begin{equation*}
\langle H \rangle_{\otimes} = 
\cup_{l \in \mathbb{N}, C\in \mcT} \ 
\langle C\otimes H \rangle_N,
\end{equation*}
where the notation $\langle - \rangle_N$ 
is recalled in Definition~\ref{<>}.

Many tensor triangulated categories $\mcT$  
are (tensor) generated by a single element.

\end{Remark}

It should be mentioned that, just like 
the nilpotency theorem 
Theorem~\ref{Hopkins-Smith nilpotency} was 
crucial in the proof of Hopkins-Smith 
theorem Theorem~\ref{Hopkins-Smith}, 
some algebro-geometric analogue of 
(Devinatz-)Hopkins-Smith nilpotency 
is crucial to prove these algebro-geometric 
analogues of the Hopkins-Smith theorem 
(see e.g. 
\cite[Th.1.1]{MR1174255}
\cite[Th.3.6,Th.3.8]{MR1436741}
).
In this direction, Hovey-Palmieri-Strickland 
\cite[5]{HPS97} developed a general theory 
how nilpotence implies 
classications of thick subcategories.

Now, the following simple consequence of 
the above theorem of Thomason 
will be used later:

\begin{Corollary}
\label{thick tensor ideal with suppph the whole space}
For a quasi-compact and 
quasi-separated scheme $X,$ any thick 
tensor ideal generated by 
a single $H \in \Dperf (X)$ with 
$\Supph(H) = | X |$ is all of 
$\Dperf (X).$

\end{Corollary}

\item
In terms of $\Dperf (X) = \Dqc (X)^c,$ 
we may refine the smashing part of the 
fundamental theorem of Hopkins, Neeman,
Thomason and others 
(Theorem~\ref{Ohkawa for D_{qc}(X)}) to
become an algebraic analogue of the 
Hopkins-Smith theorem 
(Theorem~\ref{Hopkins-Smith}), with
 an extra bonus of the validity of an 
algebraic analogue of the 
telescope conjecture.
We shall review it now, together with 
(a restatement of) 
Theorem~\ref{Ohkawa for D_{qc}(X)}.  
For the notations below, 
consult the list just before 
Theorem~\ref{The diagram for SH}.
%

\begin{Theorem} 
\label{MAIN THEOREM}
{\rm
(\cite{MR0932260}
\cite[Th.2.8,Th.3.3]{MR1174255}
\cite{MR1436741},
\cite[Cor4.6;Cor.4.13;Th.5.6
]{MR2071654}
\cite[Cor.6.8]{MR2806103}
\cite[Cor.6.8;Ex.6.9]{DS13}
\cite[Th.B]{MR3797596})
}
For a Noetherian 
scheme $X,$ we have a commutative 
diagram consisting of mutually inverse
horizontal arrows:
\begin{equation}
\label{THE DIAGRAM}
\xymatrix{
2^{| X | } 
\ar@<1mm>[rr]^-{
\left\{ Q \in \Dqc (X) \  \mid \ 
\supp (Q) \subseteq - \right\}
} & &
\ar@<1mm> [ll]^-{\supp}
\mathbb{L}( \Dqc (X) )
\\
\operatorname{Tho}( | X | )
\ar@{^{(}-_{>}}[u]
\ar@<1mm>[r]^-{\Dperf_{-}(X)} & 
\mathbb{T}\left( \Dperf(X) \right)
\ar@<1mm> [l]^-{\supp}
\ar@<1mm>[r]^-{I_X} & 
\ar@<1mm> [l]^-{C_X}
\mathbb{S}( \Dqc (X) )
\ar@{^{(}-_{>}}[u]
}
\end{equation}
Here,
\begin{itemize}
\item The upper side mutually inverse 
arrows are those in 
Theorem~\ref{Ohkawa for D_{qc}(X)},
which is the analogue of the Ohkawa
theorem and an affirmative solution of the 
Hovey Conjecture~\ref{localizing conjecture}
(ii) for $\Dqc(X).$
\item The lower left side mutually inverse
arrows are those in Thomason's 
Theorem~\ref{Thomason's theorem of Thomason sets}, which is a $\Dqc (X)$ analogue
of the Hopkins-Smith 
Theorem~\ref{Hopkins-Smith}:

\begin{Remark} 
\label{THE REMARK}
The above commutative diagram 
\eqref{THE DIAGRAM} encapsulates our story; 
starting with Ohkawa's theorem in 
$\Ho,$ we then move on to the 
$\Dqc$ analogue, encountering the fundamental theorem
 of Hopkins, Neeman, Thomason and 
others; then going back to $\Ho^c$ to 
appreciate the Hopkins-Smith thick 
category theorem, and then, 
moving back again to the 
$\Dqc^c$ analogue, we 
disvover the above fantastic 
Theorem~\ref{MAIN THEOREM}.

In fact, the commutative diagram 
\eqref{THE DIAGRAM} is a
$\Dqc^c \subset \Dqc$ analogue 
of the commutative diagram 
\eqref{THE DIAGRAM-SH} for 
$\Ho_{(p)}^c \subset \Ho_{(p)}.$ 
%
%
%
Thus the underlying 
message here is to extend the 
commutative diagrams of 
\eqref{THE DIAGRAM} and 
\eqref{THE DIAGRAM-SH} to other 
triangulated categories.
There is a paper of Iyenger-Krause 
\cite{IK12} in this direction, 
and this is exactly the theme of our 
Homework in the intoruction.

\end{Remark}

\item 
The mutually inverse arrows at the bottom right of the diagram yield a positive solution of the 
telescope conjecture 
(Theorem~\ref{telescope conjecture}(iv)) 
by \cite[Cor.6.8]{MR2806103}
\cite[Th.B]{MR3797596}).

\end{itemize}

\end{Theorem}

\item However, the analogue of
\eqref{BL of TC} for  $\Dperf$ does
not hold in general, for
$\Ld j^* : \Dperf(X) \to \Dperf (U)$
is not surjective in general.
Still, as was noticed by Thomason-Trobaugh 
\cite{MR1106918}, there is a similar 
equivalence as soon as we apply the 
thick closure $(-)^{\widehat{}}$:
\footnote{
Let us recall the following related 
result in the setting of abelian 
category of quasi-coherent sheaves, 
which should go back at 
least to Gabiriel (see e.g. 
\cite[Prop.3.1]{Rou10}):
%
$
\Coh (X) \big/ \Coh_Z (X) 
\
\xrightarrow[\cong]{\overline{j^*}} 
\
\Coh (U).
%
$
where the left hand side is the abelian
quotient category in the sense of 
Gabriel, Grothendieck, Serre.
}
\footnote{
The following interesting 
 historical account 
on the difficulty of 
generalizing 
statements in $\Dqc$ 
\eqref{BL of TC}
\eqref{smashing Bousfield localization of Dqc(X)}:
\begin{equation*}
\begin{cases}
\Dqc (X) \big/ \left( {\Dqc} \right)_Z (X) 
\
\xrightarrow[\cong]{\overline{\Ld j^*}}
\
\Dqc (U)
\\
L 
= 
\Rd j_* \Ld j^* 
=
\left( \Rd j_*\mathcal{O}_U \right)
\otimes_{\mathcal{O}_X}^{\Ld} - :
 \ 
\Dqc (X) \ \to \ 
\Dqc (X) \big/ \left( {\Dqc} \right)_Z (X) 
\
\xrightarrow[\cong]{\overline{\Ld j^*}}
\
\Dqc (U)
\xrightarrow{\Rd j_*}
\Dqc (X)
\end{cases}
\end{equation*}
and the precursor in the setting of 
abelian categories reviewed in 
footnote 27:
\begin{equation*}
\QCoh (X) \big/ \QCoh_Z (X) 
\
\xrightarrow[\cong]{\overline{j^*}} 
\
\QCoh (U)
\end{equation*}
to the setting of $\Dperf,$ 
has been communicated to the author 
by Professor Neeman:

\begin{quote}
{\it
... But the right adjoints 
$j_* : \QCoh (U) \to \QCoh (X)$ and 
$\Rd j_* : \Dqc (U) \to \Dqc (X)$ fail to 
preserve the finite subcategories 
$\Coh (-)$ and $\Dperf (- ).$ 
For these categories some work is needed. Especially in the case of $\Dperf ( - )$;
for a long time all that was known was that 
$\Ld j^* : \Dperf (X) \to \Dperf(U)$ 
isn't surjective
on objects, hence the natural map 
\begin{equation*}
\frac{ \Dperf (X) }{ \Ker (\Ld j^*) } \ 
\longrightarrow \Dperf (U) \ 
\end{equation*}
couldn't be an equivalence. So the assumption was that this map had to be worthless.

Thomason's ingenious insight was that the old counterexamples were a red herring.
Up to idempotent completion this map is an equivalence, and in particular induces an
isomorphism in higher $K$-theory. This of course required proof. Thomason gave a rather
involved proof, following SGA6, and I noticed that the proof simplifies and generalizes
when one uses the methods of homotopy theory.

It was an amusing role reversal: Thomason, the homotopy theorist, had the brilliant
idea but gave a clumsy proof using the techniques of algebraic geometry, while I, the
algebraic geometer, simplified the argument with the techniques of homotopy theory.

}

\end{quote}

}

\begin{Theorem}[Thomason's localization 
theorem]
\label{Thomason's localization 
theorem}
Under the situation of \eqref{BL of TC}, 
i.e. let $X$ be a 
quasicompact and quasiseparated schume, 
$Z = X \setminus  U \subset X,$ 
the complement of a quasi-compact Zariski 
open immersion $j : U \hookrightarrow X,$ 
we have a triangulated embedding
\begin{equation*}
\Dperf (X) \big/ \left( {\Dperf} \right)_Z (X) \ \subset \ \Dperf (U),
\end{equation*}
which yields an equivalence upon applying
the thick closure:
\begin{equation}
\label{IC for Perf}
\left(
\Dperf (X) \big/ \left( {\Dperf} \right)_Z (X) 
\right)^{\widehat{}}
\
\xrightarrow[\cong]{\Ld j^*}
\
\Dperf (U).
\end{equation}
\end{Theorem}

In applications, we sometime have to take 
care of elements in 
$\left( {\Dperf} \right)_Z (X).$ 
Then we wonder if 
they are in the image of 
$\Rd i_* \Dperf (Z )$ or not.  
Now, Rouquier \cite{MR2434186} gave an 
affirmative answer for a weaker question 
in the coherent setting:

\begin{Theorem}
\label{Rouquier's Lemmma 7.1}
{\rm \cite[Lem.7.40]{MR2434186}}
Let $X$ be a separated noetherian scheme and 
 $Z$ be its closed subscheme 
given by the ideal sheaf $\mcI$ of $\mcO_X.$
For $n \in \mathbb{N},$ let $Z_n$ be the 
 closed subscheme of $X$ with ideal sheaf 
$\mcI^n$ and $i_n : Z_n\to X$ the corresponding immersion. Then,
\begin{equation*}
\forall Q\in \left( {\Dcoh} \right)_Z (X),
\qquad
\ \exists n \in \mathbb{N}, 
\exists P_n \in \Dcoh (Z_n)\  \text{s.t.}\quad
Q = \Rd {i_n}_* P_n.
\end{equation*}

\end{Theorem}

While the original proof of 
Theorem~\ref{Thomason's localization 
theorem} 
given in \cite{MR1106918} is purely 
algebro geometric in the spirit of SGA6, 
Neeman \cite[Th.2.1]{MR1191736} gave a proof  
from a general triangulated category 
theoretical point of view,  
in the homotopy theoretical spirit of 
Bousfield, Ohakawa, and others, 
building upon Corollary~\ref{Nee92b, Lem.1.7} \cite[Lem.1.7]{MR1191736}:

\begin{Theorem}[Neeman's generalization of 
Thomason's localization theorem]
\label{Neeman's generalization of Thomason's localization theorem}
Let $\mcT$ 
be 
a compactly generated 
triangulated category, 
generated by a set 
$K$ consisting of compact 
objects in $\mcT.$ 
For a subset $S \subseteq K,$ 
set
$\mcS$ 
be
 the smallest 
localizing triangulated subcategory containing $S.$ 
Then, the canonical sequence of 
triangulated cagtegories
\begin{equation} 
\label{canonical sequence of triangulated categories}
\mcS \to \mcT \to \mcT/\mcS
\end{equation}
induces another sequence of 
triangulated cagtegories of 
compact objects
\begin{equation}
\label{induced sequence of compact objects}
\mcS^c \to \mcT^c \to \left(\mcT/\mcS\right)^c,
\end{equation}
which induces an equivalence
\begin{equation}
\label{induced equivalence}
\mcS^c = \mcS \cap \mcT^c,
\end{equation}
a fully faithful embedding
\begin{equation}
\label{fully faithful embedding}
\mcT^c/\mcS^c \to (\mcT/\mcS)^c,
\end{equation}
and, although it may fail to induce
an equivalence 
$\mcT^c/\mcS^c \xrightarrow{\cong} 
\left(\mcT/\mcS\right)^c,$ it does
induce an equivalence upon 
applying the thick closure:
\begin{equation}
\label{equivalence upon thick closure}
\left(
\mcT^c/\mcS^c 
\right)^{\widehat{}}
\ 
\xrightarrow{\cong} 
\ 
\left(\mcT/\mcS\right)^c.
\end{equation}

\end{Theorem}


\begin{proof}
{\rm (i)}  The first triangulated functor 
in \eqref{induced sequence of compact objects} is an easy consequence of 
Proposition~\ref{compact objects observation}.
The second triangulated functor 
in \eqref{induced sequence of compact objects} is induced by the smashing
Bousfield localization functor
 $\mcT \to 
\mcT/\mcS,$ which preserves arbitrary coproducts 
Theorem~\ref{Nee92b, Lem.1.7}. 
Then for $c\in \mcT^c, 
t_{\lambda} \in \mcT\ (\lambda\in\Lambda ),$
regarding $\mcT/\mcS$ as the full subcategory of $L$-local objects, we evaluate as follows:
\begin{equation*}
\begin{split}
&\ 
\Hom_{\mcT/\mcS}\left(Lc, 
\oplus_{\lambda\in\Lambda}Lt_{\lambda}\right)
= \Hom_{\mcT}\left(Lc, 
\oplus_{\lambda\in\Lambda}Lt_{\lambda}\right)
\overset{\text{$L$:\ smashing}}=
 \Hom_{\mcT}\left(Lc,
L( \oplus_{\lambda\in\Lambda} t_{\lambda} )
\right)
\\
&=  \Hom_{\mcT}\left(c,
L( \oplus_{\lambda\in\Lambda} t_{\lambda} )
\right) \overset{\text{$L$:\ smashing}}= 
\Hom_{\mcT}\left(c,
\oplus_{\lambda\in\Lambda}Lt_{\lambda}\right)
\overset{c:\text{ compact}}=
\oplus_{\lambda\in\Lambda}
\Hom_{\mcT}\left(c,Lt_{\lambda}\right)
\\
&= \oplus_{\lambda\in\Lambda}
\Hom_{\mcT}\left(Lc,Lt_{\lambda}\right)
= \oplus_{\lambda\in\Lambda}
\Hom_{\mcT/\mcS}\left(Lc,Lt_{\lambda}\right),
\end{split}
\end{equation*}
which implies $Lc$ is also compact.\newline
On the other hand, Krause \cite{MR2681709} 
gave a conceptually simple, though more involved,  proof of 
the existence of  
\eqref{induced sequence of compact objects}, 
applying the following easy observation 
\cite[Lem.5.4.1.(1)]{MR2681709}, 
which goes back at least to 
\cite[Th.5.1]{MR1308405} where the 
converse, i.e. compactness 
preservation of $F$ $\implies$ 
small coproducts preservation of $G,$ 
is also 
shown under the additional 
compact generation
assumption of $\mcT$:

\begin{quote}
For any pair of adjoint triangulated 
functors $\xymatrix{ 
\mcT \ar@<1mm>[r]^F & 
\ar@<1mm>[l]^G  \mcU
}$ such that $G$ preserves small coproducts,
\\
$F$ preserves compactness. 
\\
\\
$\because$ ) In fact, for any $c\in \mcT^c, 
u_{\lambda} \in \mcU\ (\lambda\in\Lambda),$
\begin{equation*}
\begin{split}
\Hom_{\mcU}(Fc,\oplus_{\lambda}u_{\lambda})
&=
\Hom_{\mcT}(c,G(\oplus_{\lambda}u_{\lambda}))=
\Hom_{\mcT}(c,\oplus_{\lambda}G(u_{\lambda}))
\\
&=
\oplus_{\lambda}\Hom_{\mcT}(c,G(u_{\lambda}))=\oplus_{\lambda}\Hom_{\mcU}(Fc,u_{\lambda}).
\end{split}
\end{equation*}
\end{quote}
Now, \eqref{induced sequence of compact objects} is induced from 
\eqref{canonical sequence of triangulated categories} by applying this easy 
observation to the recollement given by
Proposition~\ref{smashing localization characterizations}.6.
\footnote{
This is the involved part of this proof, for
the existence of recollement there 
requires Brown representability.
}
\newline
{\rm (ii)}
To see \eqref{induced equivalence},
first note 
$\mcS^c \supset \mcS\cap \mcT^c$ is 
trivial from the definition.
Then \eqref{induced equivalence} follows 
since converse $\mcS^c \subset \mcS\cap \mcT^c$ also follows from 
\eqref{induced sequence of compact objects}.
\newline
{\rm (iii)} For \eqref{fully faithful embedding}, suffices to show the composite
\begin{equation*}
\Hom_{\mcT^c/\mcS^c}(c,c') \to 
\Hom_{(\mcT/\mcS)^c}(c,c') 
\xrightarrow{\cong}
\Hom_{(\mcT/\mcS)}(c,c') 
\overset{Th.~\ref{Nee92b, Lem.1.7}}=
\Hom_{\mcT}(c,\operatorname{hocolim} (x_n))
\end{equation*}
is an isomorphism.

For the surjectivity, take 
$\left(f : c \to 
\operatorname{hocolim} (x_n) \right) \in 
\Hom_{\mcT}(c,\operatorname{hocolim} (x_n)),$ then we can find its preimage 
$( c \xleftarrow{\bigstar} 
c\times^h_{c_n}c' \xrightarrow{
\widetilde{f}_n'
}
 c') \in 
\Hom_{\mcT^c/\mcS^c}(c,c')$ by a
straightforward contemplation 
summarized in the following
commutative diagram:
\begin{equation*}
\xymatrix{
 & c' \ar@{=}[r] \ar[d]^{\bigstar} 
 & x_0 \ar[d]
\\
c\times^h_{c_n}c' \ar[ur]^{\widetilde{f}_n'} 
\ar[d]_{\bigstar} & \exists c_n \ar[r] &
x_{n} \ar[d]
\\
c \ar@{..>}[ur]^{\exists \widetilde{f}_n}
 \ar[rr]_-f & & 
\operatorname{hocolim} (x_n)
}
\end{equation*}
Here, $c_n$ is some compact object 
so that 
arrows with $\bigstar$ have cones of 
the form  finite extension of finite 
coproducts of elements in $S,$ and  
$c\times^h_{c_n}c'$ is the homotopy 
pullback (see e.g.\cite[p.252,(1.1.2.5)]{MR1106918}). 

For the injectivity, suppose 
$(c \xleftarrow{\bigstar} \widetilde{c} 
\xrightarrow{f'}  c') \in 
\Hom_{\mcT^c/\mcS^c}(c,c')$ is sent to 
$(c \xleftarrow{\bigstar} \widetilde{c} 
\xrightarrow{0}  
\operatorname{hocolim}(x_n)
) =
0 \in \Hom_{\mcT}(c,\operatorname{hocolim} (x_n)).$ Then we can see 
$(c \xleftarrow{\bigstar} \widetilde{c} 
\xrightarrow{f'}  c') = 
( x \xleftarrow{\bigstar} 
\widetilde{c}\times^h_{c_m'}c' 
\xrightarrow{0} c' ) 
= 0 \in 
\Hom_{\mcT^c/\mcS^c}(c,c')$ 
by a
straightforward contemplation 
summarized in the 
following commutative diagram:
\begin{equation*}
\xymatrix{
\widetilde{c}\times^h_{c_m'}c' 
\ar[r]^{0} \ar[d]^{\bigstar} 
\ar@/_/[dd]_{\bigstar}
&
c' \ar[d]^{\bigstar} \ar@{=}[r]  & 
x_0 \ar[d]
\\
\widetilde{c} \ar[d]^{\bigstar} 
\ar[ur]^{f'} \ar[r]^{0} \ar[drr]^{0} & 
\exists c_m' \ar[r] & x_m \ar[d]
\\
c & &
\operatorname{hocolim} (x_n).
}
\end{equation*}
Here, $c_m'$ is some compact object 
so that 
arrows with $\bigstar$ have cones of 
the form  finite extension of finite 
coproducts of elements in $S,$ and  
$c\times^h_{c_m'}c'$ is the homotopy 
pullback \cite[p.252,(1.1.2.5)]{MR1106918}. 
\newline
{\rm (iv)} 
To see \eqref{equivalence upon thick closure}, write $\mcT = \langle K \rangle,$ and observe from the construction
of the Verdier quotient 
$\mcT \xrightarrow{F_{univ}} \mcT/\mcS$ 
that $\mcT/\mcS = \langle 
F_{univ}(K) \rangle,$ where 
$F_{univ}(K) \subseteq 
\mcT^c/\mcS^c \subseteq 
( \mcT/\mcS )^c$ 
by \eqref{induced sequence of compact objects} and 
\eqref{fully faithful embedding}.
Now apply Proposition~\ref{compact objects observation} to conclude any object $y$ of 
$( \mcT/\mcS )^c$ is a direct summand of 
a finite extension (in $( \mcT/\mcS )^c$)
of finite direct sums 
of objects in $F_{univ}(K)
\subseteq \mcT^c/\mcS^c,$ which is a 
full triangulated subcategory by 
\eqref{fully faithful embedding}.
This implies the desired equivalence 
upon thick closure 
\eqref{equivalence upon thick closure}: 
$\left(
\mcT^c/\mcS^c 
\right)^{\widehat{}}
\ 
\xrightarrow{\cong} 
\ 
\left(\mcT/\mcS\right)^c.$
\end{proof}

The following consequence of
Theorem~\ref{Thomason's localization 
theorem}
and 
Remark~\ref{localizing subcategory is thick} (iv) will be used later:

\begin{Corollary}
\label{A useful consequence of density}
Let $X$ be a Noetherian scheme, and 
$Z = X \setminus  U \subset X,$ 
the complement of a quasi-compact Zarisiki 
open immersion $j : U \hookrightarrow X.$
Then, for any $P \in \Dperf (U),$ 
there exists $H \in \Dperf (X)$ such that
\begin{equation*}
\Ld j^*H \cong P \oplus \Sigma P 
\in \Dperf (U).
\end{equation*}

\end{Corollary}

Now, to motivate Balmer's construction 
reviewed next, let us single out the 
 following slight strenghning of 
Theorem~\ref{Neeman's generalization of Thomason's localization theorem}
(and so also of 
Theorem~\ref{Thomason's localization 
theorem}):

\begin{Theorem} 
\label{strengthening of Thomason-Neeman}
Under the same assumption of 
Theorem~\ref{Neeman's generalization of Thomason's localization theorem}, 
the extrinsic thick closure equivalence 
\eqref{equivalence upon thick closure}
can be upgraded to the intrinsic 
idempotent completion
\footnote{
For the fact that the idempotent
 completion of a triangulated category has a
natural structure of a triangulated category, 
there is a proof in 
Balmer-Schlichting \cite{BS01}.
}
 equivalence:
\begin{equation}
\label{equivalence upon idempotent completion}
\left(
\mcT^c/\mcS^c 
\right)^{\sharp}
\ 
\xrightarrow{\cong} 
\ 
\left(\mcT/\mcS\right)^c.
\end{equation}
In particular, under the same assumption
of Theorem~\ref{Thomason's localization 
theorem}, we have an equivalence 
upon applying the idempotent
completion:

\begin{equation}
\label{IC for Perf - idempotent completion}
\left(
\Dperf (X) \big/ \left( {\Dperf} \right)_Z (X) 
\right)^{\sharp}
\
\xrightarrow[\cong]{\Ld j^*}
\
\Dperf (U).
\end{equation}

\end{Theorem}

To show \eqref{equivalence upon idempotent completion}, 
it suffices to show 
$\left(
\mcT^c/\mcS^c 
\right)^{\widehat{}}
\cong
\left(
\mcT^c/\mcS^c 
\right)^{\sharp}$
thanks to  
\eqref{equivalence upon thick closure}. 
For this, note from 
\eqref{fully faithful embedding} a 
fully faithful embedding
$\mcT^c/\mcS^c \to \mcT/\mcS.$ 
Here, $\mcT/\mcS$ is idempotent complete, 
because $\mcT/\mcS$ is first seen to be 
equipped with 
arbitrary small coproducts by 
Theorem~\ref{Nee92b, Lem.1.7}.2, 
Proposition~\ref{smashing localization characterizations}.5, 
Proposition~\ref{Bousfield localization in triangulated category}.5, and then 
we may apply Remark~\ref{localizing subcategory is thick}.(i) to 
find $\mcT/\mcS$ is idempotent complete. Thus, any added idempotent object of  $\left(
\mcT^c/\mcS^c 
\right)^{\sharp}$ shows up in  $\mcT/\mcS,$
but, because of $\mcT^c/\mcS^c \ \subseteq \ 
\left(\mcT/\mcS\right)^c$ and any direct
summand of a compact object is still 
compact, these added idempotent objects 
actually show up in 
$\left(\mcT/\mcS\right)^c.$
This implies the desired \eqref{equivalence upon idempotent completion}.

\item In view of 
Theorem~\ref{topological reconstruction}, 
we wonder whether the 
 spectrum $X$ is reconstructed 
from $(\Dperf(X), \otimes^{\mathbb{L}}).$ 
But, this is nothing but the theorem of 
Paul Balmer \cite{MR2196732}:

\begin{Definition} \label{tt-geometry}
For an essentially small tensor triangulated category $\mcK,$ 
we defined in Definition~\ref{tt-spectrum}
the spectrum (topological space) 
$\Spc (\mcK).$
\begin{itemize}
%
%
\item 
Here,  
motivated by 
\eqref{IC for Perf - idempotent completion}, 
we can construct a presheaf of tensor
 triangulated categories by
\begin{equation}
U \ \mapsto \ \mcK (U) :=
\left( \mcK \big/ \mcK_Z \right)^{\sharp},
\end{equation}
where $\mcK_Z := \left\{ a \in \mcK \ 
\mid \ \supp (a) \subseteq Z \right\}$ 
with $Z := X\setminus U$ and
$\supp (a) := \Spc (\mcK) \setminus
U(a) = \left\{ \mcP \in \Spc (\mcK) \ 
\mid \ a\notin \mcP \right\}.$

\item Finally, we obtain the ringed space
\begin{equation}
\Spec \left( \mcK \right) = 
\left(\Spc (\mcK), \mcO_{\mcK} \right),
\end{equation} 
as the sheafication of the preseaf of 
commutative rings
\begin{equation}
U \ \mapsto  \ \End_{\mcK (U)}
( \one ),
\end{equation}
where $\one$ is the unit object
of the tensor triangulated caegory 
$\mcK (U).$

\end{itemize}

\end{Definition}

Now Balmer's reconstruction theorem 
\cite{MR2196732} states:

\begin{Theorem} \label{Balmer's reconstruction theorem}
For a quasi-compact and quasi-separated
scheme $X,$ we have an isomorphism of
ringed spaces
\footnote{
The weaker reconstruction just as 
a topological space was already shown 
by Thomason (see 
Theorem~\ref{topological reconstruction}               
) in the course of his establishing a 
$\Dqc (X)$ analogue of the Hopkins-Smith
theorem (see 
Theorem~\ref{Thomason's theorem of Thomason sets}
Theorem~\ref{Tho=R,T}).
}

\begin{equation*}
\Spec 
\left( \Dperf (X), \otimes^{\Ld} \right)
\ \cong \ X .
\end{equation*}

\end{Theorem}

\end{itemize}

\subsection{$\Dcoh (X)$ and $\Dperf (X)$ 
determine each other}

With the concepts \lq\lq approximable\rq\rq, \lq\lq noetherian approximable\rq\rq, 
\lq\lq metric\rq\rq, 
\lq\lq preferred $t$-structure\rq\rq, 
and 
\lq\lq Cauchy sequence\rq\rq\ in a black box, 
Amnon Neeman's strategy to prove this 
may be summarized as follows:

\begin{itemize}
\item
\cite[Ex.8.4]{1806.06995}:

 Out of an \underline{approximable} triangulated
category $\mcT$ 
with a \underline{preferred $t$-structure} 
 $(\mcT^{\leq 0}, \mcT^{\geq 0}),$
we can construct a couple of 
triangulated categories $\mcS$ with 
\underline{metric}s 
:

\begin{enumerate}
\item $\mcS = \mcT^c \subset \mcT,$ and
$\mcM_i = \mcT^c\cap \mcT^{\leq -i}.$
\item $\mcS = [\mcT_c^b]^{\operatorname{op}}
,$ 
and
$\mcM_i^{\operatorname{op}} = \mcT_c^b\cap \mcT^{\leq -i}.$
\end{enumerate}

\item 
\cite[Def.1.10]{1806.06471}
For an essentially small triangulated
category $\mcS$ with a metric 
$\{ \mcM_i \},$ we define three
full subcategories 
$\mathfrak{L}(\mcS), 
\mathfrak{C}(\mcS), 
\mathfrak{S}(\mcS)$ of the category 
\begin{equation*}
\operatorname{Mod}-\mcS := 
\text{additive functors}\ 
\mcS^{\operatorname{op}} \to
\mathbb{Z}-\operatorname{Mod}.
\end{equation*}
With 
$Y : \mcS \to \operatorname{Mod}-\mcS;\  
A \mapsto Y(A) := \Hom(-,A)$ the 
\emph{Yoneda functor}, we set 
\begin{equation*}
\begin{split}
\mathfrak{L}(\mcS) &:=
\left\{ 
\underset{\longrightarrow}{\operatorname{colim}} 
Y(E_i) \in 
 \operatorname{Mod}-\mcS \ 
\big| \
E_*, \ \text{is a 
\underline{Cauchy sequence} in 
$\mcS$.}
\right\}
\\
\mathfrak{C}(\mcS) &:=
\left\{ A \in 
 \operatorname{Mod}-\mcS \ 
\big| \
\substack{
\text{
For every $j\in \mathbb{Z}$ there exists
$i\in \mathbb{Z}$ with
}
\\
\operatorname{Hom}
\left(Y(\mcM_i), \Sigma^{-j}A \right)
= 0.
}
\right\}
\\
\mathfrak{S}(\mcS) &:=
\mathfrak{L}(\mcS) \cap
\mathfrak{C}(\mcS).
\end{split}
\end{equation*}
By construction, we see 
\cite[Obs.2.3]{1806.06471}
\begin{equation*}
\mathfrak{S}(\mcS)  = 
\bigcap_{j\in \mathbb{Z}} 
\bigcup_{i\in \mathbb{N}}
\left[ Y(\Sigma^j E_i) \right]^{\perp}
\end{equation*}

Intuitively, $\mathfrak{S}(\mcS)$ 
consists of compactly supported 
objects (for contained in 
$\mathfrak{C}(\mcS)$) of the 
Cauchy completion with respect to 
the given metric inside the 
$\operatorname{Ind}$-completion given by 
 the Yoneda embedding (for contained in 
$\mathfrak{L}(\mcS).$

Appriori, it is not clear whether 
$\mathfrak{S}(\mcS)$ is triangulated 
or not.  However, Neeman proves:

\begin{Theorem}{\rm 
\cite[Def.2.10,Th.2.11]{1806.06471}}
$\mathfrak{S}(\mcS)$ becomes a triangulated 
category with the distinguished triangles
of the form
$
\underset{\longrightarrow}
{\operatorname{colim}} 
Y(
A_i \xrightarrow{f_i} 
B_i \xrightarrow{g_i} 
C_i \xrightarrow{h_i} 
\Sigma A_i
),
$
where 
$(
A_* \xrightarrow{f_*} 
B_* \xrightarrow{g_*} 
C_* \xrightarrow{h_*} 
\Sigma A_*
)$
is a Cauchy sequence of 
triangles in $\mcS.$

\end{Theorem}

\item 
\cite[Th.8.8]{1806.06995}
With the metrics as 
above, we have 
triangulated equivalences
\begin{enumerate}
\item $\mathfrak{S}(\mcT^c) = \mcT_c^b.$
\item If $\mcT$ is noetherian then
$\mathfrak{S}\left( 
[\mcT_c^b]^{\operatorname{op}} \right) =
[ \mcT^c ]^{\operatorname{op}}.$
\end{enumerate}

\item 
\cite[Ex.3.6]{1804.02240} 
The above theory 
works when $X$ is separated and quasi-compact: 
If $X$ is separated and quasi-compact,
$\mcT = \Dqc (X)$ is approximable 
with the standard $t$-structure in the preferred equivalence class.

\item 
Consequently, we obtain our 
desired result:

When $X$ is separated and quasi-compact,
we have the following:
\begin{enumerate}
\item 
$\mathfrak{S}(\Dperf (X) ) =   \Dcoh (X).$

\item If $X$ is further noetherian,
$\mathfrak{S}\left( 
 \left[ \Dcoh (X) \right]^{\operatorname{op}}
\right)
=  \left[ \Dperf (X) \right]^{\operatorname{op}}.$
\end{enumerate}

\end{itemize}

For the rest of this section, we 
explain the concepts of 
\lq\lq approximable\rq\rq, \lq\lq noetherian approximable\rq\rq,
\lq\lq metric\rq\rq, 
\lq\lq preferred $t$-structure\rq\rq, 
and 
\lq\lq Cauchy sequence\rq\rq, which were put
 in a black box in the above summary. 
We urge readers  to consult Neeman's own survey 
\cite{1806.06995} for more details
about the approximable triangulated 
categories.

Now, it is rather straightforward to 
define \lq\lq metric\rq\rq \ and 
\lq\lq Cauchy sequence\rq\rq.

\begin{Definition}
{\rm
\cite[Def.1.2]{1806.06471}
\cite[Def.8.3]{1806.06995}
}
A \emph{metric} on a triangulated
category $\mcS$ is a sequence of additive 
subcategories 
$\{ \mcM_i, i\in \mathbb{N} \},$ 
satisfying:
\begin{enumerate}
\item $\mcM_{i+1} \subset 
\mcM_i$ for every $i\in \mathbb{N}.$
\item Any $b\in \mcS,$ with a distinguished 
triangle $a \to b \to c$ s.t. 
$a,c \in \mcM_i,$ belongs to 
$\mcM_i.$
\end{enumerate}

\end{Definition}

\begin{Definition}
{\rm
\cite[Def.1.6]{1806.06471}
\cite[Def.8.5]{1806.06995}
}
A \emph{Cauchy sequence} in $\mcS,$ 
a triangulated category with a metric
$\{ \mcM_i \},$ is a sequence
\begin{equation*}
E_1 \to E_2 \to E_3 \to \cdots
\end{equation*}
such that, for any $i \in \mathbb{N},
j \in \mathbb{Z},$ there exists
$M \in \mathbb{N}$ such that,
\begin{equation*}
\operatorname{Cof} 
( E_m \to E_{m'} ) \in \Sigma^{-j}\mcM_i
\end{equation*}
for any $m' > m \geq M.$

\end{Definition}

Next, we aim at 
\lq\lq preferred $t$-structure\rq\rq, 
but we shall make a little detour for 
some later purpose.
%

\begin{Definition} 
{\rm
\cite[Rem.3.1]{1806.06995}
}
	Let $\mcA$ be a full 
subcategory of a category $\mcT.$ 
Define the full subcategories $\add\mcA$, $\Add\mcA$, and $\smd\mcA$ as follows.
	\begin{enumerate}
		\item 
Assume $\mcT$ has finite coproducts. 
$\add\mcA$ consists of all finite 
{coproducts} 
of objects in $\mcA$.
		\item 
Assume $\mcT$ has coproducts. 
$\Add\mcA$ consists of all coproducts of objects in $\mcA$.
		\item $\smd\mcA$ consists of all direct summands in $\mcT$ 
of objects in $\mcA$.
	\end{enumerate}
\end{Definition}

The following construction will 
play major roles:

\begin{Definition}
\label{<>}
{\rm
\cite[Def.3.3]{1806.06995}
\cite[Rem.0.1]{1703.04484}
}
\label{Def.3.3}
Given $\mcA \subset \mcT,$  a full
subcategory of a triangulated category, 
and possibly infinite integers $m\leq n,$ 
define the full subcategories:

\begin{enumerate}
\item $\mcA[m,n] = \cup_{i=m}^n
\mcA[-i].$
\item For $l \in \mathbb{N}$, 
define inductively 
the full subcategory $\generatedset{\mcA}_{l}^{[m,n]}$ (resp.\ $\overline{\generatedset{\mcA}}_{l}^{[m,n]}$ if $\mcT$ has coproducts) as follows.
		
\begin{enumerate}
\item $\generatedset{\mcA}_{1}^{[m,n]}
=\smd(\add\mcA{[m,n]})$\quad 
{
(resp.\
$
\overline{
\generatedset{\mcA}
}_{1}^{[m,n]}=\smd(\Add\mcA {[m,n]})$
}
			\item $\generatedset{\mcA}_{l+1}^{[m,n]}=
\smd(\generatedset{\mcA}_{1}^{[m,n]}*
\generatedset{\mcA}_{l}^{[m,n]})$\quad 
(resp.\ $\overline{\generatedset{\mcA}}_{l+1}^{[m,n]}=
\smd(\overline{\generatedset{\mcA}}_{1}^{[m,n]}
*\overline{\generatedset{\mcA}}_{l}^{[m,n]})$).
		\end{enumerate}

\item 
For the case $m= - \infty, n= \infty$ and 
$l\in \mathbb{N},$ 
following Bondal-Van den Bergh \cite{MR1996800}, 
we shall 
simply denote as follows:
\footnote{
It was Neeman's insight to notice 
surprising usefullness of introducing 
related categories
$\generatedset{\mcA}_{l}^{[m,n]}$ and  $\overline{\generatedset{\mcA}}_{l}^{[m,n]}$
as well.
}
$$
\langle \mcA \rangle_l 
:= \langle \mcA \rangle_l^{[-\infty,\infty]} \quad 
(resp.\
\overline{\langle \mcA \rangle}_l 
:= \overline{\langle \mcA \rangle}_l^{[-\infty,\infty]}
)
$$
\end{enumerate}

\end{Definition}

Whereas the above definition might look 
complicated, its major part is 
reflected in the following simpler 
definition:

\begin{Definition} 
\label{coprod}
{\rm
\cite[Def.1.3]{{1703.04484}}
}
Given $\mcA \subset \mcT,$  a full
subcategory of a triangulated category, 
and  $l \in \mathbb{N}$, 
define inductively 
the full subcategory 
$\operatorname{coprod}_l(\mcA)$ 
(resp.\ 
$\operatorname{Coprod}_l(\mcA)$  if $\mcT$ has coproducts) as follows.
		
\begin{enumerate}
\item $\operatorname{coprod}_1(\mcA)
= \add(\mcA)$ \quad 
{
(resp.\
$
\operatorname{Coprod}_1(\mcA)
= \Add(\mcA),$
}
			\item 
$\operatorname{coprod}_{l+1}(\mcA) = 
\operatorname{coprod}_1(\mcA) * 
\operatorname{coprod}_l(\mcA)$\quad
(resp.\
$\operatorname{Coprod}_{l+1}(\mcA) = 
\operatorname{Coprod}_1(\mcA) * 
\operatorname{Coprod}_l(\mcA).$
)
		\end{enumerate}

\end{Definition}

The key for Definition~\ref{coprod} 
to reflect a major part of 
Definition~\ref{<>} is the following 
elementary observation of Bondal-Van den 
Bergh \cite{MR1996800}:

\begin{Lemma} 
\label{Bondal-Van den Bergh lemma}
{\rm \cite[Lem.2.2.1]{MR1996800}}
Let $\mcA$ and $\mcB$ be full 
subcategories of a triangulated category 
with small coproducts. Then:
\begin{enumerate}
\item[{\rm (1)}] $\smd (\mcA) * \mcB \subset
\smd( \mcA * \mcB ),\quad 
\mcA * \smd( \mcB ) \subset 
\smd (\mcA*  \mcB)$;
\item[{\rm (2)}] $\smd \left( \smd (\mcA) * \mcB \right)
= \smd \left( \mcA * \smd( \mcB ) \right)
= \smd (\mcA * \mcB).$
\end{enumerate}

\end{Lemma}

To show the first inclusion of (1): 
 $\smd (\mcA) * \mcB \subset
\smd( \mcA * \mcB ),$ pick 
$x \in \smd (\mcA) * \mcB$ fitting in 
a triangle:
\begin{equation*}
s \to x \to b \quad (s\in \smd(\mcA), 
b \in \mcB ),
\end{equation*}
for which we pick $s'\in \mcT$ with 
$s\oplus s' \in \mcA$ and form a 
new triangle:
\begin{equation*}
s\oplus s' \to x\oplus s' \to b.
\end{equation*}
This shows the desired 
$x \in \smd (\mcA * \mcB).$
The second inclusion of (1): 
$\mcA * \smd( \mcB ) \subset 
\smd (\mcS* \mcB)$ is shown similarly.
Then (2) follows immediately from (1).

Using
Lemma~\ref{Bondal-Van den Bergh lemma},
we can easily prove, by 
induction on $l,$ 
the following 
transparent expression relating 
Definition~\ref{<>} with 
Definition~\ref{coprod}.

\begin{Corollary}
\label{coproduct expression of <>}
{\rm (c.f. \cite[Cor.1.11]{1703.04484})}
Given $\mcA \subset \mcT,$  a full
subcategory of a triangulated category, 
a natural number $l\in \mathbb{N},$ 
and possibly infinite integers $m\leq n,$ 
\begin{equation*}
\generatedset{\mcA}
_{l}^{[m,n]} = 
\smd \left( \operatorname{coprod}_l
\mcA
[m,n] \right),
\qquad 
\overline{
\generatedset{\mcA}
}_{l}^{[m,n]} = 
\smd \left( \operatorname{Coprod}_l
\mcA
[m,n] \right).
\end{equation*}

\end{Corollary}

%
%

The following Proposition~\ref{Nee17,Cor.1.11}
follows immediately by combining the second equality of
Corollary~\ref{coproduct expression of <>} 
and Lemma~\ref{Nee17,Lem.1.9} 
below. Philosophically Proposition 4.31 may be viewed as saying
that $\overline{\langle - \rangle}_l$ and 
$\operatorname{Coprod_l}(-)$  
are interchangeable.


\begin{Proposition}
\label{Nee17,Cor.1.11}
{\rm (c.f. \cite[Cor.1.11]{1703.04484})}
Given $\mcA \subset \mcT,$  a full
subcategory of a triangulated category, 
a natural number $l\in \mathbb{N},$ 
and possibly infinite integers $m\leq n,$ 

\begin{equation*}
\operatorname{Coprod}_l
\left( \mcA [m,n] \right) \ \subseteq \
\overline{
\generatedset{\mcA}
}_{l}^{[m,n]} 
\ \subseteq \
\operatorname{Coprod}_{2l}
\left( \mcA [m-1,n] \right).
\end{equation*}

\end{Proposition}

We include a proof of the following 
Lemma~\ref{Nee17,Lem.1.9}, 
to highlight the point at which infinite 
 coproducts are
used. Just in case the reader is wondering: 
 the finite analogue of
Proposition~\ref{Nee17,Cor.1.11} is 
false. While the inclusion 
$\operatorname{coprod}_l
\left( \mcA [m,n] \right) \ 
\subseteq \
\generatedset{\mcA}
_{l}^{[m,n]}$  
is true and easy, it isn't in general
true that 
$\generatedset{\mcA}
_{l}^{[m,n]} \ \subseteq \ 
\operatorname{coprod}_{2l}
\left( \mcA [m-1,n] \right).$
%
%

\begin{Lemma} 
\label{Nee17,Lem.1.9}
{\rm (c.f. \cite[Lem.1.9]{1703.04484})}
Let $\mcB$ a subcategory of $\mcT,$ 
a triangulated category with coproducts, 
and $l\in \mathbb{N}.$
Then
\begin{equation*}
\operatorname{Coprod}_l( \mcB ) \ 
\subseteq \ \smd \left( 
\operatorname{Coprod}_l( \mcB ) \right) \
\subseteq \ 
\operatorname{Coprod}_{2l}( \mcB [-1,0]  ).
\end{equation*}

\end{Lemma}

\begin{proof} 
The first inclusion is obvious.
For the second inclusion, recall from 
Remark~\ref{localizing subcategory is thick}(i)
that 
\begin{equation*}
\begin{split}
\forall x \in \smd \left( 
\operatorname{Coprod}_l( \mcB ) \right),
\quad 
&\exists b \in 
\operatorname{Coprod}_l( \mcB )
\ \text{and an idempotent}\ 
e : \to b,\ 
\\
&\qquad \text{s.t.} \
x = eb= \operatorname{Cone}
\left( \oplus_{\mathbb{N}} b \to
\oplus_{\mathbb{N}} b  \right).
\end{split}
\end{equation*}

From this, we obtain the following  
triangle:
\begin{equation*}
 \oplus_{\mathbb{N}} b \to 
 \oplus_{\mathbb{N}} b \to x \to
\Sigma  \left( \oplus_{\mathbb{N}} b \right),
\end{equation*}
where $\oplus_{\mathbb{N}} b \in 
\Add \left( 
\operatorname{Coprod}_l( \mcB )
\right) = \operatorname{Coprod}_l( \mcB )$ 
and so $\Sigma  \left( \oplus_{\mathbb{N}} b \right) \in 
\Sigma \operatorname{Coprod}_l
(  \mcB )
=
\operatorname{Coprod}_l
( \Sigma \mcB ).$
Thus,
\begin{equation*}
x \in \operatorname{Coprod}_l( \mcB ) *
 \operatorname{Coprod}_l
(\Sigma  \mcB ) \subseteq 
\operatorname{Coprod}_l
( \mcB \cup \Sigma \mcB) *
 \operatorname{Coprod}_l
( \mcB \cup \Sigma \mcB )
\subseteq 
 \operatorname{Coprod}_{2l}
( \mcB \Sigma \mcB ). 
\end{equation*}
\end{proof}

The constructions 
$\langle - \rangle_l$ and 
$\overline{\langle - \rangle}_l$ are 
older than 
$\operatorname{coprod}_l(-)$ and  
$\operatorname{Coprod}_l(-)$, 
and for most purposes they work just fine. 
 But there are results which become much 
 easier to prove by working with
$\operatorname{coprod}_l(-)$ and  
$\operatorname{Coprod}_l(-)$; 
for example the reader can look at
the proof of \cite[Lem.4.4]{BNP18}.
\footnote{The author is grateful to 
Professor Neeman for this reference.}
Thus one way to view the difference is to
regard
$\operatorname{coprod}_l(-)$ and  
$\operatorname{Coprod}_l(-)$ 
as technically more powerful than the older
$\langle - \rangle_l$ and 
$\overline{\langle - \rangle}_l.$

Now, in practice, 
as 
their constructions 
suggest, $\operatorname{coprod}_l \ 
\ (\text{resp.} \operatorname{Coprod}_l)$ 
are more 
tractible than 
$\generatedset{\mcA}_{l}^{[m,n]}\ 
\ (\text{resp.} 
\overline{\generatedset{\mcA}}_{l}^{[m,n]}).$
However, $\generatedset{\mcA}_{l}^{[m,n]}\ 
\ (\text{resp.} 
\overline{\generatedset{\mcA}}_{l}^{[m,n]}).$ 
occurs more frequently, for instance,

\begin{Theorem} \label{Ex.0.13}
{\rm
\cite[Th.A]{MR1974001}
 (See also \cite[Ex.0.13]{1804.02240})
}
For a triangulated category  $\mcT$ with coproducts and a compact generator 
$G \in \mcT,$ there is a unique
$t$-structure of the following form:
\begin{equation*}
\left( \mcT_G^{\leq 0}, 
 \mcT_G^{\geq 0} \right) :=
\left( 
\overline{\langle G \rangle }^{[-\infty,0]}, 
\left( 
\overline{\langle G \rangle }^{[-\infty,0]}
\right)^{\perp}[1] 
\right).
\end{equation*}

\end{Theorem}

\begin{Definition} 
{\rm
\cite[Def.7.3, Rem.7.4]{1806.06995}
}
\begin{enumerate}
\item Two $t$-structures 
$\left( \mcT_1^{\leq 0}, 
 \mcT_1^{\geq 0} \right)$ and
$\left( \mcT_2^{\leq 0}, 
 \mcT_2^{\geq 0} \right)$ are called
\emph{equivalent}, if there exists 
$A\in \mathbb{N}$ with
\begin{equation*}
\mcT_1^{\leq -A}  \subset
\mcT_2^{\leq 0} \subset
\mcT_1^{\leq A} .
\end{equation*}

\item For a triangulated category  $\mcT$ with coproducts and a compact generator, 
a $t$-structure 
$\left( \mcT^{\leq 0}, 
 \mcT^{\geq 0} \right)$ 
is in the 
\emph{preferred equivalence class }if it 
 is equivalent to
$\left( \mcT_G^{\leq 0}, 
 \mcT_G^{\geq 0} \right)$ 
for some compact generator $G$ 
(in fact, for every compact generator).

\end{enumerate}

\end{Definition}


The importance of \lq\lq preferred 
equivalence class\rq\rq \ is that 
$\mcT^-, \mcT^+,$ and $\mcT^b,$ recalled in
the next definition, are independent of 
the particular representative 
 $(\mcT^{\leq 0}, \mcT^{\geq 0})$ 
in the preferred equivalence class
\cite[Fact.0.5.(iii)]{1806.06995}:

\begin{Definition}
{\rm
\cite[Def.7.5, Def.7.6]{1806.06995}
}
\begin{enumerate}
\item Given a $t$-structure 
 $(\mcT^{\leq 0}, \mcT^{\geq 0}),$
we have the usual subcategories:
\begin{equation*}
\mcT^- = \cup_n \mcT^{\leq n},\quad
\mcT^+ = \cup_n \mcT^{\geq n},\quad
\mcT^b = \mcT^-\cap \mcT^+.
\end{equation*}

\item For a triangulated category  $\mcT$ with coproducts and a compact generator, 
choose a $t$-structure 
$(\mcT^{\leq 0}, \mcT^{\geq 0})$ 
in the preferred equivalence class, define
the full subcategories $\mcT_c^-$ and 
$\mcT_c^b$ as follows:
\begin{equation*}
\mcT_c^- 
:= \left\{ F \in \mcT \ 
\bigg| \ 
\substack{
\text{For any $n\in \mathbb{N}$ 
there exists a triangle}
\\
E \to F \to D \to E[1]
\\
\text{with $E$ compact and $D \in
\mcT^{\leq -n-1}$ }
}
\right\},\qquad
\mcT_c^b 
:= \mcT^b\cap \mcT_c^-
\end{equation*}

\end{enumerate}

\end{Definition}

Intuitively, $\mcT_c^-$ is the closure, 
with respect to the metric 
$\mcM_i = \mcT^{\leq -i},$ 
of $\mcT^c.$ 

$\mcT_c^-$ and $\mcT_c^b$ 
in the above definition
do not depend on the choice of compact 
generator $G$ and are both intrinsic
\cite[Rem.7.7,Fact.0.5.(iv)]{1806.06995}.

Now we are ready state 
the fundamental concepts of 
\lq\lq approximable\rq\rq \ and \lq\lq noetherian (approximable)\rq\rq \ 
triangulated categories:

\begin{Definition}
{\rm
\cite[Def.0.21]{1804.02240}
\cite[Def.4.1]{1806.06995}
}
A triangulated category $\mcT$ 
with coproducts is called 
\underline{\emph{approximable}} if
there exits 
a compact generator $G \in \mcT,$ 
a $t$-structure $\left( \mcT^{\leq 0},
\mcT^{\geq 0} \right),$ and $A \in
\mathbb{N}$ such that
\begin{enumerate}
\item $G[A] \in \mcT^{\leq 0}$ and
$\Hom ( G[-A], \mcT^{\leq 0} ) = 0.$

\item For every object 
$F\in \mcT^{\leq 0},$
there exists a triangle 
\begin{equation*}
E \to F \to D \to E[1],
\end{equation*}
with $D\in \mcT^{\leq -1}$ and
$E \in 
\overline{\langle G\rangle}_A^{[-A,A]}.$

\end{enumerate} 

\end{Definition}

From the definition, we find  
for any approximable
triangulated category $\mcT,$ 
the closure, with respect to the metric 
$\mcM_i = \mcT^{\leq -i},$ of \ 
$\bigcup_n 
\overline{\langle G\rangle}_n^{[-n,n]}$ 
is nothing but $\mcT^{-}.$ Thus we may 
intuitevely say every object in 
$\mcT^{-}$ may be \lq\lq Taylor approximable\rq\rq \ regarding 
$
{\langle G\rangle}_n^{[-n,n]}$  
as consisting of 
\lq\lq Taylor polynomials of $G$ 
of degree $\leq n.$\rq\rq \ 
\cite[Dis.0.1,Rem.02]{1806.06995}.

\begin{Definition}
\label{noetherian triangulated category}
{\rm \cite[Def.5.1]{1806.06471}
\cite[Not.8.9]{1806.06995}}
Suppose $\mcT$ is a triangulated category with coproducts, and assume it has
a compact generator $G$ with 
$\Hom(G,\Sigma^iG) = 0$ for $i \gg 0.$
We declare $\mcT$ to be 
\underline{\em noetherian} if there
exists $N \in \mathbb{N}$ and a 
$t$--structure $\left( \mcT^{\leq 0}, 
 \mcT^{\geq 0} \right)$ in the preferred
 equivalence class, s.t.
\begin{equation*}
\forall X\in \mcT_c^{-},\ \exists 
\ \text{triangle} \
A \to X \to B\ \text{s.t.} \
A \in \mcT_c^{-}\cap \mcT^{\leq 0},\ 
B \in \mcT_c^{-}\cap \mcT^{\geq -N} =
\mcT_c^b\cap \mcT^{\geq -N}.
\end{equation*}

\end{Definition}

\begin{Remark} 
\label{noetherian remark}
{\rm (i)} 
The noetherian hypothesis is somewhat 
 weaker than the assumption that there 
 exists a $t$--structure in the preferred
 equivalence class which restricts to a 
 $t$--structure on $\mcT_c^-.$
\newline
{\rm (ii) \cite[Fac.0.23,Exa.3.6]{1804.02240}}
For a quasicompact and separated scheme 
$X,$ the standard $t$-structure on 
$\mcT = \Dqc (X)$ 
is in 
the preferred equvalence class.
Suppose furrther that 
\underline{$X$ is noetherian}, then
$\mcT_c^{-} = \DD^{-}_{\mathrm{coh}},$
the category of bounded-above complexes 
of coherent sheaves, and so, 
the standard $t$--structure, which is in
the preferred equivalence class, on 
$\mcT = \Dqc (X)$ restricts to 
a $t$--structure on $\mcT = \Dqc (X).$
This implies $\Dqc (X)$ becomes 
noetherian in the sense of 
 Definition~\ref{noetherian triangulated category}, provided $X$ is noetherian 
and separated.  This is the origin of 
the terminology \lq\lq noetherian\rq\rq \ 
of Definition~\ref{noetherian triangulated category}.
\newline
{\rm (iii)}
\underline{\em WARNING!} The  
\lq\lq noetherian\rq\rq \ triangulated 
category
of Definition~\ref{noetherian triangulated category} is nothing to do with 
the \lq\lq Noetherian\rq\rq \ 
stable homotopy category of 
\cite[Def.6.0.1]{HPS97}.

For instance, for the case of 
$\mcT = \mcS \mcH,$ the stable homotopy
category of spectra, it is easy to see 
$\mcT_c^-$ consists of those spectra $X$ 
whose homotpy group $\pi_i(X)$ is a
 finitely generated abelian groups for 
each $i$ and vanishes for $i \ll 0.$
Thus, the standard $t$-structure, 
which is obviously in the preferred 
equivalence class, restricts to a
$t$--structure on $\mcT_c^-.$
This implies $\mcS \mcH$  is 
noetherian in the sense of 
 Definition~\ref{noetherian triangulated category} 
\cite[Fac.0.23]{1804.02240}.

 On the other hand, 
$\mcS \mcH$ is clearly NOT a 
Noetherian stable homotopy category in 
the sense of \cite[Def.6.0.1]{HPS97}, 
for the graded ring of the stable homotopy 
category of spheres $\pi_* S^0$ is 
not a Noetherian graded commutative ring, 
which can be easily seen by applying the
Nishida nilpotency, the precursor of 
(Devinatz-)Hopkins-Smith nilpotency.
 \end{Remark}

Then we have the following 
somewhat straightforward result 
to produce examples of approximable 
triangulated categories:

\begin{Proposition}
{\rm \cite[Ex.3.3]{1804.02240}}

If $\mcT$ has a
compact generator $G,$ such that 
$\Hom(G,\Sigma^iG) = 0$ for all $i > 0,$ 
then $\mcT$ is approximable.
Just take the $t$–-structure
$\left( \mcT_G^{\leq 0}, 
 \mcT_G^{\geq 0} \right)$ of 
Theorem~\ref{Ex.0.13} 
with $A=1.$

\end{Proposition}

From this, we immediately see 
the stable homotopy category 
$\mcS \mcH$ is approximable.
(actually noetherian, as was remarked 
in Remark~\ref{noetherian remark}(iii)).

Our  principal example of approximable 
triangulated categories is supplied 
by the following theorem:

\begin{Theorem}
{\rm 
\cite[Ex.3.6]{1804.02240}
}

Let $X$ be a quasicompact, 
separated
\footnote{
Unlike 
\eqref{BL of TC} and  
Theorem~\ref{Thomason's localization theorem}, 
the general case (where $X$ is quasicompact and quasiseparated) is still open -- see
\cite[Just above Lem.3.5]{1804.02240}.
%
%
}
 scheme. Then the category 
$\Dqc(X)$ is approximable.
(actually noetherian if $X$ is further 
noetherian, as was remarked 
in Remark~\ref{noetherian remark}(ii)).

\end{Theorem}

The proof is very involved and 
we urge readers  to consult Neeman's 
original paper 
\cite{1804.02240}.

For now, we shall record the following 
application of approximability:

\begin{Corollary}
\label{18d,Lem.6.5}
{\rm \cite[Lem.6.5]{1806.06995}
\cite[Th.0.18]{1703.04484}
}
Let $X$ be a quasicompact, separated scheme, 
let $G \in \Dqc (X)$ be a compact
generator, and let $u : U \to X$ be an
open immersion with $U$ quasicompact. Then
\begin{equation*}
\exists n\in \mathbb{N}\ \text{s.t.} \quad  
\Rd u_*\mcO_U \in 
\overline{\langle G\rangle}_n^{[-n,n]} 
 \subset \Dqc (X).
\end{equation*}

\end{Corollary}

\begin{proof}[Outline of the proof of 
Corollary~\ref{18d,Lem.6.5} using 
approximability presented in 
\cite
{1806.06995}:]
\ 
\begin{description}
\item[\underline{Step 1}]
$\exists l\in \mathbb{N}$ s.t. 
$\Hom \left( \Rd u_*\mcO_U, 
\Dqc (X)^{\leq -l} \right)=0.$

\item[\underline{Step 2 (This is where the approximability
 of $\Dqc(X)$ is used!)}]
By 
he approximability
 of $\Dqc(X)$,
\footnote{
There is some subtlety here.
See e.g. 
\cite[footnote 4 in Proof of Lem.5; 
Sketch 7.19.(i)]{1806.06995}
}
$\exists n\in \mathbb{N}$ and a triangle:
\begin{equation*}
E \to \Rd u_* \mcO_U  
\to D
\end{equation*}
with $D\in \Dqc(X)^{\leq -l}$ and 
$E\in \overline{\langle G\rangle}_n^{[-n,n]}.$

\item[\underline{Step 3}]
From Step 1 and Step 2, the map 
$\Rd u_* \mcO_U \to D$ in Step 2 is $0,$ 
which implies $\Rd u_* \mcO_U$ is a 
direct summand of $E\in 
\overline{\langle G\rangle}_n^{[-n,n]},$ 
as desired. 
\end{description}
\end{proof}

For details about the approximable 
triangulated categories. consult Neeman's 
own survey 
\cite{1806.06995}.

\section{Strong generation in derived categories of schemes}



In the previous section, we saw 
$\Dperf(X)$ and $\Dcoh(X)$ carry rich 
 information  and are intimately related to each
other. In this section, we would like to 
investigate the important 
\lq\lq strong generation\rq\rq \ 
property,  
in the sense of Bondal and Van den Bergh 
\cite{{MR1996800}}, 
for $\Dperf(X)$ and $\Dcoh(X),$ 
via 
approximable triangulated 
category techniques.





For this purpose, we have to start with 
what we mean by a \lq\lq generator\rq\rq of 
$\Dperf (X)$ and $\Dcoh (X),$ because  
our previous definition of 
a generator in 
Definition~\ref{generation when there are small coproducts} only works for 
triangulated categories with small coproducts, which $\Dperf (X)$ and $\Dcoh (X)$ are not.

\begin{Definition}
{\rm
\cite[Expl.5.4]{1806.06995}
}
Let $G$ be an element of a triangulated category 
$\mcS.$ Then, in the notation of Definition~\ref{Def.3.3},
\begin{enumerate}
\item $G$ is called a 
\emph{classical generator} if 
$\mcS = \cup_n \langle G\rangle_n^{[-n,n]}.$
\item $G$ is called a 
\emph{strong generator} if 
there exists an integer $l > 0$ with
$\mcS = \cup_n \langle G\rangle_l^{[-n,n]}.$
In this case, $\mcS$ is called 
\emph{strongly generated}.

\end{enumerate}

\end{Definition}

With this opportunity, let us record the 
following important concept intimately related 
to the above definition:

\begin{Definition} 
\label{Rouquier dimension}
{\rm
\cite[Def.3.2]{MR2434186}
}
The \emph{Rouquier dimension} of a 
triangulated category $\mcS,$ 
denoted by $\dim \mcS,$ 
is the smallest $d$ for which there exists $G\in\mcS$ with
$\mcS=\cup_n \langle G\rangle_{d+1}^{[-n,n]}.$

\end{Definition}

\begin{Remark}
{\rm (i)} 
Rouquier {\rm \cite{MR2434186}} 
proved the following 
properties of the Rouquier dimension of
$\Dcoh (X):$
\begin{itemize}
\item
{\rm \cite[Prop.7.9]{MR2434186}}
For a smooth quasiprojective
scheme $X$ over a field, we have 
$\dim \Dcoh (X) \leq  2 \dim X.$

\item 
{\rm \cite[Prop.7.16]{MR2434186}}
For  a reduced separated scheme $X$ of finite type over a field,
$\dim \Dcoh (X) \geq \dim X.$

\item 
{\rm \cite[Th.7.17]{MR2434186}}
For a smooth affine scheme $X$ of 
finite type over a field,
$\dim \Dcoh (X) = \dim X.$

\end{itemize}
{\rm (ii)} 
 For a sample of examples of 
Rouquier dimension in affine case, see 
{\rm \cite{IT14}\cite{DT15a} \cite{DT15b} }
for instance.

\end{Remark}

On the other hand, 
Neeman deduces strong generation of 
$\Dperf(X)$ and $\Dcoh(X)$ from 
some properties of $\Dqc (X)$:

{



\begin{Definition}
\label{strongly compactly or boundedly generated}
	Let $X$ be a separated scheme.
	\begin{enumerate}
		\item $\Dqc(X)$ is called \emph{strongly compactly generated} if there exists $G\in\Dperf(X)$ and and integer $l>0$ with
$\Dqc (X) = 
\overline{
\langle G\rangle
}^{(-\infty,\infty)}_l.$

\item  $\Dqc(X)$ is called \emph{strongly boundedly generated} if there exists $G\in\Dcoh(X)$ and and integer $l>0$ with
$\Dqc (X) = 
\overline{
\langle G\rangle
}^{(-\infty,\infty)}_l.$

\end{enumerate}

\end{Definition}

\begin{Remark}
\label{remark on strongly compactly or boundedly generated}
From Proposition~\ref{Nee17,Cor.1.11}, 
we may replace the required equality
$\Dqc (X) = 
\overline{
\langle G\rangle
}^{(-\infty,\infty)}_l$ showing up twice in
Defintion~\ref{strongly compactly or boundedly generated} with more tractible
$\Dqc(X) =
 \operatorname{Coprod}_l
\left(G(-\infty,\infty)\right)$ 
(of course, $l$ here is a doubling of 
old $l.$).

\end{Remark}

\begin{Theorem} \label{GenerationCondition}

Let $X$ be a separated scheme.
	\begin{enumerate}
		\item\label{GenerationCondition1} 
{\rm \cite[Proof of Lem.2.2]{1703.04484} }
If $\Dqc(X)$ is strongly compactly
generated, then 
$\Dperf(X)$ is strongly generated.
%
		\item\label{GenerationCondition2}
{\rm \cite[Proof of Lem.2.7]{1703.04484} }
 Suppose $X$ is noetherian. 
If $\Dqc(X)$ is strongly boundedly
generated, then 
$\Dcoh(X)$ is strongly generated.
%
	\end{enumerate}
\end{Theorem}

To prove these claims, the following 
 observation is crucial:

\begin{Lemma} 
\label{1.8.(i)+2.6.}
\begin{enumerate}
\item[{\rm (i)}] {\rm \cite[Prop.1.8.(i)]{1703.04484}}
Let $\mcT$ be a triangulated category with 
coproducts, and let $\mcB$ be a 
subcategory of $\mcT^c.$ Then, 
for any $l\in \mathbb{N},$
\begin{equation*} 
\mcT^c \cap 
\operatorname{Coprod}_l( \mcB ) \ 
\subseteq \ 
\smd \left( \operatorname{coprod}_l
(\mcB) \right).
\end{equation*}

\item[{\rm (ii)}]  {\rm \cite[Lem.2.6]{1703.04484}}
Let $X$ be a noetherian scheme, and let 
$G$ be an object in $\Dcoh (X).$ 
Then, 
for any $l\in \mathbb{N},$
\begin{equation*} 
\Dcoh (X) \cap 
\operatorname{Coprod}_l( 
G(-\infty,\infty)
 ) \ 
\subseteq \ 
\smd \left( \operatorname{coprod}_{2l}
(
G(-\infty,\infty)
) \right).
\end{equation*}

\end{enumerate}
\end{Lemma}

Of course, we are going to apply (i) 
with 
$$\mcT = \Dqc (X),\quad \mcB = 
G(-\infty,\infty ) \subseteq 
\mcT^c = \Dperf (X).$$ 
Then (i) becomes
$$\Dperf(X) \cap 
\operatorname{Coprod}_l( 
G(-\infty,\infty)
 ) \ 
\subseteq \ 
\smd \left( \operatorname{coprod}_{l}
(
G(-\infty,\infty)
) \right),$$
a clear analogue of (ii).

However, the point is that we can not
prove (ii) with a generality like (i).
In fact, while the proof of (i) is 
somewhat straightforward, the proof of (ii)
is more involved. 
For instance 
(see \cite[Proof of Lem.2.4]{1703.04484}), 
the \lq\lq phantom ideal\rq\rq \ $\mcI,$ 
consisting of those maps $f: x \to y$ 
such that any composite $\Sigma^iG 
\to x \xrightarrow{f} y$ vanishes for 
any $i \in \mathbb{Z}$ and any map  
$\Sigma^iG \to x$ is studied carefully, 
resorting Christensen's phantom map 
theory:

\begin{Theorem}
{\rm \cite[Th.1.1]{Chr98}}
Suppose $(\mcP,\mcI)$ is a 
\underline{\em projective class} \ 
of a triangulated category $\mcT,$ i.e.
$\mcP$ is a collection of objects in $\mcT,$
$\mcI$ is a collection of maps in $\mcT,$ 
such that

\begin{itemize}
\item $\mcP-\text{\rm null}=\mcI,$ where 
$\mcP-\text{\rm null}$ is the collection of 
\lq\lq $\mcP$-phantom maps\rq\rq 
, i.e. those
 maps $x \to y$ such that the composite 
$p\to x\to y$ is zero for all objects 
$p \in \mcP$ and all maps $p\to x.$
(This condition makes $\mcI$ an ideal.)

\item $\mcI-\text{\rm proj}=\mcP,$ where
$\mcI-\text{\rm proj}$ is the collection
of all objects $p$ such that the composite
$p \to x\to y$ is zero for all maps
$x\to y$ in $\mcI$ and all maps $p\to x.$

\item For any object $x\in \mcT,$ there 
exists a triangle
$p \to x \to y$ with $p\in \mcP$ and 
$x\to y$ in $\mcI.$

\end{itemize}

Then, for any $n\in \mathbb{N},$
$( \mcP_n, \mcI^n )$
is also a projective class, where
$\mcI^n$ is the $n$-th power of 
the \lq\lq phantom ideal\rq\rq \ $\mcI,$ 
and $\mcP_n =  \langle \mcP \rangle_n,$ 
whhich is by defined inductively 
analogous to 
Definition~\ref{<>}:
\begin{equation*}
\langle \mcP \rangle_1 = \mcP,\quad
\langle \mcP \rangle_{l+1} = 
\smd \left( \langle \mcP \rangle_1 * 
\langle \mcP \rangle_l \right).
\end{equation*}

\end{Theorem}
But, we also need some algebro-geometric 
input also to prove (ii) (see \cite[Lem.2.5]{1703.04484} \cite[Th.4.1]{LN07}).

Anyway, assuming Lemma~\ref{1.8.(i)+2.6.},
the proof of Theorem~\ref{GenerationCondition} becomes
straightforward:

\begin{proof}[Proof of 
Theorem~\ref{GenerationCondition} 
assuming  Lemma~\ref{1.8.(i)+2.6.}:
]
In both cases, assuming the respective 
assumption on $\Dqc (X),$ together with 
Remark~\ref{remark on strongly compactly or boundedly generated}, the claims follow 
as follows:
\begin{equation*}
\begin{split}
\Dperf (X) &= \Dperf (X) \cap \Dqc (X) =
\Dperf(X) \cap 
\operatorname{Coprod}_l( 
G(-\infty,\infty)
 ) 
\\ 
&
\subseteq \ 
\smd \left( \operatorname{coprod}_{l}
(
G(-\infty,\infty)
) \right)
\subseteq \cup_n \langle G\rangle_l^{[-n,n]}.
\\
\Dcoh (X) &= \Dcoh (X) \cap \Dqc (X) =
\Dcoh(X) \cap 
\operatorname{Coprod}_l( 
G(-\infty,\infty)
 ) 
\\ 
&
\subseteq \ 
\smd \left( \operatorname{coprod}_{2l}
(
G(-\infty,\infty)
) \right)
\subseteq \cup_n \langle G\rangle_{2l}^{[-n,n]}.
\end{split}
\end{equation*}

\end{proof}

}

\subsection{Strong generation of $\Dperf (X)$}

From Theorem~\ref{GenerationCondition} 1, 
we search for situations when 
$\Dqc(X)$ becomes strongly compactly
generated:

\begin{Theorem}[Max Kelly 
\cite{Kel65}]
\label{Max Kelly}
	Suppose $X=\Spec R$ is affine. Then $\Dqc(X)$ is strongly compactly generated if and only if $R$ is of finite global dimension.
\end{Theorem}



\begin{Theorem}[Bondal--Van den Bergh \cite{MR1996800}]\label{BondalVandenBergh}
	Let $X$ be smooth scheme of finite type over a field $k$. Then $\Dqc(X)$ is strongly compactly generated.
\end{Theorem}


Theorem~\ref{BondalVandenBergh} has recently been improved by Orlov as a characterization
of the strong generation of $\Dperf (X)$:

\begin{Theorem} {\rm (Orlov 
\cite[Th.3,27]{MR3545926}
)}
Let $X$ be a separated noetherian scheme of finite Krull dimension over an arbitrary
fieled $k.$ Assume that the square $X\times X$ 
is noetherian too. 
Then the following conditions are
equivalent:
\begin{enumerate}
\item $X$ is regular;
\item $\Dperf (X)$ is strongly generated.
\end{enumerate}

\end{Theorem}

It was
this paper of Orlov \cite{MR3545926} which 
motivated Neeman to develop his theory of
approximable triangulated category 
(see e.g. \cite[p.6, the paragraph before Rem.0.10]{1703.04484}).

In fact, the approximability of $\Dqc (X)$
allowed Neeman to prove the following statement
by reducing to the Kelly's old theorem 
in a straightforward way, 
i.e. by induction on the
number of open affines covering $X$:

\begin{Theorem}{\rm (Neeman 
\cite[Th.2.1]{1703.04484})}
\label{SCG criterion}
	Let $X$ be a quasi-compact separated scheme. 
If $X$ can be covered by open affines $\Spec R_{i}$ with $R_{i}$ of finite global dimension, then $\Dqc(X)$ is strongly compactly generated.
\end{Theorem}

\begin{proof}{(Outline of a proof of 
Theorem~\ref{SCG criterion} following 
\cite[Sketch.6.6]{1806.06995})}
\ Proceed as follows:

\begin{itemize}
\item Write $X = \cup_{1\leq i\leq r}U_i$ 
with $u_i : U_i = \operatorname{Spec}(R_i),$ by assumption.

\item By induction on $r$ unsing the
Mayer–Vietoris sequence 
\cite[Prop.5.10]{MR2434186}
(as in the proof 
given in 
\cite[Proof of Theorem 2.1]{1703.04484}), 
we find
\begin{equation}
\hspace{-4mm}
\begin{split}
\label{MayerVietoris}
&\quad \Dqc (X) = 
\\
&\underbrace{
\bigg( \add \left[ \cup_{i=1}^r 
\Rd {u_i}_* \Dqc (U_i) \right] \bigg)  * 
\bigg( \add \left[ \cup_{i=1}^r 
\Rd {u_i}_* \Dqc (U_i) \right] \bigg)  * 
\cdots *
\bigg( \add \left[ \cup_{i=1}^r 
\Rd {u_i}_* \Dqc (U_i) \right] \bigg) }
_{\text{$r$ copies}}.
\end{split}
\end{equation}

\item By a minor variant of Max Kelly's 
Theorem~\ref{Max Kelly}, 
\begin{equation}
\label{by Max Kelly}
\exists l\in \mathbb{N},\ \text{s.t.}\quad 
1\leq \forall i\leq r,\quad 
\Dqc (U_i) = 
\overline{\langle \mcO_{U_i}\rangle}
^{(-\infty,\infty)}_l.
\end{equation}

\item From Corollary~\ref{18d,Lem.6.5} 
(recall
\underline{\em
 this is where the 
approximability of $\Dqc (X)$ was 
exploited}), 
\begin{equation}
\label{by approximability}
\exists n\in \mathbb{N}\ \text{s.t.}\quad 
1\leq \forall i\leq r,\quad 
\Rd {u_i}_* \mcO_{U_i} \in 
\overline{\langle G\rangle}^{[-n,n]}_n 
\subset \Dqc (X).
\end{equation}

\item From \eqref{by Max Kelly} and 
\eqref{by approximability}, 
\begin{equation*}
\Rd {u_i}_* \Dqc (U_i) = \Rd {u_i}_*
\left[
\overline{\langle \mcO_{U_i}\rangle}
^{(-\infty,\infty)}_l
\right]
\subset 
\overline{\langle 
\Rd {u_i}_*
\mcO_{U_i}\rangle}
^{(-\infty,\infty)}_l
\subset 
\overline{\langle G\rangle}
^{[-\infty,\infty]}_{ln},
\end{equation*}
and so
\begin{equation} \label{last key}
\add \left[ \cup_{i=1}^r 
\Rd {u_i}_* \Dqc (U_i) \right] \ \subset \
\overline{\langle G\rangle}
^{[-\infty,\infty]}_{ln},
\end{equation}

\item From \eqref{MayerVietoris} and \eqref{last key},
we obtain the desired strong compact 
generation of $\Dqc (X)$:
\begin{equation*}
\Dqc (X) = \overline{\langle G\rangle}
^{[-\infty,\infty]}_{lnr},
\end{equation*}

\end{itemize}


\end{proof}

Now, Neeman proves his main theorem on 
strong generation of $\Dperf (X)$:


\begin{Theorem}{\rm (Neeman 
\cite[Th.0.5]{1703.04484}
\cite[Th.6.1]{1806.06995})}
\label{Th.6.1-1806.06995}
	Let $X$ be a quasi-compact separated scheme. Then $\Dperf(X)$ is strongly generated if and only if $X$ can be covered by open affines $\Spec R_{i}$ with $R_{i}$ of finite global dimension.
\end{Theorem}

\begin{proof}

\begin{description}

\item[\underline{\lq\lq if\rq\rq \ part}]
This is immediate 
from Theorem~\ref{SCG criterion} and 
Theorem~\ref{GenerationCondition}(1).

\item[\underline{\lq\lq only if\rq\rq \  part}] \cite[Rem.0.10]{1703.04484}
By Thomason-Trobaugh \cite{MR1106918} 
recalled in Theorem~\ref{Thomason's localization theorem} and \eqref{IC for Perf - idempotent completion}, we have an 
equivalence upon idempotent completion:
 \begin{equation*}
\left(
\Dperf (X) \big/ \left( {\Dperf} \right)_Z (X) 
\right)^{\sharp}
\
\xrightarrow[\cong]{\Ld j^*}
\
\Dperf (U).
\end{equation*}
Thus, if $G\in \Dperf(X)$ is a strong 
generator, then so is 
$\Ld j^* G \in \Dperf (U).$
Now the strong generation of an affine 
$U = \operatorname{Spec} (R)$ forces 
$R$ to be of finite global dimension, 
as is shown in \cite[Prop.7.25]{MR2434186}.

\end{description}
\end{proof}

\subsection{Strong generation of $\Dcoh (X)$}


Here, we start with a nice theorem of Rouquier:

\begin{Theorem} {\rm (Rouquier 
\cite[Th.7.39]{MR2434186}) }
\label{Rouquiter}
	Let $X$ be a scheme of finite type over a perfect field $k$. Then $\Dqc(X)$ is strongly boundedly generated, and $\Dcoh (X)$ is 
strongly generated.
\end{Theorem}

To go further, let us recall:
\begin{itemize}
\item
the canonical map
$\Dperf (X) \to \Dcoh (X)$ is an isomorphism when 
$X$ is smooth over a field, and in this case,
the strong generation of $\Dcoh (X) \cong \Dperf (X)$ 
is already discussed in the previous subsection.
\item
the Verdier quotient
$\DD_{sg}(X) = \Dcoh (X) / \Dperf (X)$ reflects 
singular information of $X.$  
\end{itemize}

Thus, we must take care of singular property of $X.$  
However, while Theorem~\ref{Th.6.1-1806.06995}
is easy and classical in the case where $X$ is affine, this problem is 
{\it neither easy nor classical for affine $X.$}
See \cite[H.S..6.12]{1806.06995} for more on this 
point.
\footnote{
In fact, when $X$ is affine, 
strong generation of $\Dqc(X)$ has been proved by
Iyengar and Takahashi \cite{IT16} 
under different hypotheses, and using quite
different techniques, from Neeman's 
Theorem~\ref{SBG criterion}.
And they give examples where strong generation
fails; see \cite{IT16} and references therein.
}

Now, for this purpose, Neeman turned his attention to 
 de Jong's alteration:
\footnote{
(Gabber's strengthening \cite{Gab05} of)
de Jong's alteration is now widely used 
in the Morel-Voevodsky motivic stable
homotopy theory. see e.g. \cite{Kel13}
\cite{HKO17}.
For an introductory review of de Jong's 
alteration, consult Oort's \cite{Oor98} 
for instance.
}

\begin{Definition} {\rm 
\cite{deJ96}\cite{deJ97}
\cite{Oor98}
\cite[Remi.0.13]{1703.04484}}
Let $X$ be a noetherian scheme. 
A \emph{regular alteration} of $X$ is a proper,
surjective morphism $f : Y \to X,$ so that
\begin{enumerate}
\item 
$Y$ is regular and finite dimensional.
\item 
There is a dense open set $U \subset X$ over which 
$f$ is finite.
\end{enumerate}

\end{Definition}

Now, Neeman proves:

\begin{Theorem}
{\rm (Neeman 
\cite[Th.2.3]{1703.04484})}
\label{SBG criterion}

	Let $X$ be a noetherian scheme, and assume every closed subscheme $Z\subset X$ admits a regular alteration. Then $\Dqc(X)$ is strongly boundedly generated.
\end{Theorem}
{

\begin{proof}{(Outline of a proof of 
Theorem~\ref{SBG criterion} following 
\cite[Proof that Theorem 2.3 follows from Theorem 2.1 ]{1703.04484}
)}
\footnote{
This proof does not directly use the of 
 approximability of $\Dqc(X),$ the approximability
enters only indirectly, when we appeal to Theorem 5.10. What we want to highlight here, following a strong suggestion of 
Professor Neeman, 
is the pivotal role that the 
homotopy-theoretical ideas of Bousfield, 
Ohkawa, Hopkins-Smith and many others play in the reduction."

}:
This is proved in the following order:
%
\ 
\begin{itemize}
\item 
Suppose there is a counterexample $X$ to 
Theorem~\ref{SBG criterion}(SBG criterion)
.  Since $X$ is 
noetherian, we may choose a minimal closed 
subscheme $Z \subset X$ which does not 
satisfy Theorem~\ref{SBG criterion}(SBG criterion).

\item
Replacing $X$ by $Z,$ may assume all 
proper closed subschemes $Z \subset X$ 
satisfy Theorem~\ref{SBG criterion}(SBG criterion).

\item 
To prove Theorem~\ref{SBG criterion}(SBG criterion) for $X,$ we may assume it is 
\underline{\em reduced}: for, let 
$j : X_{\text{\rm red}} \to X$ be the 
inclusion of the reduced part of $X,$ 
and let $\mcJ$ be the corresponding 
ideal sheaf with $\mcJ^n=0.$
Then, expressing any $C\in \Dqc (X)$ by 
a complex of quasi-coherent sheaves, we 
obtain a filtration
\begin{equation*}
0 = \mcJ^nC \subset  \mcJ^{n-1}C \subset
\cdots \subset \mcJ C  \subset C,
\end{equation*}
with $\mcJ^jC/\mcJ^{j+1} \in 
\Rd j_* \Dqc (X_{\text{\rm red}}) \ 
(0\leq \forall j\leq n-1 ).$ Then, as in 
\cite[7.3]{MR2434186}, we find: 
\begin{equation*}
C \ \in \ 
\left[ \Rd j_* \Dqc (X_{\text{\rm red}})
\right]^{* n} =
\underbrace{
\left[ \Rd j_* \Dqc (X_{\text{\rm red}})
\right] *
\left[ \Rd j_* \Dqc (X_{\text{\rm red}})
\right] * \cdots *
\left[ \Rd j_* \Dqc (X_{\text{\rm red}})
\right]
}_{\text{$n$}}.
\end{equation*}
So, it suffices to prove the strong 
bounded gneration 
$$\Dqc (X_{\text{red}} )
= \operatorname{Coprod}_{\widetilde{N}}
\left( \widetilde{G} ( - \infty, \infty ) 
\right)$$ 
for some $\widetilde{N} \in 
\mathbb{N}$ and some $\widetilde{G} \in 
\Dcoh (X_{\text{red}} ),$ for then 
we would get:
\begin{equation*}
\begin{split}
\Dcoh (X) &\subseteq 
\left[ \Rd j_* \Dqc (X_{\text{\rm red}})
\right]^{* n} = 
\big[
\Rd j_* 
\operatorname{Coprod}_{\widetilde{N}}
\left( \widetilde{G} ( - \infty, \infty ) 
\right)
\big]^{* n}
\\
&\subseteq 
\big[
\operatorname{Coprod}_{\widetilde{N}}
\left( 
( \Rd j_* \widetilde{G} )
 ( - \infty, \infty ) 
\right)
\big]^{* n}
=
\operatorname{Coprod}_{\widetilde{N} n}
\left( 
( \Rd j_* \widetilde{G} )
 ( - \infty, \infty ) 
\right),
\end{split}
\end{equation*}
where $\Rd j_* \widetilde{G} \in 
\Dcoh (X )$ by Theorem~\ref{proper morphism preserves coherence}.
So, the strong bounded generation of 
$\Dcoh (X)$ would follow.

\item Now that we may assume $X$ is reduced, 
we may apply de Jong's regular alteration to $X$:
\begin{equation*}
\xymatrix{
Y \ar[rr]^f_-{\text{proper \& surjective}}
 & &  X  
\\
f^{-1}(U) \ar@{^{(}-_{>}}[u]
 \ar[rr]^{f|_{f^{-1}(U)}}_-{\text{finite \& flat}} 
& &
\exists U 
\ar@{^{(}-_{>}}[u]_{\text{dense open}}
}
\end{equation*}
where we may apply  
Theorem~\ref{SCG criterion}(SCG criterion)
to $Y,$ because $Y$ is finite-dimensional,
separated and regular: 
Here, let us consider 
$\Rd f_* \left( \mcO_Y \oplus 
\Sigma \mcO_Y \right) \in \Dcoh (X)$ 
(see Theorem~\ref{proper morphism preserves coherence}).  Then,
\begin{itemize}
\item 
Since $f|_{f^{-1}(U)}$ is 
finite, flat and surjective, the 
 restriction to $U$ of the object
$\Rd f_*  \mcO_Y \in \Dqc (X)$ is a 
nowhere vanishing vector bundle on $U.$ 
In particular,
%
%
\begin{equation} 
\label{restriction to U is perfect}
\left( \Ld j^* \Rd f_* \mcO_Y \right)
\oplus \Sigma 
\left( \Ld j^* \Rd f_* \mcO_Y \right)
=
\Ld j^* \Rd f_*  \left( 
\mcO_Y \oplus 
\Sigma  \mcO_Y
\right) \ 
\in \
\Dperf (U).
\end{equation}

\item Then, we can 
\newline
\underline{\em
apply 
Corollary~\ref{A useful consequence of density}, a corollary of 
Thomason's localization theorem 
}
(Theorem~\ref{Thomason's localization theorem}), to \eqref{restriction to U is perfect} to
find some $H \in \Dperf (X)$ such that
\begin{equation} 
\label{local isomorphism}
\Ld j^* H  \xrightarrow{\cong} \Ld j^*  
\Rd f_* \left( \mcO_Y \oplus 
\Sigma  \mcO_Y \right) \ 
\in \
\Dperf (U).
\end{equation}

\end{itemize}

\item To the local isomorphism 
\eqref{local isomorphism}, applying the ajoint isomorphism
$$
\Hom_{\Dqc (U)}\left(  
\Ld j^* H,\ \Ld j^*  
\Rd f_* \left( \mcO_Y \oplus 
\Sigma  \mcO_Y \right) \right) \cong 
\Hom_{\Dqc (X)}\left(  
H,\  \Rd j_* \Ld j^*  
\Rd f_* \left( \mcO_Y \oplus 
\Sigma  \mcO_Y \right) \right),
$$ 
we obtain a map 
\footnote{
\underline{WARNING!} \ 
In \cite[Proof that Theorem 2.3 follows from Theorem 2.4]{1703.04484}, Neeman 
concluded the existence of an honest map 
$H \to \Rd f_* 
\left( \mcO_Y \oplus \Sigma \mcO_Y\right)$ 
corresponding to \eqref{local isomorphism}. 
However, this is quite problematic, and 
usually, such an honest map 
$H \to \Rd f_* \mcO_Y \oplus \Sigma 
\Rd f_* \mcO_Y$ does not exist.
Thus, some sort of patch is needed.
The \lq\lq patch\rq\rq \ presented above 
 was  communicated to the author by 
Professor Neeman, and the author replaced
his own patch, which concentrates on 
$\widetilde{R}$ (see \eqref{The homotopy pullback}), 
with Professor Neeman's 
\lq\lq patch\rq\rq , which concentrates 
on $\widetilde{H}$ 
(see \eqref{The homotopy pullback}), 
because Professor 
Neeman's patch delivers a simple 
 message how to read 
 \cite[Proof that Theorem 2.3 follows from Theorem 2.4]{1703.04484}:
just replace $H$ with $\widetilde{H}$ and 
pretend the map $\widetilde{\psi}' : 
\widetilde{H} \to \Rd f_* 
\left( \mcO_Y \oplus 
\Sigma  \mcO_Y \right)$
obtained in \eqref{The homotopy pullback}
as our \lq\lq honest map\rq\rq \
$H \to \Rd f_* 
\left( \mcO_Y \oplus 
\Sigma  \mcO_Y \right),$ and then, 
just proceed as is written in 
\cite[Proof that Theorem 2.3 follows from Theorem 2.4]{1703.04484}.

According to Professor Neeman, this 
leap and omission of justification 
is standard. So, the reader
is required to come up with this kind of
patch spelled out in terms of elementary
Bousfield (or Miller's finite) localization
instantaneously 
at the top of his or her head.
Thus, homotopy theoretical insight is 
prerequisite to read Professor Neeman's 
papers!
}
\begin{equation}
\label{psi map}
\psi :  H \to  \Rd j_* \Ld j^*  
\Rd f_* \left( \mcO_Y \oplus 
\Sigma  \mcO_Y \right).
\end{equation}

\begin{itemize}
\item  Recall, 
since   
$\left( {\Dqc} \right)_Z (X)$ is 
compactly generated (\cite[Th.6.8]{MR2434186}), we can 
apply Miller's finite localiztion 
Theorem~\ref{Nee92b, Lem.1.7} to 
form the Verdier quotient 
with the equivalence \eqref{BL of TC}:
\begin{equation}
\label{BL of TC - RE}
\Dqc (X) \big/ \left( {\Dqc} \right)_Z (X) 
\
\xrightarrow[\cong]{\overline{\Ld j^*}}
\
\Dqc (U),
\end{equation}
and that $\Rd j_*\Ld j^*$ which shows 
up in the target of the 
$\psi$ map \eqref{psi map}
can be interpreted as the 
Bousfield localization, as in 
\eqref{smashing Bousfield localization of Dqc(X)}, which is consequently 
expressed by a  mapping telescope 
$\operatorname{\bf hocolim}$ 
as Miller's finite localization 
(Theorem~\ref{Nee92b, Lem.1.7}).
Then, cosider the following pair of maps:
\begin{equation} \label{pair of maps}
H \xrightarrow{\psi}   
\Rd j_* \Ld j^*  
\Rd f_* \left( \mcO_Y \oplus 
\Sigma  \mcO_Y \right) = 
\operatorname{hocolim} (R_n) 
\xleftarrow[\text{canonical map}]{c} R_0 =
\Rd f_* \left( \mcO_Y \oplus 
\Sigma  \mcO_Y \right).
\end{equation}

\item The both maps in \eqref{pair of maps} 
are local isomorphism, i.e. isomorphisms 
when restricted $U.$ This is trivial for 
the canonical map (which is the Bousfield 
localization) and the claim for $\psi$ 
follows from the local isomorphism 
\eqref{local isomorphism}.

\item  Since $H \in \Dperf (X) = 
\Dqc (X)^c$ is compact, arguing as in 
Proposition~\ref{compact objects observation} and its comments below, 
we may factorize the pair of maps 
\eqref{pair of maps} as follows:

\begin{equation}
\label{factorize the pair of maps}
\xymatrix{
H \ar[dr]_-{\exists \widetilde{\psi}} 
\ar[r]^-{\psi}  &   
\operatorname{hocolim}(R_n)  & 
\Rd f_* \left( \mcO_Y \oplus 
\Sigma  \mcO_Y \right)
\ar[l]_-{
c}
\ar[dl]^-{\exists \widetilde{c}}
\\
& \exists \widetilde{R} \ar[u]^{\iota} & 
},
\end{equation}
where:
\begin{itemize}
\item $\widetilde{R}$ is obtained from 
$\Rd f_* \left( \mcO_Y \oplus 
\Sigma  \mcO_Y \right) \ \in \
\Dcoh (X)$ 
via $\widetilde{c}$ 
by a finite step 
extensions of finite coproducts of 
elements in  $\Dperf (X).$
Thus, we have a triangle of the 
following form:
\begin{equation}
\label{triangle 1}
\Rd f_* \left( \mcO_Y \oplus 
\Sigma  \mcO_Y \right) 
\xrightarrow{\widetilde{c}} 
\widetilde{R} \to Q'\quad 
(Q' \in (\Dperf)_Z (X), 
\widetilde{R} \in \Dcoh (X) )
\end{equation}

\item From \eqref{triangle 1}, we see 
$\widetilde{c}$ is a local isomorphism, 
then, since $c$ is also a local 
isomorphism, $\iota$ is a local isomorphism 
as well from the right hand side
commutative diagram of 
\eqref{factorize the pair of maps}.

Then, since $\phi$ is also a local 
isoorphism, from the left hand side 
commutative diagram of 
\eqref{factorize the pair of maps}, 
we find $\widetilde{\psi}$ is also 
a local isomorphism.
Thus, we have a triangle of the 
following form:
\begin{equation}
\label{triangle 2}
Q'' \to H 
\xrightarrow{\widetilde{\psi }} 
\widetilde{R} \quad 
(Q'' \in (\Dcoh)_Z (X) )
\end{equation}

\end{itemize}

\item Take the homotpy pullback 
$\widetilde{H}$ of the pair of maps
$H \xrightarrow{\widetilde{\psi}}  
\widetilde{R} \xleftarrow{\widetilde{c}} 
\Rd f_* \left( \mcO_Y \oplus 
\Sigma  \mcO_Y \right)$ obtained in 
\eqref{factorize the pair of maps}:

\begin{equation}
\label{The homotopy pullback}
\xymatrix{
& \widetilde{H} := 
H\times^h_{\widetilde{R}} 
\Rd f_* \left( \mcO_Y \oplus 
\Sigma  \mcO_Y \right)
\ar[dl]_-{\widetilde{c}'}
\ar[dr]^-{\widetilde{\psi}'}
&
\\
H \ar[dr]_-{\widetilde{\psi}} 
&   &  
\Rd f_* \left( \mcO_Y \oplus 
\Sigma  \mcO_Y \right)
\ar[dl]^-{\widetilde{c}}
\\
& \widetilde{R} &
}
\end{equation}
where:

\begin{itemize}
\item From \eqref{triangle 1},
the homotopy pullback diagram
\eqref{The homotopy pullback} 
and 
$H \in \Dperf (X),$ we 
have a triangle of the following form:

\begin{equation}
\label{triangle 1'}
\widetilde{H} 
\xrightarrow{\widetilde{c}'} H \to Q'
\quad 
(Q' \in (\Dperf)_Z (X),\ 
H, \widetilde{H} \in \Dperf(X)
 )
\end{equation}

\item From \eqref{triangle 2} and
the homotopy pullback diagram
\eqref{The homotopy pullback}, 
we 
have a triangle of the following form:

\begin{equation}
\label{triangle 2'}
Q'' \to \widetilde{H} 
\xrightarrow{\widetilde{\psi}'} 
\Rd f_* \left( \mcO_Y \oplus 
\Sigma  \mcO_Y \right)
\quad 
(Q'' \in (\Dcoh)_Z (X) )
\end{equation}

\end{itemize}

\end{itemize}

\item Concerning the homological support 
$\Supph (\widetilde{H} )$ 
of $\widetilde{H} 
\overset{\eqref{triangle 1'}}\in 
\Dperf (X)$, we see:
\begin{itemize}
\item $\Supph (\widetilde{H} )$ is closed,  
because $\widetilde{H} \in \Dperf (X)$ 
implies $\mcH^{\bullet}\widetilde{H}$ 
is of finite type as an $\mcO_X$-module, 
and so we may apply \cite[Lem.17.9.6]{Stack} for instance.

\item \ \vspace{-5mm}
\begin{equation*}
\begin{split}
\Supph \left( \widetilde{H} \right)
\bigcap U 
&\overset{\eqref{triangle 1'}}= 
\Supph \left( H \right)
\bigcap U
\overset{\eqref{local isomorphism}}= 
\Supph \left(
\Rd f_* \left( \mcO_Y \oplus 
\Sigma  \mcO_Y \right)
\right)
\bigcap U
\\
&\overset{\text{direct summand}}{\supseteqq}
\Supph \left(
\Rd f_*  \mcO_Y 
\right)
\bigcap U = U,\ \text{a dense open of $X.$}
\end{split}
\end{equation*}
where the last equality follows from 
the fact $\Rd f_*  \mcO_Y$ restricted 
to $U$ is a nowhere vanishing vector
bundle.

\end{itemize}

Thus the homological support 
$\Supph \left( \widetilde{H} \right)$ 
is whole 
$X.$
%
Then, we can 
\newline
\underline{\em 
apply
Corollary~\ref{thick tensor ideal with suppph the whole space}, a corollary 
of Thomason's theorem of Thomason sets} 
\newline
(Theorem~\ref{Thomason's theorem of Thomason sets}) to conclude that, 
$\langle \widetilde{H}\rangle_{\otimes},$ the
 tensor ideal generated by $\widetilde{H},$ is the whole
$\Dperf (X),$ which obviously contains 
$\mcO_X.$ Then, applying 
Remark~\ref{thick tensor ideal generated by a single element} and 
Proposition~\ref{Nee17,Cor.1.11}, we may pick some 
$C\in \Dperf(X)$ and $L\in \mathbb{N}$ such 
that
\begin{equation}
\mcO_X \in \langle C\otimes \widetilde{H}\rangle_L
\subseteq 
\operatorname{Coprod}_{2L}
\left( \left(C\otimes \widetilde{H} \right)(-\infty,\infty) \right).
\end{equation}
Consequently, for any 
$D \in \Dqc (X),$ 
\begin{equation}
\label{important relation for D}
D = D\otimes
 \mcO_X 
\ \in \ \langle D\otimes C\otimes \widetilde{H}\rangle_L
\ \subseteq \ 
\operatorname{Coprod}_{2L}
\left( \left(D\otimes C\otimes \widetilde{H}\right)(-\infty,\infty) \right).
\end{equation}

\item 
Having 
\eqref{important relation for D}
in mind, we apply 
$D\otimes C\otimes -$ to 
\eqref{triangle 2'}
to
obtain the following triangles:
\begin{equation}
\label{tensored triange}
D\otimes C\otimes Q'' \to
D\otimes C\otimes \widetilde{H} \to
D\otimes C\otimes 
\Rd f_* \left( \mcO_Y \oplus 
\Sigma  \mcO_Y \right)
\end{equation}
where 
$\Rd f_* \left( \mcO_Y \oplus 
\Sigma  \mcO_Y \right) \ \in \
\Dcoh (X),\ 
Q''\in ( \Dcoh )_Z(X).$

\item For $Y,$ obtained by de Jong's 
regular alteration, we may apply 
Theorem~\ref{SCG criterion} to conclude 
its strong compact generation.
Thus, $\exists G\in \Dperf (X),\ 
\exists N\in \mathbb{N},$ s.t. 
$
\Dqc (Y ) = \operatorname{Coprod}_N
\left( G(-\infty,\infty) \right).
$ Hence,
\begin{equation*}
\Ld f^*
( D\otimes C )
\otimes
\left(
\mcO_Y \oplus 
\Sigma  \mcO_Y 
\right)
 \ \in \ \Dqc (Y ) =
\operatorname{Coprod}_N
\left( G(-\infty,\infty) \right) \qquad 
( G \in \Dperf (X) )
\end{equation*}
Consequently, by the projection formula, 
\begin{equation}
\label{Rf*V is generated by Rf*G}
\begin{split}
D\otimes C\otimes
\Rd f_*
\left( \mcO_Y \oplus \Sigma  \mcO_Y \right)
=
\Rd f_*
\bigg(
\Ld f^*
( D\otimes C )
\otimes
\left(
\mcO_Y \oplus 
\Sigma  \mcO_Y 
\right)
\bigg)
\\
 \in \ \Rd f_*
\operatorname{Coprod}_N
\left( G(-\infty,\infty) \right)
\ \subseteq \ 
\operatorname{Coprod}_N
\left( (\Rd f_*G)(-\infty,\infty) \right)
\end{split}
\end{equation}
where $\Rd f_*G \in \Dcoh (X)$ 
by Theorem~\ref{proper morphism preserves coherence}.

\item For $
Q'' \in (\Dcoh)_Z(X)$ in 
\eqref{tensored triange}, 
we may apply Rouquier's 
Theorem~\ref{Rouquier's Lemmma 7.1}
to find 
$n\in \mathbb{N}, P_{n}
\in \Dcoh (Z_{n})$
%
s.t. 
\begin{equation} 
\label{Q''}
Q'' = \Rd {i_{n}}_*P_{n}
\qquad ( P_{n}
\in \Dcoh (Z_{n}) ). 
\end{equation}

\item For 
$Z_n,$ whose 
underlying space is equal to that of 
the proper closed subscheme
$Z$ of $X$ from their constructions in
Theorem~\ref{Rouquier's Lemmma 7.1}, 
we may apply Theorem~\ref{SBG criterion} 
by inductive assumption 
to conclude their strong bounded 
generations. Thus, 
 $\exists G'' \in \Dcoh (Z_{n} ), 
\exists M \in \mathbb{N}$ s.t. 
%
$\Dqc (Z_{n} ) = 
\operatorname{Coprod}_{M}
(G''(-\infty,\infty)).$
%
Hence, 
\begin{equation}
\label{P}
\Ld i_{n}^*(D\otimes C)\otimes
P_{n} \ \in \ 
\Dqc (Z_{n} ) = 
\operatorname{Coprod}_{M}
(G''(-\infty,\infty)) \qquad
( G'' \in \Dcoh (Z_{n} ) )
\end{equation}
%
Consequently, by the projection formula,
\begin{equation}
\label{Ri*P is generated by Ri*G}
\begin{split}
&D\otimes C\otimes
Q'' = 
D\otimes C\otimes 
\Rd {i_{n}}_* P_{n}
= 
\Rd {i_{n}}_* \bigg( 
\Ld i_{n}^*(D\otimes C)\otimes
P_{n} 
\bigg)
\\
\in \ 
&\Rd {i_{n}}_*
\operatorname{Coprod}_{M}
(G''(-\infty,\infty)) 
\ \subseteq \ 
\operatorname{Coprod}_{M}
( (\Rd {i_{n}}_*G'')(-\infty,\infty)) 
%
\end{split}
\end{equation}
%
%
%
%
where 
$\Rd {i_{n}}_*G'' \in \Dcoh (X)$ 
by Theorem~\ref{proper morphism preserves coherence}.

\item 
\ 
From 
\eqref{triangle 2'}
\eqref{Rf*V is generated by Rf*G}
\eqref{Ri*P is generated by Ri*G}, we find
\footnote{
In Neeman's corresponding calculation  
\cite[1st paragraph in p.24]{1703.04484},
the extension length of 
$\operatorname{Coprod}$ was doubled 
to be $2(M+N)$ rather than 
$M+N$ given  in  \eqref{widetildeH is generated by}.
However, the author does not see such a 
need, and so, the author opted to present 
as in \eqref{widetildeH is generated by}. 
%
}
\begin{equation}
\label{widetildeH is generated by}
\begin{split}
&\quad \  D\otimes C \otimes 
\widetilde{H} \ 
\in \ 
\operatorname{Coprod}_{M}
( (\Rd {i_{n}}_*G'')(-\infty,\infty)) 
*
\operatorname{Coprod}_N
\left( (\Rd f_*G)(-\infty,\infty) \right) 
\\
&\subseteq \ 
\operatorname{Coprod}_{M}
( (\Rd f_*G \oplus \Rd {i_{n}}_*G'')(-\infty,\infty)) 
*
\operatorname{Coprod}_N
\left( (\Rd f_*G
\oplus \Rd {i_{n}}_*G'')
(-\infty,\infty) \right)
\\
&\subseteq \ 
\operatorname{Coprod}_{M+N}
\left( (\Rd f_*G
\oplus \Rd {i_{n}}_*G'')
(-\infty,\infty) \right),
\end{split}
\end{equation}
where $\Rd f_*G \oplus 
\Rd {i_{n}}_*G'' \in \Dcoh (X).$
%
%
%

\item Finally, from 
\eqref{important relation for D}
\eqref{widetildeH is generated by}
we see 
for any $D\in \Dqc (X),$
\begin{equation}
\label{final relation}
\begin{split}
D &\overset{\eqref{important relation for D}}
\in \operatorname{Coprod}_{2L}
\left( (D\otimes C\otimes \widetilde{H} )
(-\infty,\infty) \right) 
\\
&\overset{\eqref{widetildeH is generated by}}\subseteq
\operatorname{Coprod}_{2L}
\bigg(
\big(
\operatorname{Coprod}_{M+N}
\left( (\Rd f_*G
\oplus \Rd {i_{n}}_*G'')
(-\infty,\infty) \right)
\big)
(-\infty,\infty)  
\bigg)
\\
&\subseteq 
\operatorname{Coprod}_{2L(M+N)}
\left( (\Rd f_*G
\oplus \Rd {i_{n}}_*G'')
(-\infty,\infty) \right),
%
%
\end{split}
\end{equation}
where $
\Rd f_*G
\oplus \Rd {i_{n}}_*G''
\in 
\Dcoh (X).$
Thus, we have obtained the desired
\begin{equation*}
\Dqc (X) = 
\operatorname{Coprod}_{
2L(M+N)}
\left( 
\Rd f_*G \oplus \Rd {i_{n}}_*G'' )
(-\infty,\infty) \right),
\end{equation*}
%
%
%
which shows the strong bounded 
generation of $\Dqc (X)$ for 
$\Rd f_*G
\oplus \Rd {i_{n}}_*G''
\in 
\Dcoh (X).$


\end{itemize}
%
%

\end{proof}

From Theorem~\ref{SBG criterion}
and Theorem~\ref{GenerationCondition}
(2), we obtain Neeman's main theorem 
on strong generation of $\Dcoh (X)$:

\begin{Theorem}{\rm (Neeman 
\cite[Th.0.15]{1703.04484}
\cite[Th.6.11]{1806.06995})}
\label{Th.0.15}
Let $X$ be a noetherian scheme, and assume every closed subscheme $Z\subset X$ admits a regular alteration. Then $\Dcoh(X)$ is strongly generated.

\end{Theorem}

From 
\cite{deJ96}\cite{deJ97} and 
\cite{Nay09}, we see any $X,$ which is 
separated and essentially of finite
type over a separated excellent scheme $S$ 
of dimension $\leq 2,$ satisfies 
the assumptions of 
Theorem~\ref{SBG criterion} and 
Theorem~\ref{Th.0.15}.
Thus, Theorem~\ref{SBG criterion}
and Theorem~\ref{Th.0.15} 
generalize 
Rouquier's Theorem~\ref{Rouquiter}.

}

For more details about strong generations of
$\Dperf (X)$ and $\Dcoh (X),$ consult Neeman's 
original article \cite{1703.04484} and the survey
\cite{1806.06995}.

\end{document}